\RequirePackage{fix-cm}
\documentclass[12pt]{amsart}
\usepackage{amsmath,amssymb}
\usepackage{graphicx}
\usepackage{mathptmx}      
\usepackage{latexsym}
\usepackage{xcolor}

\newtheorem{proposition}{Proposition}
\newtheorem{theorem}{Theorem}
\newtheorem{definition}{Definition}
\newtheorem{corollary}{Corollary}
\newtheorem{lemma}{Lemma}
\newtheorem{remark}{Remark}

\newcommand{\Iup}{\overline{I}}
\newcommand{\He}{\mathcal{H}}

\newcommand{\Me}{\mathcal{M}}

\newcommand{\Se}{\mathcal{S}}
\newcommand{\Ee}{\mathcal{E}}
\newcommand{\Fe}{\mathcal{F}}

\begin{document}

\title{Conditional Non-Lattice Integration, Pricing and Superhedging}



\author{C. Bender}
\address[C. Bender]{Department of Mathematics\\
Saarland University\\
Campus E 2 4,
66123 Saarbr\"ucken
Germany}
\email{bender@math.uni-sb.de}
\author{S.E.Ferrando}
\address[S.E.Ferrando]{Department of Mathematics\\
Ryerson University\\
350 Victoria Street}
\email{ferrando@ryerson.ca}
\author{A.L. Gonzalez}
\address[A.L. Gonzalez]{Department of Mathematics\\
Mar del Plata National University\\
Mar del Plata, Funes 3350, 7600 Argentina}
\email{algonzal@mdp.edu.ar}


\date{Received: date / Accepted: date}

\maketitle

\begin{abstract}
Closely motivated by financial considerations, we develop an integration theory which is not classical i.e. it is not necessarily associated to a measure. The base space, denoted by $\mathcal{S}$ and called a trajectory space, substitutes the set $\Omega$ in probability theory and provides
a fundamental structure via conditional subsets $\mathcal{S}_{(S,j)}$ that allows the definition of conditional integrals. The setting is a natural by-product of no arbitrage assumptions that are used to model financial markets and games of chance (in a discrete infinite time framework). The constructed conditional integrals can be interpreted as required investments, at the conditioning node, for hedging an integrable function, the latter characterized a.e. and in the limit as we increase the number of portfolios used. The integral is not classical due to the fact that the original vector space of portfolio payoffs is not a vector lattice. In contrast to a classical  stochastic setting, where price processes are associated to conditional expectations (with respect to risk neutral measures), we uncover a theory where prices are naturally given by conditional non-lattice integrals. One could then study analogues of classical probabilistic notions in such non-classical setting, the paper stops after defining trajectorial martingales the study of which is deferred to future work.



\vspace{.1in}
\noindent
\subjclass[AMS Classification]~~~{91G20, 60G42, 28C05, 60G48, 60G17}
\end{abstract}

\section{Introduction}

Consider a measurable space $(\Omega, \mathcal{F})$, where $\mathcal{F}= \{\mathcal{F}\}_{0 \leq n \leq T}$ is a discrete time filtration, providing the unfolding of the possible prices of tradable assets in a financial market model. It is well known that  if $\Omega$ has finite cardinality and under a no-arbitrage assumption, prices of attainable functions (which are the ones given by portfolio payoffs) are given  by martingale measures. Non attainable payoffs, which are present in a general incomplete market model, are related to such measures via a superhedging duality result (where superhedging prices
of non attainable payoffs are obtained by considering a supremum over martingale measures of the payoff's expectation). When $\Omega$  has arbitrary cardinality the story retelling the above results is more complex. The standard approach (and in particular the original literature  as in \cite{dalang} and \cite{schachermayer}), starts with a ``physical'' probability measure or a collection of such. Then the martingale measures (used for pricing)  appear as subsidiary to the physical measures. This restriction, namely, pricing measures being dependent on apriori given physical measures, has been essentially removed in recent theoretical results (e.g. \cite{acciaio} and \cite{burzoni}). In this literature, martingale pricing measures
emerge in order to establish a superhedhing duality result, the relationship
between superheding prices and martingale measures is delicate and technical. In this respect, we point out
to \cite{burzoni} where it is shown that one needs to restrict to a subset of $\Omega$ in order to obtain equality in the superheding duality theorem.

In a discrete infinite time and model-free framework,
we take a different perspective to the approaches mentioned above. Suppose we have access to portfolio trading with a riskless bank account and an asset. Consider a set $\mathcal{S}$ (a structured version of the abstract set $\Omega$) prescribing future scenarios for the asset's prices and  assume no knowledge of their occurring probability. We refer to such an approach as a trajectory model and develop a theory that starts with portfolio payoffs (elementary functions) and their initial values (elementary integrals) and construct conditional integrals based on such an elementary data plus the trajectory space. In contrast to the previous description of prices given by
expectations (resulting from no arbitrage considerations and as a mean to obtain superhedging duality results), we develop a closely linked theory between prices and integrals.
The extension of  prices from an elementary setting to an integral acting on integrable functions is performed analogously to  the classical theory (Daniell=Lebesgue integration) but, by necessity, it requires a space of functions that is not a vector lattice. In this way, prices of integrable functions are  naturally given by  non-lattice integrals in such a way that their meaning is associated to prices by hedging, the latter provided
by a limit argument. Namely, integrable functions are given by suitable limits of simple portfolios. The constructions also include null events and conditional integrals. One can interpret the results as showing that, for the basic aspects of a pricing theory, classical expectations and their additional structure (martingales, null sets, etc) are not necessarily inherent to the elementary setting. In fact, there is a more basic non-lattice
integration that relates to prices in close analogy to how the construction of  Lebesgue's theory of integration connects to volumes and extensive quantities. The usual notion of probability is founded in terms of occurrence of events, in our pricing setting, the latter
notion is replaced with superhedging by means of portfolios. The non-lattice integral
captures the additivity property of prices that remains without the lattice property.

Independently of any financial motivation, the paper develops an integration theory
with does not require a probability measure and that relies on a
set of trajectories  $\mathcal{S}$. $\mathcal{S}$ could be visualized as the skeleton of the path space of a stochastic process reflecting basic indeterminism in the structure  of its conditional spaces $\mathcal{S}_{(S,j)}$. The latter are the key to develop  conditional integrals with
a distinct flavour than the classical theory.
We develop this conditional theory of integration with some detail and
expect to explore in the future some results from the classical theory that may persist in the proposed setting. Such a non vector lattice theory of integration
was developed, in an abstract setting, in \cite{leinert} and refined in \cite{konig}. In the particular trajectorial setting of our paper
we had to extend the abstract theory in order  to introduce conditional integrals. The abstract theory depends on two analytical properties which we denote by $(L)$ and $(K)$, they were introduced in \cite{leinert} and \cite{konig} respectively. Property $(K)$  implies  $(L)$  and the latter is the analogue of Daniell's continuity property
in a non-lattice setting. We provide conditions for the validity of $(L)$ and  $(K)$; the proof of property $(L)$ is a key result of our paper and shows how this analytic necessary condition for non-lattice integration results from weak no-arbitrage properties of the space $\mathcal{S}$ and what is the impact
of working with an unbounded time for the validity of $(L)$. Novel concepts, emerging from the trajectory based approach are required in the proof of property $(L)$, this approach was initiated in \cite{ferrando}.  Alternative
proofs of $(L)$ and $(K)$ are presented in \cite{bender}, a technical
difference between the setting of our paper and \cite{bender} is that the latter works mainly  with the set of all possible non-anticipative
portfolios. Availability of certain portfolios does change the
class of null functions, i.e. null sets in our setting are directly related to financial
consideration, this is in contrast to a standard stochastic case.
For reasons of space we have deferred to \cite{bender} for detailed comparisons between the setting of the present paper  and the standard stochastic approach to financial mathematics.

We rely on some informal language used in financial mathematics but
which is used here only
as an aid to motivate definitions and results. The basic object in our developments is a  \emph{trajectory set} $\mathcal{S}$; a \emph{trajectory} $S \in \mathcal{S}$ is an infinite sequence of real numbers
$S = (S_0, S_1, \ldots)$ with common initial value $S_0 = s_0$.
 We also require a set
$\mathcal{H}$; elements $H \in \mathcal{H}$ are called  \emph{portfolios} and are infinite sequences of functions: $H=(H_0, H_1, \ldots)$, $H_i:\mathcal{S} \rightarrow \mathbb{R}$ which are required to be {\it non-anticipative} i.e., whenever $S, S' \in \mathcal{S}$,
with $S_k = S_k', 0 \leq k \leq i$, it then follows that $H_i(S) = H_i(S')$. Such sets $\mathcal{H}$ are called  (non-anticipative) portfolio sets.  $\mathcal{S}$ takes over the role of the abstract set $\Omega$ in stochastic models, we show that elementary properties
of $\mathcal{S}$, along with the possibility to trade, imply the existence of a natural integral functional.

In such general (discrete time) setting, we show how to define price processes
for functions $f$ defined on $\mathcal{S}$ and that are hedged a.e. in the limit, as the number of portfolios growths. The price processes are defined by conditional trajectorial integrals and behave analogously to  martingale processes. In other words, there is an endogenous pricing methodology that owes its existence to a worst case point of view (where every trajectory counts) as contrasted to a probabilistic point of view.

Our application of the theory of no-lattice integration to a pricing setting is new,
we briefly place it in context with the literature. The reference \cite{kassberger} (see also \cite{biagini})  presents an axiomatic approach to  price functionals
and represent them by means of  a classical (i.e. measure based) expectation. Analogously, we follow a related route in that we start with an elementary integral
 representing the prices of simple portfolios that trade on the underlying asset. The prices of such portfolios are, by definition, the capital required to set up the portfolio. Starting with this subspace domain for our pricing operator
 is a key technical difference with the approaches in \cite{kassberger}
 and \cite{biagini}. These papers start with larger domains for their
 pricing operator, those spaces of functions are required to satisfy the crucial property of being vector lattices. In particular, their settings need a large set of financial instruments for which prices are available. It is easy to see, as we show in  Section \ref{integralConstruction}, that this hypothesis does not hold if one restricts to the simplest of settings consisting  of a one dimensional, dynamically traded asset like the one we are considering. In technical terms, the mentioned papers rely on the classical
 theory of integration of the Daniell-Lebesgue type that requires to start with an elementary integral on a set of functions that is a vector lattice. An extended discussion
 of the absence of the lattice property is in \cite{bender}.

Our first task is to construct conditional trajectorial integrals (or {\it conditional integrals} for short). The nonconditional integral emerges as a special case, this is in contrast to the classical approach. Given a finite initial trajectory segment $(S_0, \ldots, S_n)$, conditional integrals act on an appropriate function $f$ providing (essentially) the initial capital necessary to superhedge $f$.  The superhedging becomes exact via a limiting process that increases the number of portfolios used in the superposition. The process is analogous to Eudoxus' method of exhausting the area of a circle by means of covering polygons. The developments make key use of a notion of a.e. given by a countable subadditive norm related to the integral. An event (i.e., a subset of $\mathcal{S}$) is null if its characteristic function can be superhedged with an arbitrarily small initial investment. These events are naturally considered as unlikely as they resemble a lottery that costs arbitrarily little to purchase while its payoff could provide a yield arbitrarily larger than its original buying price.

The paper is organized as follows; Section \ref{integralConstruction}
introduces the trajectorial setting, the elementary integral $I_j$ and the operator $\overline{I}_j$ giving the notion of null functions.
Section \ref{conditionalIntegrals} defines, what would be, the conditional integral operator $\overline{\sigma}_j$; the key condition $(L_j)$ is also introduced. Section \ref{LjConditions} provides conditions on $\mathcal{S}$ for the validity of
$(L_j)$. Section \ref{IntegrabilityCharacterization} provides the definition of conditionally integrable function and an extensive study of integrable functions and classical integration results. In particular,
 in Corollary \ref{restrictedPositive} and  Proposition \ref{characterizationLjK}, we prove the classes of functions $\mathcal{L}_j$ and $\mathcal{L}^K_j$ to be nonnegative and arbitrary integrable functions respectively. We also prove a Beppo-Levi Theorem for $\mathcal{L}_j$.  With the additional
property $(K_j)$, Theorem \ref{punctualCharacterization} establishes that the previously introduced sets of integrable functions actually characterize the integrable functions. A monotone convergence theorem is also proven
along with norm completion characterization of the integrable functions as well as a general Beppo-Levi theorem. Section \ref{basicPropertiesOfConditionalIntegrals} establishes the tower property of conditional integrals which is then used in Section \ref{definitionsOfTM} to illustrate the notion
of trajectorial martingale. Section \ref{conclusions}
summarizes and indicates future work.
Appendix \ref{someProof1} provides several proofs left out in the main sections. Finally, Appendix \ref{proofOfK} provides the intermediate results and arguments to establish
property $(L_j)$.  Under stronger hypothesis than the ones used to obtain $(L_j)$, Theorem \ref{proofOfK_j}, which relies on \cite{bender}, establishes property $(K_j)$. Examples illustrating special cases and the need for the required hypothesis to establish $(L_j)$ are provided. More clarifying examples are also available in \cite{bender}.

\section{Elementary Conditional Integral}  \label{integralConstruction}

The reference \cite{konig} provides an axiomatic approach to a generalized theory of integration, in particular, the integral is not necessarily associated to a measure. K\"{o}nig's  work  resulted from an effort to refine the framework introduced originally by Leinert in \cite{leinert}. Both references dispense with the requirement that the relevant space of functions is a vector lattice, a property which is needed  to fully develop  Daniell's approach to integration.

This theory of integration crucially depends on a property (introduced in \cite{leinert}) we have labelled $(L_j)$ which plays the analogue role of continuity at $0$ in Daniell's theory of integration.
We establish property $(L_j)$ under sufficient weak conditions on the space $\mathcal{S}$ in Section \ref{LjConditions} and  Appendix \ref{proofOfK}. There is also a related, stronger, property (introduced in \cite{konig}), that we have labeled $(K_j)$, which allows for a stronger theory of integration as gives access to some additional convergence theorems and a characterization of the integrable functions. Sufficient conditions to establish $(K_j)$ are in Appendix \ref{proofOfK}.

We work on a specific (trajectorial) setting as in \cite{ferrando}, which is a particular case of K\"{o}nig/Leinert's. On the other hand we concentrate on developing conditional integrals, only unconditional integrals are developed in the mentioned references. That being said, each of our conditional integrals is a special case of the (unconditional) integral of K\"{o}nig/Leinert, it is in
the study on how all these conditional integrals ``hang together" that
properties of $\mathcal{S}$ are put to use.
We indicate the results from K\"{o}nig that we borrow from and at the same time extend. We stay close to the notation in \cite{konig} but depart from it at some points (we provide remarks when we do so).

\subsection{Trajectorial Setting}
\begin{definition}
Given a real number $s_0$, a \emph{trajectory set}, denoted by $\mathcal{S}= \mathcal{S}(s_0)$, is a subset of
\[
\mathcal{S}_{\infty}(s_0) =\{S=(S_i)_{i\ge 0}: S_i\in \mathbb{R},~ S_0=s_0\}.
\]
We make fundamental use of the following \emph{conditional spaces}; for $S \in \Se$ and $j \geq 0$  set:
\begin{equation} \nonumber
\Se_{(S,j)}\equiv\{\tilde{S} \in \Se: \tilde{S}_i= S_i, ~~ 0 \le i \le j\}.
\end{equation}
\end{definition}
In practice, the coordinates $S_i$ are multidimensional in order to allow for multiple sources of uncertainty (\cite{ferrando2}), for simplicity we restrict to $S_i \in \mathbb{R}$. One can also extend the theory to allow for several traded assets $S_i^k$.

Notice $\Se_{(S,0)} = \Se$ and, if $\tilde{S}\in \Se_{(S,j)}$, then $\Se_{(\tilde{S},j)}=\Se_{(S,j)}$. Moreover for $j \le k$ it follows that $\Se_{(S,k)}\subset \Se_{(S,j)}$. On the other hand, for any fixed $j\ge 0$, $\Se$ is the disjoint union of $\Se_{(S,j)}$, since if $S'\notin\Se_{(S,j)}$ then $\Se_{(S,j)}$ and $\Se_{(S',j)}$ are disjoint. For simplicity, sometimes we will refer to the space $\mathcal{S}_{(S,j)}$ through a pair $(S, j)$ with  $S\in\Se$ and $j\ge 0$ which will be called  a {\it node}; {\it local} properties are relative to a given node.

\vspace{.1in} The other basic component are the \emph{portfolios} defined as follows.
\begin{definition} For any fixed $S\in \Se$ and $j\ge 0$, $\He_{(S,j)}$ will be a set of sequences of functions $H = (H_i)_{i \geq j}$,
where $H_i: \Se_{(S,j)} \rightarrow \mathbb{R}$ are non-anticipative in the following sense: for all $\tilde{S},\hat{S}\in\Se_{(S,j)}$ such that $\tilde{S}_k = \hat{S}_k$ for $j \le k\le i$, then $H_i(\tilde{S}) = H_i(\hat{S})$ (i.e. $H_i(\tilde{S}) = H_i(\tilde{S}_0,\ldots, \tilde{S}_i)$).
\end{definition}

Notice that, in general, $\He_{(S,j)}$ does not include all possible sequences of available non-anticipative functions. Instead of indicating what properties of $\He_{(S,j)}$ are required for
particular results we  list below a minimal set of properties which we are applying in this paper  to develop a general theory of conditional non-lattice integration for pricing and superheding:

\begin{itemize}
\item[(H.1)] The sets $\He_{(S,j)}$ are assumed to  be positive cones i.e. $\alpha \He_{(S,j)} + \He_{(S,j)} \subseteq \He_{(S,j)}$ for all $ \alpha \geq 0$.

\item[(H.2)] The portfolios $H^c=(H^c_i)_{i \geq j}$ where $H^c_i$ are constant valued $-1,\;0$ or $1$, on $\Se_{(S,j)}$, are assumed to belong to $\He_{(S,j)}$. So the null portfolio ($\mathbf{0}=(H^{\mathbf{0}}_i\equiv 0)_{i \geq j}$) belongs to $\He_{(S,j)}$.

\item[(H.3)]  If $(H_i)_{i \geq j} \in \mathcal{H}_{(S,j)}$ and $k \geq j$,
then $((H_i)|_{\mathcal{S}_{(S,k)}})_{i \geq k} \in \mathcal{H}_{(S,k)}$.

\item[(H.4)] Let $(H_i)_{i \geq j} \in \mathcal{H}_{(S,j)}$,\;  if $G_i \equiv H_i$ \; for $j \leq i \leq k$\; and  \;$G_i = 0$ for $i >k$, then \\ $G \in \He_{(S,j)}$ as well.
\end{itemize}

Conditions (H.1)--(H.4) are supposed for the rest of this paper without further notice.

\vspace{.1in}$H \in \He_{(S,j)}$ may be referred to as a \emph{ conditional portfolio}.

For a node $(S,j)$, $H\in\He_{(S,j)}$, $V\in \mathbb{R}$ and $n\ge j$ we define the following functions $\Pi_{j,n}^{V, H}: \mathcal{S}_{(S,j)} \rightarrow \mathbb{R}$:
\begin{equation}\label{Pi_notation}
\Pi_{j,n}^{V, H}(\tilde{S}) \equiv V+\sum_{i=j}^{n-1}H_i(\tilde{S})~ \Delta_i \tilde{S},~~~~~~~\tilde{S} \in \mathcal{S}_{(S,j)},
\end{equation}
where $\Delta_i \tilde{S} = \tilde{S}_{i+1}- \tilde{S}_i,~~i\ge j$.

\vspace{.1in}
It will also be convenient to define \emph{global portfolios}. For a fixed $j \geq 0$,  $\He_j$ denotes the set of sequences of functions $H\equiv(H_i)_{i\ge j}$ with $H_i: \Se \rightarrow \mathbb{R}$, where for each $S\in\Se$ there exists $G \in \mathcal{H}_{(S,j)}$ such that $H_i(\tilde{S})= G_i(\tilde{S}) \;\;\forall \tilde{S}\in \Se_{(S,j)}~~\mbox{and}~~i \geq j$. A global portfolio $H$  could be characterized by indicating that its restriction to $\Se_{(S,j)}$ belongs to $\He_{(S,j)}$. We can also define $\He \equiv (\He_j)_{j\ge 0}$ and the pair $(\Se,  \He)$ which may, occasionally, be referred to as a \emph{market}.

\subsection{Elementary Vector Spaces $\mathcal{E}_j$ and Elementary Conditional Integrals $I_j$}\label{Elementary Spaces}

In the sequel, being $\mathcal{A}$ a set of real valued functions, $\mathcal{A}^+$ will denote the set of its non-negative elements, and $\mathcal{A}^-$ the set of non-positive ones.

\vspace{.05in} For a fixed node $(S,j)$  define
\begin{equation} \nonumber 
\mathcal{D}_{(S,j)}= \{f = \Pi_{j, n}^{V, H}: H \in \He_{(S,j)},~~V \in \mathbb{R}~~~\mbox{and}~~~n \in \mathbb{N}\}.
\end{equation}
 Observe that the constant functions on $\mathcal{S}_{(S,j)}$ belong to $\mathcal{D}_{(S,j)}$; from the positive cone assumption for the sets $\He_{(S,j)}$, the spaces
\begin{equation} \nonumber 
\mathcal{E}_{(S,j)} \equiv \mathcal{D}_{(S,j)}^+ - \mathcal{D}_{(S,j)}^+,
\end{equation}
are vector spaces. Elements of $\mathcal{E}_{(S,j)}$ are called {\it elementary functions}. Let also define
\begin{equation} \label{conditionalElementarySpace}
\mathcal{D}_j = \{f:\Se\rightarrow \mathbb{R}: f|_{\Se_{(S,j)}}\in\mathcal{D}_{(S,j)} \;\; \forall S\in\Se \}\quad \mbox{and}\quad \mathcal{E}_j \equiv \mathcal{D}_j^+ - \mathcal{D}_j^+.
\end{equation}
Then, the spaces $\mathcal{E}_j$ are vector spaces as well and its elements are called {\it global elementary functions}. We will use the notation $\mathcal{E} \equiv \mathcal{E}_0$.
The functions in $\mathcal{D}_{(S,j)}$, and in consequence those in $\Ee_{(S,j)}$, will be considered to be in $\mathcal{D}_j$, and $\mathcal{E}_j$ respectively, by considering the embedding of the spaces given by extending the respective functions by $0$ outside $\Se_{(S,j)}$. We may refer to global elementary functions by simply  {\it elementary functions} and let the context
to provide any disambiguation that may be required.


If $f\in \mathcal{D}_j$ then  $f|_{\Se_{(S,j)}}= \Pi_{j, n^S}^{V^S, H^S}$ with $H^S \in \He_{(S,j)}$, $V^S \in \mathbb{R}$ and $n^S \in \mathbb{N}$, we then have the expression:
\[
f(S)=V(S)+\sum_{i=j}^{n(S)-1}H_i(S)\Delta_iS,\quad \mbox{with}\quad H\in\He_j,
\]
where $V(S)=V^S$, $n(S)=n^S$ and for $\tilde{S}\in \Se_{(S,j)}$, $H_i(\tilde{S})=H^S(\tilde{S})$. Naturally, the notation used to for the functions introduced in  (\ref{Pi_notation}) is extended to $f \in \mathcal{D}_j$ as $f = \Pi_{j, n}^{V, H}$ where here $n$ and $V$ are functions of $S$ depending only on $S_0,...,S_j$, i.e. constants on $\Se_{(S,j)}$. This notation will also be further extended to the case that $H$ is a sequence of non-anticipative functions, not necessarily in $\He_j$, as in Corollary
\ref{representationUniqueness}.

In the case  that $\He_j$ is a vector space we have $\mathcal{E}_j^+ = \mathcal{D}_j^+$.
Alternatively to working with the positive cone assumption on $\mathcal{H}_{(S,j)}$, one could proceed with the stronger assumption that $\He_j$ is a vector space and work with $\mathcal{E}_j= \mathcal{D}_j$.

\vspace{.1in} Naturally, we assume that the {\it elementary integral} $I_j: \mathcal{D}_j^+ \rightarrow \mathbb{R}^+$:
\begin{equation} \nonumber
I_j f(S) = V(S)~~ \mbox{for}~~~~ f = \Pi_{j, n}^{V, H} \in \mathcal{D}_j^+,
\end{equation}
is a well defined, linear isotone operator for every $j\geq 0$.

In view of Corollary \ref{representationUniqueness} below, this property is easily seen to be equivalent to the following notion of local $0$-neutrality of $\mathcal{S}$. Therefore, the condition of $\mathcal{S}$ being locally $0$-neutral, will be considered a standing assumption throughout the paper (but we do comment on changes to our constructions
that are needed to weaken the dependency on the $0$-neutral hypothesis).


\begin{definition} [$0$-Neutral nodes] \label{localDefinitions}
Given a trajectory space $\Se$ and a node $(S,j)$:
\begin{itemize}
\item $(S,j)$ is called a \emph{$0$-neutral node} if
\begin{equation} \label{cone}
\sup_{\tilde{S} \in \Se_{(S, j)}}~~~~ (\tilde{S}_{j+1} - S_{j}) \geq 0 \quad \mbox{and} \quad
\inf_{\tilde{S} \in \Se_{(S, j)}} ~~(\tilde{S}_{j+1} - S_{j}) \leq    0.
\end{equation}
\end{itemize}
\noindent $\Se$ is called \emph{locally $0$-neutral} if (\ref{cone}) holds at each node $(S, j)$.
\end{definition}
A stronger condition is given by the following definition.
\begin{definition} [Up-Down nodes] \label{localDefinitions}
Given a trajectory space $\Se$ and a node $(S,j)$:
\begin{itemize}
\item $(S,j)$ is called an \emph{up-down node} if
\begin{equation} \label{upDownProperty}
\sup_{\tilde{S} \in \Se_{(S, j)}}~~ (\tilde{S}_{j+1} - S_{j}) >    0\quad \mbox{and}\quad
\inf_{\tilde{S} \in \Se_{(S, j)}} ~~(\tilde{S}_{j+1} - S_{j}) <    0.
\end{equation}
\end{itemize}
\noindent $\Se$ is called \emph{locally up-down} if (\ref{upDownProperty}) holds at each node $(S, j)$.
\end{definition}

\begin{definition} [Flat nodes] \label{localDefinitions}
Given a trajectory space $\Se$ and a node $(S,j)$:
\begin{itemize}
\item $(S,j)$ is called a \emph{flat node} if
\begin{equation} \label{flat}
\sup_{\tilde{S} \in \Se_{(S, j)}}~~ (\tilde{S}_{j+1} - S_{j}) =    0 =
\inf_{\tilde{S} \in \Se_{(S, j)}} ~~(\tilde{S}_{j+1} - S_{j}).
\end{equation}
\end{itemize}
\noindent $(S, j)$ is called an arbitrage free node if (\ref{upDownProperty}) or (\ref{flat}) hold, otherwise it
is called an arbitrage node. An arbitrage node $(S, j)$ is said to be of type I if there exist $\hat{S}\in\Se_{(S, j)}$
such that $\hat{S}_{j+1} = S_j$, otherwise it is said to be of type II.

\vspace{.1in} \noindent $\Se$ is called \emph{locally free of arbitrage} if each node $(S, j)$ is arbitrage free.
\end{definition}
Note that in a financial context the above notions assume that the real numbers $S_i$ are
represented in units of the traded numeraire (i.e. in "discounted" units).

The following simple Lemma from \cite{ferrando} gives a basic procedure to construct particular useful trajectories.

\begin{lemma} \label{scapingSequence} Assume $\mathcal{S}$ is locally $0$-neutral, $n_0 \geq 0$, and let $F=\{F_i\}_{i\ge n_0}$ be a sequence of non-anticipative functions and $\epsilon>0$. Then, for any $S \in\Se$ there exists a sequence of trajectories $\{S^n\}_{n \ge n_0}$ with $S^{n_0}=S$ 
such that for every $n > n_0$,\\ $S^{n}\in \Se_{(S^{n-1},n-1)}\subset \Se_{(S,n_0)}$ and
\begin{equation} \nonumber 
F_{i}(S^n) \Delta_i S^n < \frac{\epsilon}{2^{i+1}},\; n_0\le i \le n-1,
\end{equation}
and so,
\begin{equation}\label{sequenceInequality}
\sum_{i=n_0}^{n-1} F_{i}(S^n) ~~\Delta_i S^n < \sum_{i=n_0}^{n-1} \frac{\epsilon}{2^{i+1}}.
\end{equation}
\end{lemma}


\begin{remark}  \label{contrarianAtArbitrageNodeOfTypeI}
If at any point in the construction, a node $(S^n, n)$ 
is an arbitrage node of type I we could choose, without loss of generality, $S^{n+1} \in \mathcal{S}_{(S^n, n)}$ such that   $\Delta_{n} S^{n+1}= S_{n+1}^{n+1}- S_{n}^{n+1}=0$.
\end{remark}

\begin{corollary}[$I_j$ is Well Defined]\label{representationUniqueness}
Let $\Se$ be locally $0$-neutral, and  $F, G$ be sequences of non-anticipative functions. The following holds at an arbitrary node $(S,j)$:
\begin{itemize}
\item If $\Pi_{j, n}^{V, F} \geq 0~~\mbox{on}~~\mathcal{S}_{(S,j)}$ then
$\Pi_{j, k}^{V, F}(S) \geq 0 ~\mbox{on}~~\mathcal{S}_{(S,j)}$ for all $j \leq k \leq n$, and so  $V(S) \geq 0$ as well. In particular, this is valid for $\Pi_{j, n}^{V, F}\in\mathcal{D}_j$.

\item If $\;\;\Pi_{j, m}^{U, G}\le \Pi_{j, n}^{V, F}$ then $U(S)\le V(S)$. Consequently $I_j$ is well defined i.e. if $f=\Pi_{j, n^f}^{V^f, H^f}, g=\Pi_{j, n^g}^{V^g, H^g} \in \mathcal{D}^+_j$ and $f|_{\mathcal{S}_{(S,j)}}= g|_{\mathcal{S}_{(S,j)}}$ then $I_j f(S)=V^f(S)=V^g(S)= I_jg(S)$ (see also \cite[Proposition 4]{ferrando}).
\end{itemize}
\end{corollary}
\begin{proof} Assume that for some $j\le k\le n$ there exist $\tilde{S}\in\mathcal{S}_{(S,j)}$ such that $\Pi_{j, k}^{V, F}(\tilde{S})= \delta  < 0$. From Lemma \ref{scapingSequence}, for $n_0=k$, $\tilde{S}$, and $\epsilon= \mbox{\small$-\frac {\delta}2$}$ there exists $S^n\in\Se_{(\tilde{S},k)}$, such that $\sum\limits_{i=k}^{n-1} F_i(S^n)\Delta_i S ^n  < \mbox{\small$-\frac {\delta}2$}$, it then follows the contradiction
\[0\le V(\tilde{S}) + \sum_{i=j}^{n-1} F_i(S^n) \Delta_i S ^n \leq \mbox{\small$\frac {\delta}2$}< 0. \]
\end{proof}

Given $(\Se,\He)$  with $\Se$ locally $0$-neutral, we can then define
$I_j$ on $\mathcal{E}_j$ by:\\ $I_j f = I_j f_1- I_j f_2$ with
$f= f_1- f_2$, $f_1, f_2 \in \mathcal{D}_j^+$. It follows that
$I_j$ is well defined and linear on $\mathcal{E}_j$. By the mentioned embedding of $\mathcal{E}_{(S,j)}$ into $\mathcal{E}_j$, we also can write $I_jf$ for $f\in \mathcal{E}_{(S,j)}$.

\vspace{.1in}
In the general case when $\mathcal{S}$ is locally $0$-neutral
we can  see that $\mathcal{E}$ is not a vector lattice. For example,
consider  $f(S) = (S_1- S_0) \in \mathcal{E}$ and
assume $|f(S)| \in \mathcal{E}$ as well. Therefore $|f(S)|= |S_1- S_0| = V_G + \sum_{i=0}^{n-1} G_i(S) ~\Delta_i S~\forall~S$ for some $\{G_i\}$, a sequence of non-anticipative functions. A similar reasoning as in Corollary
\ref {representationUniqueness} implies  $|f(S)|= |S_1- S_0| = V_G +  G_0(S) ~(S_1- S_0)~\forall~S$. Clearly, it is easy to construct
an example with $S_1$ taking values $S_1^k$, $k=1,2,3$, such that
it is impossible to solve for $V_G$ and $G_0$ in the previous equality
(for further comments see \cite{bender}). In the general case, it follows that $f \in \mathcal{E}$ does not imply $|f| \in \mathcal{E}$ and so the latter is not a vector lattice. This indicates that the well known Daniell-type of integral construction can not be used as it requires the lattice property. The same conclusion applies to each $\mathcal{E}_j$.

\begin{remark}
The property of $\Se$ being locally $0$-neutral is a weak assumption and it is
required for $I_j$ to be well defined. Nonetheless, the hypothesis should  be relaxed in order to incorporate arbitrage
nodes of type II which may not be $0$-neutral. This possible extension of our setting is briefly noted at the end of Section \ref{LjConditions}.
\end{remark}

\subsection{Almost Everywhere Notions}\label{a.e. section}
Let $Q$ denote the set of all functions from $\mathcal{S}$ to $[-\infty, \infty]$ and $P \subseteq Q$ denotes the set of non-negative functions. The following conventions are in effect: $0 ~\infty =0$, $\infty + (- \infty) = \infty$, $u - v \equiv u +(-v) ~\forall~u, v \in [-\infty, \infty]$, and $\inf \emptyset = \infty$.

We define next the \emph{conditional norm operator} \; $\overline{I}_j : P \rightarrow \mathcal{D}^+_j$. 

\begin{definition}\label{Iup_definition}
For a general $f:\mathcal{S} \rightarrow [0, \infty]$ (i.e. $f \in P$) and a given node $(S, j)$ define
\begin{equation} \nonumber
\overline{I}_j f (S)\equiv  \inf_{\{H^m\} \subseteq \He_{(S,j)}} \{\sum_{m \geq 1} V^m: ~~f(\tilde{S}) \leq  \sum_{m \geq 1} f_m(\tilde{S})\;\;\forall \tilde{S} \in \Se_{(S,j)},~f_m = \Pi_{j, n_m}^{V^m, H^m}\in \mathcal{D}_{(S,j)}^+\},
\end{equation}
we will use the notation $\overline{I}f \equiv \overline{I}_0f$.
We also set, for a general $f:\mathcal{S} \rightarrow [-\infty, \infty]$ (i.e. $f \in Q$) and a given node $(S, j)$:
\begin{equation} \nonumber
||f||_j(S) \equiv \overline{I}_j|f|(S)~~\mbox{and}~~||f|| \equiv ||f||_0(S).
\end{equation}
Notice that $\overline{I}_j f (S)\ge 0$, and $\overline{I}_j f (S)= \overline{I}_j f (S_0, \ldots, S_j)$, i.e. $\overline{I}_j f (\cdot)$ is constant on $\mathcal{S}_{(S,j)}$.
\end{definition}\label{conditional_norm}
\noindent
We will refer to $||\cdot||_j(S)$ as the {\it conditional norm(s)}.

The following proposition establishes the properties: isotony (i.e. order preserving), positive homogeneity and countable subadditivity of $\overline{I}$.

\begin{proposition}[see proof in Appendix \ref{someProof1}] \label{propertiesOfIBarra}
Let $(S,j)$ be a fixed node:
\begin{enumerate}
\item Consider $f, g \in P$, if $ f(\tilde{S}) \leq g(\tilde{S}) \;\; \forall \tilde{S} \in \Se_{(S,j)}$, then $\overline{I}_j f(S) \leq \overline{I}_j g(S)$.

\item Consider $f \in P$ and $\alpha \in \mathbb{R}^+$ then
$\overline{I}_j (\alpha f) (S) = \alpha  \overline{I}_j f(S)$.

\item  Consider $g, g_k \in P$, for $k \geq 1$ satisfying $g(\tilde{S}) \leq
\sum_{k \geq 1} g_k(\tilde{S}) \;\; \forall \tilde{S} \in \Se_{(S,j)}$   then
\begin{equation} \nonumber
\overline{I}_j g(S) \leq \sum_{k \geq1 }\overline{I}_j g_k(S).
\end{equation}
\item $\overline{I}_jf\le I_jf$, for $f\in \mathcal{D}_j^+$. Moreover $\overline{I}_j0=0$.
\end{enumerate}
\end{proposition}

\begin{definition}[\bf{Conditional and unconditional a.e. notions}]\label{nullObjects} For $g:\Se\rightarrow [-\infty,\infty]$, $g$ is a \emph{null function} if $\|g\|=0$, a subset $E\subset\Se$ is a \emph{null set} if $\|\mathbf{1}_E\|=0$. Similarly, a subset $E\subset\Se$ is a \emph{full set} if $\|\mathbf{1}_E\|=1$.  A property that holds in the complement of a null set, is said to hold ``almost everywhere" (a.e.).

We now provide the conditional versions of the previous a.e. notion. Given a node $(S,j)$ a function ~~$g$~~is a \emph{conditionally null function at} $(S,j)$  if:
\begin{equation} \nonumber
\|g\|_j(S)=0.
\end{equation}
 A subset $E\subset\Se$ is a \emph{conditionally null set at $(S,j)$} if $\|\mathbf{1}_E\|_j(S)=0$. Given a node $(S,j)$, a property is said to hold conditionally a.e. at $(S,j)$
(or equivalently: the property holds ``a.e. on $\mathcal{S}_{(S,j)}$") if the subset of  $\Se_{(S,j)}$ where it does not hold
is a conditionally null set at $(S,j)$. In particular, the latter definition applies to $g=f$ a.e. on $\mathcal{S}_{(S,j)}$.
\end{definition}

All appearing equalities and inequalities are valid for all points in the spaces where the functions are defined unless qualified by an explicit a.e.


\begin{proposition}[see proof in Appendix \ref{someProof1}]\label{leinertTheorem} Consider $f,g:\Se\rightarrow [-\infty,\infty]$ and a fixed node $(S,j)$, then:
\begin{enumerate}
\item $\|g\|_j(S) = 0$ iff $g = 0$  a.e. on $\mathcal{S}_{(S,j)}$.
\item $\|g\|_j(S) < \infty$ ~ then ~ $|g| < \infty$ a.e. on $\mathcal{S}_{(S,j)}$.
\item If $|f| \leq |g|$ a.e. on $\mathcal{S}_{(S,j)}$ then $\|f\|_j(S) \leq \|g\|_j(S)$. Therefore, if $|f| = |g|$ a.e. on $\mathcal{S}_{(S,j)}$ then $\|f\|_j(S) = \|g\|_j(S)$.
\item The countable union of conditional null sets is a conditionally null set.

\item For $f \in P$ and $0 \leq j \leq k$: $0 \leq \overline{I}_j(\overline{I}_k f) \leq \overline{I}_j f$.
Therefore, if $g \in Q$ is conditionally null
at $\mathcal{S}_{(S,j)}$ then $\overline{I}_k(|g|)=0 $ a.e. on $\mathcal{S}_{(S,j)}$.
\end{enumerate}
\end{proposition}

\vspace{.1in}

\section{Conditional Integrals}\label{conditionalIntegrals}

\vspace{.1in} Here we introduce the \emph{conditional integration operator} $\overline{\sigma}_j : Q \rightarrow \mathcal{D}_j$. 

\begin{definition}\label{cond_integ_def} For a node $(S,j)$ of $(\Se,\He)$ and a general $f \in Q$,
\begin{equation} \nonumber
\overline{\sigma}_j f(S) \equiv  \inf_{\{H^m\} \subseteq \He_{(S,j)}} \{\sum_{m \geq 0}V^m: ~~f(\tilde{S}) \leq  \sum_{m \geq 0} f_m(\tilde{S}) \;\;\forall \tilde{S} \in \Se_{(S,j)},
\end{equation}
\begin{equation} \nonumber
~~~~f_0 = \Pi_{j, n_0}^{V^0, H^0}\in \mathcal{D}_{(S,j)}^-, f_m= \Pi_{j, n_m}^{V^m, H^m} \in \mathcal{D}_{(S,j)}^+ ~\mbox{for}~ m \geq 1\}.
\end{equation}
Define also $\underline{\sigma}_j f(S) \equiv -\overline{\sigma}_j(-f) (S)$.
We will use the notation $\overline{\sigma}f \equiv \overline{\sigma}_0 f$.
\end{definition}
Under the assumption of $0$-neutrality of $\mathcal{S}$, $I_j$ is well defined hence in the Definitions \ref{Iup_definition} and \ref{cond_integ_def}, $V^m=I_jf_m(S)$ when $f_m$ is considered
as an element of $\mathcal{D}_j^+$ or $\mathcal{D}_j^-$.  It then follows from Corollary \ref{representationUniqueness} that  if we replace the elementary functions $f_m \in \mathcal{D}_{(S,j)}^+$ for global elementary functions $f_m \in \mathcal{D}_j^+$ in Definitions \ref{Iup_definition} and \ref{cond_integ_def}, the values of $\overline{I}_jf(S)$ and $\overline{\sigma}_jf(S)$, respectively, are unchanged. This fact will be used implicitly in several computations.

\begin{remark}
\noindent Note that $\overline{\sigma}_j f(S)= \overline{\sigma}_j f(S_0, \ldots, S_j)$, also that, if $\mathcal{S}$ is $0$-neutral,   $\sum_{m \geq 0}V^m = \sum_{m \geq 0}I_jf_m(S)$.
\end{remark}

The next proposition summarizes some properties of $~\overline{\sigma}_j$, the corresponding properties for $\underline{\sigma}_j$ follow through (and so are not stated).

\begin{proposition}[1.4 i, 1.3 Eigenchaften vi, vii, in \cite{konig}; see proof in Appendix \ref{someProof1}]\label{konig} Unless indicated otherwise, consider $f, g \in Q$ and let $(S,j)$ denote a generic node.
\begin{enumerate}
\item $\overline{\sigma}_j (f+g)(S) \le \overline{\sigma}_j f(S) +
 \overline{\sigma}_j g(S)$.

\item $\overline{\sigma}_j f(S) \le \overline{I}_jf(S)$ if $f\in P$. \label{sigma_le_Iup}

\item $\overline{\sigma}_j f(S) \le I_jf(S)$ if $f\in \mathcal{E}_j$.

\item If  $f \leq
g ~a.e.$ on $\mathcal{S}_{(S,j)}$, then
$\overline{\sigma}_j f(S) \leq \overline{\sigma}_j g(S)$. Therefore
If  $f =
g ~a.e.$ on $\mathcal{S}_{(S,j)}$, then
$\overline{\sigma}_j f(S) = \overline{\sigma}_j g(S)$.

\item $\overline{\sigma}_j(g  f)(S)= g(S)~\overline{\sigma}_j(f)(S)$ ~~~~ if ~~~~ $g(S)=g(S_0,\ldots,S_j)>0$.

\item $\overline{\sigma}_j f(S) \le \overline{\sigma}_j(|f-g|)(S)+ \overline{\sigma}_jg(S)$.

\item \label{negative_sigma} Either: a) $\overline{\sigma}_j0=0$ or b) $A \equiv \{\overline{\sigma}_j0 <0\} \neq \emptyset$ and  $\overline{\sigma}_j f(S) = \pm \infty$ for all $S \in A$ and for all $f \in Q$.


\item If $\overline{\sigma}_jf(S)<\infty$ then $f<\infty$ conditionally a.e. at $(S,j)$.

\item (Bemerkung 1.5 in \cite{konig}) If $\overline{\sigma}_j 0=0,~~~\overline{\sigma}_j f(S) < \infty$ ~~ and ~~ $\underline{\sigma}_j f(S) > -\infty$, then:  \\ $\overline{\sigma}_j f(S) > -\infty$ ~and~ $\underline{\sigma}_j f(S) <  \infty$ ~and $|f|< \infty$ a.e. on $\mathcal{S}_{(S,j)}$.
\end{enumerate}
\end{proposition}

We now introduce several properties that will be proven equivalent in Proposition \ref{Properties_L}.

The property ($L_j$) below generalizes a non-conditional version from  \cite{leinert} and will be called (conditional) \emph{continuity from below}.
\begin{equation}\label{Lk}
(L_j):\quad f \leq \sum_{m \geq 1} f_m, ~~~~f \in \mathcal{E}_j, f_m \in \mathcal{D}_j^+~~~ \implies I_jf \leq \sum_{m \geq 1} I_j f_m.
\end{equation}
It will be convenient to localize $(L_j)$ to a given node, more precisely
we will say that the property $(L_{(S,j)})$ holds at node $(S,j)$
if whenever the left hand side of (\ref{Lk}) is replaced by  $f \leq \sum_{m \geq 1} f_m,$ (the inequality holding on $\mathcal{S}_{(S,j)}$), $ ~~~~f \in \mathcal{E}_{(S, j)}, f_m \in \mathcal{D}_{(S, j)}^+~$, then the right hand side of (\ref{Lk}) holds at $S$. Note that although most of the results involving property $(L_j)$ are globally stated, it can be seen from the proofs that they are valid pointwise.

$(L_j)$ is essential for the theory of integration we develop;
under appropriate conditions we establish $(L_j)$ in Section \ref{LjConditions}.
 Under $(L_j)$ we obtain $\overline{I}_j f= I_j f$
for $f \in \mathcal{D}^+_j$.

For the sake of the next proposition we also introduce:
\begin{equation}\nonumber
(L_{j, +}):\quad f \in \mathcal{E}_j~~~~ \implies  I_j f \le \overline{I}_j f^+, 
\end{equation}
as  well as:
\begin{equation}\nonumber 
(L_{j}^+):\quad \quad f \leq \sum_{m \geq 1} f_m, \;\; f , f_m \in \mathcal{D}_j^+~~~~ \implies  I_j f \le \sum_{m \geq 1} I_j f_m.
\end{equation}
Notice that $(L_j) \implies (L_j^+)$ is immediate.

\begin{proposition}[see proof in Appendix \ref{someProof1}] \label{Properties_L}
The following items are equivalent.
\begin{enumerate}
\item $\overline{\sigma}_j 0 = 0$.
\item Property $(L_j)$.
\item Property $(L_{j, +})$.
\item \label{L_j}$\underline{\sigma}_j f = I_j f= \overline{\sigma}_j f$ for $f\in \mathcal{E}_j$.
\item $|I_j f| \le \|f\|_j$ for $f\in \mathcal{E}_j$.
\item $I_j f = \overline{I}_j f$ for $f\in \mathcal{D}^+_j$.
\item Property $(L_{j}^+)$.
\end{enumerate}
\end{proposition}

\begin{corollary}\label{good_properties} Assume $\overline{\sigma}_j0=0\;$ for a fixed $j\geq 0$. Then for $f \in Q$ and $S\in\Se$:
\begin{enumerate}
\item  $0 \le \overline{\sigma}_j f + \overline{\sigma}_j(-f)$ and $\underline{\sigma}_j f \leq \overline{\sigma}_j f$.
\item (Bemerkung 1.4 in \cite{konig})  $|\overline{\sigma}_j f| \leq \overline{\sigma}_j |f|$.
\item If $f=0\;\;a.e.~~~on~~~\mathcal{S}_{(S,j)}$ then $\underline{\sigma}_j f(S) = \underline{\sigma}_j |f|(S) = 0 = \overline{\sigma}_j |f|(S)= \overline{\sigma}_j f(S)$, which in turn implies $\underline{\sigma}_j f^+(S)=\underline{\sigma}_j f^-(S) = 0 = \overline{\sigma}_j f^-(S) = \overline{\sigma}_j f^+(S)$.
\item $0 \leq g~a.e.$ on $\mathcal{S}_{(S,j)}$ then $0 \leq \underline{\sigma}_j g(S)$.
\end{enumerate}
\end{corollary}
\begin{proof}
For item $(1)$, $0 \le \overline{\sigma}_j f + \overline{\sigma}_j(-f)$ follows from Proposition \ref{konig}, item $(1)$, and $\overline{\sigma}_j0=0$. Furthermore, if $|\overline{\sigma}_j(-f)|< \infty$, it follows that $\underline{\sigma}_jf =-\overline{\sigma}_j(-f) \leq  \overline{\sigma}_j f$. While if $\overline{\sigma}_j(-f)=\infty$, then $\underline{\sigma}_jf =-\infty\le\overline{\sigma}_jf$. Finally if $\overline{\sigma}_j(-f)=-\infty$, then $0 \le \overline{\sigma}_j f + \overline{\sigma}_j(-f)$ implies $\overline{\sigma}_jf =\infty =-\overline{\sigma}_j(-f)=\underline{\sigma}_jf$.


For item $(2)$, from $-|f|\le f\le |f|$ it follows that  $\overline{\sigma}_j(-|f|)\le \overline{\sigma}_jf\le \overline{\sigma}_j|f|$. From item $(1)$:
$-\overline{\sigma}_j|f|\le \overline{\sigma}_j(-|f|)$
and so $|\overline{\sigma}_j f| \leq \overline{\sigma}_j |f|$.\\
For item $(3)$; from Proposition \ref{konig} item $(\ref{sigma_le_Iup})$: $0 \le \overline{\sigma}_j |f|(S)\le \overline I_j |f|(S)=0$ (where the last equality holds by hypothesis). Then $\overline{\sigma}_j f(S)=0$ by item $(2)$, the latter  applied to $-f$ gives $\underline{\sigma}_j f(S)=0$. Moreover, since also $|f|=0\;\;a.e.~~~on~~~\mathcal{S}_{(S,j)}$ 
the first chain of equalities hold.  For the remaining statements, it is enough to observe that $0\le f^-,f^+\le|f|$ which gives $f^-,f^+ =0\;\;a.e.~~~on~~~\mathcal{S}_{(S,j)}$ and so the same previous reasoning that we used for $f$ applies to $f^+$ and $f^-$ as well.\\
For item $(4)$, from Proposition \ref{konig} item $(4)$, it follows that
$0 = \underline{\sigma}_j 0(S) \leq \underline{\sigma}_j g(S) $.
\end{proof}

\section{Conditions for validity of property $(L_j)$}\label{LjConditions}
Item (\ref{negative_sigma}) of Proposition \ref{konig} shows that property $(L_j)$ is crucial to have an appropriate conditional integral defined by $\overline{\sigma}_j$. In order to obtain conditions for its validity, some new concepts need to be introduced.

Given a sequence $\{S^n\}_{n\ge 0} \subseteq \Se_{(S,j)}$ satisfying
\begin{equation}  \nonumber 
~~~S^n_i = S^{n+1}_i \;\; 0 \leq i \leq n, ~~\forall ~n,~~~~~~~~
\end{equation}
define
\begin{equation}  \nonumber
\overline{S}= \{\overline{S}_i\}_{i \geq 0}~~\mbox{by}~~~\overline{S}_i \equiv S^i_i.~~~
\mbox{We will use the notation}~~\overline{S}= \lim_{n \rightarrow \infty} S^n.
\end{equation}
Notice that $\overline{S}_i = S_i,~0 \leq i \leq j$ given that $S^i \in \mathcal{S}_{(S,j)}$. Moreover
\begin{equation} \label{limitSequence}
\overline{S}_i = S^n_i, ~~~0 \leq i \leq n,  ~~\forall ~n \geq 0.
\end{equation}
Examples of such sequences are obtained from Lemma \ref{scapingSequence}.

\vspace{.1in}
Let $\overline{\Se_{(S,j)}}$ be the set of such $\overline{S}$,
clearly $\Se_{(S,j)} \subseteq \overline{\Se_{(S,j)}}$ given that for $\tilde{S} \in \Se_{(S,j)}$ we can take $\tilde{S}^n = \tilde{S}$ for all $n \geq 0$. We say that $\mathcal{S}_{(S,j)}$ is {\it complete} if
$\mathcal{S}_{(S,j)} = \overline{\mathcal{S}_{(S,j)}}$.

\begin{remark}
The process of going from $\mathcal{S}_{(S,j)}$ to   $\overline{\mathcal{S}_{(S,j)}}$
does not alter the properties of nodes (being $0$-neutral, no arbitrage, etc), see \cite{ferrando} for some details.
Moreover, portfolios $\Pi_{j,n}^{V, H}$ can be extended, in an obvious way, from acting on $\mathcal{S}_{(S,j)}$
to act on  $\overline{\mathcal{S}_{(S,j)}}$; we will freely make use of this fact below without further comments.
\end{remark}

\begin{definition}[Reversed Fatou Property]
We will say that $\mathcal{S}_{(S,j)}$ satisfies the Reversed Fatou Property (RFP, for short) if:
for any $f_0 \in \mathcal{E}_{(S, j)}, f_m= \Pi^{V^m, H^m}_{j, n_m} \in \mathcal{D}^+_{(S, j)}$, $m \geq 1$,
satisfying $\sum_{m \geq 1} V^m(S) < \infty$, and $\overline{S} \in \overline{\mathcal{S}_{(S,j)}}$  there exists
at least one $\{S^n\}_{n \geq 0} \subseteq \mathcal{S}_{(S,j)}$ such that
$\overline{S} = \lim_{n \rightarrow \infty} S^n$ and
\begin{equation} \label{rFP}
\sum_{m \geq 0} \limsup_{n \rightarrow \infty} f_m(S^n) \geq
\limsup_{n \rightarrow \infty} \sum_{m \geq 0}  f_m(S^n).
\end{equation}
\end{definition}
We can substract
$\sum_{m \geq 0} V^m(S)$ from both sides in (\ref{rFP});
the resulting alternative inequality is equivalent to (\ref{rFP});
we will use this observation at some points below.

\begin{remark} \label{completenessImpliesRFP}
$\sum_{m \geq 0} \limsup_{n \rightarrow \infty} f_m(S^n)=
\sum_{m \geq 0} \limsup_{n \rightarrow \infty} f_m(S^n_0, \ldots, S^n_{n_m})= \\
\sum_{m \geq 0}  f_m(S^n_0, \ldots, S^{n_m}_{n_m})=
\sum_{m \geq 0}  f_m(\overline{S})$.
Thus, Fatou's lemma in conjunction with (\ref{rFP}) implies
$\lim_{n \rightarrow \infty} \sum_{m \geq 0}  f_m(S^n)=
\sum_{m \geq 0}  f_m(\overline{S})$.

Clearly (\ref{rFP})
holds if $\overline{S} \in \mathcal{S}_{(S,j)}$ by taking $S^n \equiv \overline{S}$ for all $n$, hence (\ref{rFP}) is a condition
on elements of $\overline{\mathcal{S}_{(S,j)}} \setminus \mathcal{S}_{(S,j)}$. In particular, if $\mathcal{S}_{(S,j)}$ is complete
it then satisfies the RFP. An example in Appendix \ref{proofOfK} shows that $RFP$ is weaker than completeness.
\end{remark}

 Note that RFP is a local (i.e. nodewise) property which depends on the trajectory set and on the sets of admissible portfolios. We may also consider the following global version (i.e. defined only at the initial node $(S,0)$) which is independent of the portfolio sets:
 \begin{definition}
  We will say that $\mathcal{S}$ satisfies the Global Reversed Fatou Property (GRFP, for short) if for every sequence $(V^m)_{m\geq 1}$ of nonnegative reals, $(n_m)_{m\geq 1}$ of nonnegative integers, and $(H^m)_{m\geq 1}$ of nonanticipating sequences $H^m_i:\mathcal{S} \rightarrow \mathbb{R}$ such that $\Pi^{V^m, H^m}_{0, n_m}\geq 0$
  the following holds:
If  $\sum_{m \geq 1} V^m < \infty$ and $\overline{S} \in \overline{\mathcal{S}}$ then there exists
at least one $\{S^n\}_{n \geq 0} \subseteq \mathcal{S}$ such that
$\overline{S} = \lim_{n \rightarrow \infty} S^n$ and
\begin{equation} \label{grFP}
\sum_{m \geq 1} \limsup_{n \rightarrow \infty} f_m(S^n) \geq
\limsup_{n \rightarrow \infty} \sum_{m \geq 1}  f_m(S^n).
\end{equation}
 \end{definition}

\begin{remark}\label{rem:GRFP}
 We will show in Proposition \ref{validityOfHypothesisInLjTheorem} that completness of $\mathcal{S}$ implies GRFP which in turn yields RFP at any node $(S,j)$ for any choice of the portfolio sets satisfying (H.1)--(H.4).
\end{remark}

\begin{lemma} \label{priceMinusInfinity}
Given a trajectory set $\mathcal{S}$ consider  a  node $(S,j)$, $j \geq 0$, then:
\begin{itemize}
\item If $(S,j)$ is a type II arbitrage node:
\begin{equation}\nonumber 
\overline{\sigma}_jf (S) = - \infty~\mbox{for any}~~f \in Q.
\end{equation}
\item If $\overline{\sigma}_j0(S)=0$ then, $(S,j)$ is $0$-neutral.
\end{itemize}
\end{lemma}
\begin{proof}
We may consider the case when $\tilde{S}_{j+1} > S_j$ for all $\tilde{S} \in \mathcal{S}_{(S,j)}$. Take then, for all $m \geq 1$: $H^m_j(\tilde{S})=1$ and $H^m_i(\tilde{S})=0$
for all $i >j$, $V^m =0$. Also, $H^0_i=0$ for all $i \geq j$; then, for any $V^0 \in \mathbb{R}$:
\begin{equation} \nonumber
f(\tilde{S}) \leq V^0 + \infty = V^0 + \sum_{m \geq 1} H^m_j(\tilde{S})~ \Delta_j \tilde{S}~~\hspace{.1in} \mbox{holds for any}~~\tilde{S} \in \mathcal{S}_{(S,j)}~~\mbox{and}~~f \in Q,
\end{equation}
the first claim then follows.

Assume now $\overline{\sigma}_j0(S)=0$, therefore $(S,j)$ can not be a type II node, this implies that it is a $0$-neutral node.
\end{proof}
Because of the above lemma, we will exclude explicitly
nodes of type II in upcoming results.
We will show that by excluding type II nodes we can then prove
$\overline{\sigma}_j0(S)$ at all nodes $(S,j)$.

\begin{theorem}[Establishing property $(L_j)$] \label{proofOfL}
Assume $\mathcal{S}_{(S,j)}$ satisfies the RFP for each $S \in \mathcal{S}$ and each $j \geq 0$ and nodes are of no-arbitrage type or arbitrage nodes of type I. Then, the property $(L_{(S,j)})$ holds  at each node $(S,j)$.
\end{theorem}
Note that in view of Remark \ref{rem:GRFP} the conclusion of Theorem \ref{proofOfL} is valid, if $\mathcal{S}$ has no arbitrage nodes of type II and satisfies GRFP or completeness.

The proof of Theorem \ref{proofOfL}  is postponed to Appendix \ref{proofOfLj} as it requires some technical steps.


\vspace{.1in}
\noindent
\fbox {
    \parbox{5.5in}
    {
    {\it Standing assumption $(L_j)$: in the remaining of the paper  ({\it but not necessarily in the appendices}) it will be assumed that} $\overline{\sigma}_j(0)=0$ for all $j \geq 0.$
    }
    }

\vspace{.15in}
Lemma \ref{priceMinusInfinity} shows that the standing assumption prevents the existence of type II arbitrage nodes while at the same time  ensures $0$-neutrality.

At the expense of technical complications, it is possible to generalize our results to include arbitrage type II nodes. They lead to null sets and because of this
fact, the above standing assumption could
be weakened to $(L_j)$ holding a.e.

Besides property $(L_j)$, our theory of integration will require another fundamental property, based on a homologous condition introduced by K\"onig \cite{konig}, that we have labelled $(K_j)$ and will be discussed at due time. $(K_j)$ plays a role in later developments of the theory and hence its introduction is delayed. Section \ref{establishingKj}
in Appendix \ref{proofOfK} provides conditions for the validity of $(K_j)$ for all $j \geq 0$.

\section{Integrable Functions and Their Characterization}\label{IntegrabilityCharacterization}

\vspace{.15in}
The following are our conditional integrable functions:

\begin{equation} \label{conditional_integrable}\nonumber
\mathcal{L}_{\overline{\sigma}_j} \equiv \{f:\mathcal{S} \rightarrow [- \infty, \infty]:~~\overline{\sigma}_j f - \underline{\sigma}_jf=0 ~~a.e.\}.
\vspace{.1in}
\end{equation}
Due to the fact that by Corollary \ref{good_properties} item $(2)$, $0 \leq \overline{\sigma}_j f(S) + \overline{\sigma}_j(-f)(S)=\overline{\sigma}_j f(S) - \underline{\sigma}_jf(S)$ holds for all $S$, we have: $f \in \mathcal{L}_{\overline{\sigma}_j}$ if and only if $ \overline{\sigma}_j f - \underline{\sigma}_jf \leq 0~~ a.e.$ which, in turn, implies $\overline{\sigma}_j f \leq \underline{\sigma}_j f$ a.e., moreover, the statements are equivalent if and only if $\underline{\sigma}_j f$ and $\overline{\sigma}_j f$ are finite a.e. Therefore, $f \in \mathcal{L}_{\overline{\sigma}_j}$ implies $- \infty < \overline{\sigma}_j f(S) < \infty~$ and $- \infty < \underline{\sigma}_jf(S) < \infty~$, each set of inequalities holding a.e. Proposition \ref{konig}, item $(9)$ implies
that for each such $S$ we have $- \infty < f < \infty$
a.e. on $\mathcal{S}_{(S,j)}$.

\begin{remark}
For the special case of $j=0$ one can drop the qualifier a.e. given that
$ \overline{\sigma}_0 f - \underline{\sigma}_0f \leq 0~~ a.e.$
if and only if  $ \overline{\sigma}_0 f (S_0)- \underline{\sigma}_0f(S_0) \leq 0$.

At this stage of the developments we could have defined  $\mathcal{L}_{\overline{\sigma}_j}$ by requiring $\overline{\sigma}_j f(S) - \underline{\sigma}_j f(S)= 0~~\mbox{for all}~~S$. The addition of the a.e. in the actual definition plays a role in Corollary \ref{integrableImpliesConditionalIntegrable} and hence the
extra flexibility was added to the definition of $\mathcal{L}_{\overline{\sigma}_j}$.

\end{remark}

$f \in \mathcal{L}_{\overline{\sigma}_j}$ will be called a (conditionally) {\it integrable function}, for such a function
we set the notation:
\begin{equation} \nonumber
\int_j f \equiv \overline{\sigma}_j f,
\end{equation}
we have defined our conditional integral everywhere on $\mathcal{S}$
to avoid the introduction of classes of equivalence of functions defined a.e. That being said, we emphasize  that for $f \in \mathcal{L}_{\overline{\sigma}_j}$,  $- \infty < \int_jf(S)= \underline{\sigma}_jf(S)= \overline{\sigma}_jf(S)< \infty$ a.e. We also write:
\begin{equation} \nonumber
\int_j f \dot{=} \underline{\sigma}_j f,~~~~~ ~\mbox{where we will use $\dot{=}$ to denote equality a.e.}
\end{equation}
We will mix the equivalent uses of a.e. and  $\dot{=}$
as we see most convenient for display purposes. The case $j=0$ is denoted by  $\int f \equiv \int_0 f$  which is a constant defined everywhere on $\mathcal{S}$.

The next lemma pursues the linear property
of the integral at a point $S$ where
the necessary integrability conditions hold. This shows that the said property is local in the specified sense.

\begin{lemma} \label{linearityOfSigma} Let $f, g \in Q$ and consider
a fixed  $S\in\Se$. If:  $\overline{\sigma}_j f(S) - \underline{\sigma}_j f(S) = 0 = \overline{\sigma}_j g(S) - \underline{\sigma}_j g(S)$, then all the involved quantities are finite and
\[
\hspace{-.67in}(a)\quad\overline{\sigma}_j(cf)(S)= c\overline{\sigma}_jf(S)=c\underline{\sigma}_jf(S)=\underline{\sigma}_j(cf)(S)\quad\forall c\in \mathbb{R}, 
\]
\[
\hspace{.2in} (b)\quad\overline{\sigma}_j(f+g)(S)= \overline{\sigma}_jf(S)+\overline{\sigma}_jg(S)= \underline{\sigma}_jf(S)+\underline{\sigma}_jg(S)=\underline{\sigma}_j(f+g)(S). 
\]
\end{lemma}
\begin{proof}
The finiteness claims follow from our conventions in the first paragraph of Section \ref{a.e. section}. We then see that the  hypotheses imply that $\overline{\sigma}_j f(S)= \underline{\sigma}_j f(S)$ and $\overline{\sigma}_j g(S)= \underline{\sigma}_j g(S)$.

For (a), if $c=0$ or $c=-1$ the result is clear. For $c>0$ it follows from item $(5)$ of Proposition \ref{konig}, from where, if $c<0$
\[\overline{\sigma}_j(cf)(S)= \overline{\sigma}_j(-c(-f))(S)= -c\overline{\sigma}_j(-f)(S)= c\overline{\sigma}_jf(S).\]
(b) holds from
\[
\overline{\sigma}_jf(S)+\overline{\sigma}_jg(S)=\underline{\sigma}_jf(S)+\underline{\sigma}_jg(S)\le \underline{\sigma}_j[f+g](S)\le \overline{\sigma}_j[f+g](S)\le\overline{\sigma}_jf(S)+\overline{\sigma}_jg(S).
\]
\end{proof}

\begin{remark} For $f\in Q$ and $S\in\Se$ fixed,  $- \infty < \overline{\sigma}_j f(S)=\underline{\sigma}_j f(S)<\infty$ is equivalent to $\overline{\sigma}_j f(S) - \underline{\sigma}_j f(S) = 0$. This will be used implicitly several times to apply Lemma \ref{linearityOfSigma} or to justify that $f\in \mathcal{L}_{\overline{\sigma}_j}$.
\end{remark}

\begin{proposition}[K\"{o}nig's Behauptung 2.1]\label{Konig2.1}
\noindent
\begin{enumerate}
\item $\mathcal{E}_j\subset \mathcal{L}_{\overline{\sigma}_j}$.

\item If $f \in \mathcal{L}_{\overline{\sigma}_j}$ and $[f=g$ a.e. at $(S,j)]$ holds a.e.
then:\\ $g \in \mathcal{L}_{\overline{\sigma}_j}$ and
$\int_j f \dot{=} \int_j g$ a.e.

\item Consider $f \in \mathcal{L}_{\overline{\sigma}_j}$ and $c \in \mathbb{R}$. Then $c f \in \mathcal{L}_{\overline{\sigma}_j}$ and:
 $\int_j c f = c \int_j f$ if $c \geq 0$ and $\int_j c f = (-c) \int_j (-f) \dot{=}~ c~ \int_j  f $ if $c < 0$.

\item If $f, g \in \mathcal{L}_{\overline{\sigma}_j}$, then $f+g \in \mathcal{L}_{\overline{\sigma}_j}$ and $\int_j  (f+g) \dot{=}  \int_j f+ \int_j g$.
\end{enumerate}
\end{proposition}
\begin{proof}
Assertion $(1)$ follows from Proposition \ref{Properties_L} item $(4)$. \\
For assertion $(2)$, notice that Proposition \ref{konig} item $(4)$ implies
$\overline{\sigma}_jf(S)= \overline{\sigma}_j g(S)$ and this holds a.e. (in $S$). The same reasoning applied to $-f, -g$ gives the validity of
$\underline{\sigma}_jf(S)= \underline{\sigma}_j g(S)$ a.e. Given that $\overline{\sigma}_jf(S) - \underline{\sigma}_jf(S)=0$ holds a.e, it follows that $\overline{\sigma}_jg(S) - \underline{\sigma}_jg(S)=0$ holds a.e. and $g \in \mathcal{L}_{\overline{\sigma}_j}$ is then established.
Notice that $\int_j f \equiv \overline{\sigma}_j f \dot{=}  \overline{\sigma}_j g
\equiv \int_j g$ which completes the proof of $(2)$.

For $(3)$, $f \in \mathcal{L}_{\overline{\sigma}_j}$
gives $\overline{\sigma}_j f -\underline{\sigma}_j f\dot{=}0$ then, $(a)$ from Lemma \ref{linearityOfSigma} implies  $\overline{\sigma}_j(c f) - \underline{\sigma}_j(c f)\dot{=}0$. Hence $c f \in \mathcal{L}_{\overline{\sigma}_j}$.
 For $c \geq 0$ $\int_j c f=c \int_j f$ follows from Proposition \ref{konig} item $(5)$ and our standing assumption, on the other hand if $c <0$: $\int_j cf = \overline{\sigma}_j cf = |c|~ \overline{\sigma}_j(-f)= -c \int_j (-f)=
c \underline{\sigma}_j f \dot{=}~ c \overline{\sigma}_j f = ~c \int_j f$.

Finally, to establish $(4)$, $f, g \in \mathcal{L}_{\overline{\sigma}_j}$
gives $\overline{\sigma}_j f - \underline{\sigma}_j f\dot{=}0$ and $\overline{\sigma}_j g - \underline{\sigma}_j g\dot{=}0$ then, $(b)$ from Lemma \ref{linearityOfSigma} implies  $\overline{\sigma}_j(f+g) - \underline{\sigma}_j(f+g) \dot{=} 0$ and so $(f+g) \in \mathcal{L}_{\overline{\sigma}_j}$.
Moreover, the same result also gives $\int_j (f+g)= \overline{\sigma}_j (f+g)\dot{=}\overline{\sigma}_j f + \overline{\sigma}_j g = \int_j f +\int_j g $.
\end{proof}

\begin{remark} If  $[g=0\;\;a.e.~\mbox{at}~ (S,j)]$ holds a.e.,
item $(2)$ above gives $g\in \mathcal{L}_{\overline{\sigma}_j}$
and $\int_j g \dot{=} 0$. Moreover,
 Corollary \ref{good_properties} item $(3)$ implies also that $|g|, g^+, g^-\in \mathcal{L}_{\overline{\sigma}_j}$ and that all the conditional integrals are zero a.e.
\end{remark}

\begin{definition} For $ j \geq 0$ given, define the following set of functions:
\begin{equation} \nonumber
\mathcal{M}_j\equiv \{v = \sum_{m=1}^{\infty} f_m,~~~~f_m \in \mathcal{D}_j^+~~~~\forall~m \geq 1,~~~~~[\sum_{m=1}^{\infty} I_jf_m < \infty,~a.e.]\}.
\end{equation}
$\mathcal{M}\equiv \mathcal{M}_0$ (denoted by $M(I|\mathcal{E})$ in \cite{konig}). Moreover for fixed $(S,j)$ define
\begin{equation} \nonumber
\mathcal{M}_j(S)\equiv \{v = \sum_{m=1}^{\infty} f_m,~~~~f_m \in \mathcal{D}_j^+ \;\; \forall~m \geq 1,~~~~~\sum_{m=1}^{\infty} I_jf_m(S) < \infty\}.
\end{equation}
\end{definition}

Therefore  $\mathcal{M}_{0} = \mathcal{M}_0(S)$,  valid  for all $S$.

Observe that $\mathcal{D}_j^+\subset \cap_{S\in \Se}\mathcal{M}_j(S)\subset \Me_j$. In particular $v\in \Me_j$ iff there exist $F\subset \Se$ with null complement such that $v\in \cap_{S\in F}\mathcal{M}_j(S)$.


\begin{proposition}  \label{firstComputationOfCI}
Let $v = \sum_{m=1}^{\infty} f_m,~~~~f_m \in \mathcal{D}_j^+$, then:
\begin{equation} \label{oneInclusion}
 \underline{\sigma}_j v= \overline{\sigma}_j v =  \overline{I}_j v = \sum_{m \geq 1} I_j f_m.
\end{equation}
Therefore, if $v \in \mathcal{M}_j$ (i.e. $\sum_{m \geq 1} I_j f_m < \infty, ~a.e.$): $\mathcal{M}_j \subseteq \mathcal{L}_{\overline{\sigma}_j}^+$ and
$\int_j v = \sum_m \int_j f_m$.

Moreover,
given  $v \in \mathcal{M}_j(S)$, for any $\epsilon>0$, there exist  $h\in \mathcal{D}_j^+,~u\in\mathcal{M}_j(S)$, both depending on $S$ and $\epsilon$, satisfying: $\overline{I}_ju(S) \leq \epsilon$ and $v=h+u$.
\end{proposition}

\begin{proof}
We argue as follows, for fixed $L\ge 0$
\begin{equation}  \label{chainOfInequalities}
\sum_{m = 1}^L \underline{\sigma}_j f_m \leq \underline{\sigma}_j (\sum_{m= 1}^L f_m) \leq \underline{\sigma}_j v \leq \overline{\sigma}_j v \leq \overline{I}_j v \leq \sum_{m \geq 1} \overline{I}_j f_m \leq \sum_{m \geq 1} I_j f_m
\end{equation}
where the first inequality holds by superadditivity of $\underline{\sigma}_j$, the second one by isotony, the third one from Corollary \ref{good_properties} item $(1)$ and the rest by previous reasoning.

From (\ref{chainOfInequalities}) we can conclude (\ref{oneInclusion}): item  $(4)$ in Proposition \ref{Properties_L} gives $\sum_{m = 1}^L I_j f_m =\sum_{m = 1}^L \underline{\sigma}_j f_m$, then (\ref{oneInclusion}) follows by taking $ L \rightarrow \infty$.

From (\ref{oneInclusion}), if $v \in \mathcal{M}_j$ then $0 \leq \underline{\sigma}_j v= \overline{\sigma}_j v<\infty~~a.e.$, thus $\overline{\sigma}_j v-\underline{\sigma}_j v=0~~a.e.$ and being no-negative $v\in \mathcal{L}_{\overline{\sigma}_j}^+$, consequently last equality in (\ref{oneInclusion}) and Proposition \ref{Properties_L} item $(4)$ gives $\int_j v = \sum_m \int_j f_m$.


Finally, if  $v \in \mathcal{M}_j(S)$ then $v = \sum\limits_{m=1}^{\infty} f_m$ with $f_m \in \mathcal{D}_j^+$ for $m \geq 1$ and $\sum\limits_{m=1}^{\infty} I_jf_m(S) < \infty$. Thus for a given $\epsilon >0$, there exist $m_{\nu}$ such that $\sum\limits_{m>m_v}^{\infty}I_jf_m(S) \leq \epsilon$, then  $v=h+u$, where
$h = \sum\limits_{m=1}^{m_v} f_m\in \mathcal{D}_j^+$ and $u=\sum\limits_{m>m_v}^{\infty}f_m$ ~~ satisfies ~~ $\overline{I}_ju(S)=\sum\limits_{m>m_v}^{\infty}I_jf_m(S)\le \epsilon$, hence $u \in \mathcal{M}_j(S)$.
\end{proof}

\begin{remark} \label{movingMass}
\noindent If $f\in \Ee_j + \Me_j(S) -\Me_j(S)$ then, from the last result of Proposition \ref{firstComputationOfCI}, for any $\epsilon > 0$, $f$ can be written as $f=h+v-u$  with $h \in \Ee_{(S, j)}, u, v \in  \Me_j(S)$, and  $\overline{I}_j u(S),~\overline{I}_j v(S) \le \epsilon$. See Theorem \ref{punctualCharacterization}.
\end{remark}

\begin{corollary}\label{linearityE+M-M}
For $u, v \in \mathcal{M}_j(S), h \in \mathcal{E}_j$:
\begin{enumerate}
\item $\overline{I}_j (u + \alpha v)(S) = \overline{\sigma}_j (u + \alpha v)(S) =  \overline{\sigma}_j u(S) + \alpha \overline{\sigma}_j v(S)= \overline{I}_j u(S) + \alpha \overline{I}_j v(S),~$ whenever $\alpha \geq 0$.
\item $\underline{\sigma}_j(h+ v-u)(S) = \overline{\sigma}_j(h+v-u)(S) = I_j h(S) + \overline{I}_j v(S) - \overline{I}_j u(S)$.
\end{enumerate}
Therefore $\Ee_j + \Me_j -\Me_j \subset \mathcal{L}_{\overline{\sigma}_j}$ and
\begin{equation} \label{linearEverywhere}
\int_j (h+v-u)\dot{=} \int_j h + \int_j v- \int_j u.
\end{equation}
\end{corollary}
\begin{proof}
(\ref{oneInclusion}) in Proposition \ref{firstComputationOfCI} implies the first and last equalities in $(1)$ given that $u, ~v,\\ ~u + \alpha v\in\Me_j(S)$. The second equality in $(1)$ follows from Lemma \ref{linearityOfSigma} by taking $f \equiv u$, $g \equiv v$;
these results are applicable given that by Proposition \ref{firstComputationOfCI} and hypothesis $0 \leq \underline{\sigma}_j u(S)=\overline{\sigma}_j u(S)<\infty$ so $\overline{\sigma}_j u(S)-\underline{\sigma}_j u(S)=0$  (similarly for $v$).

To derive $(2)$, by Proposition \ref{Properties_L} item $(4)$ $\overline{\sigma}_j h(S)=\underline{\sigma}_j h(S)$, and $(\ref{oneInclusion})$ in Proposition \ref{firstComputationOfCI}\; $\overline{\sigma}_j u(S)=\underline{\sigma}_j u(S)$ and $\overline{\sigma}_j v(S)= \underline{\sigma}_j v(S)$. Since in the three cases the involved values are finite, $h,v,u$ satisfy the hypothesis of Lemma \ref{linearityOfSigma}, from where
$\overline{\sigma}_j (h+v)(S)=\underline{\sigma}_j (h+v)(S)$, so for the same reason $h+v$ satisfy the referred hypothesis. Consequently
\[
\overline{\sigma}_j (h+v-u)(S)=\overline{\sigma}_j (h+v)(S) - \overline{\sigma}_j u(S)=\overline{\sigma}_j h(S)+ \overline{\sigma}_j v(S) - \overline{\sigma}_j u(S).
\]
From where, by Proposition \ref{Properties_L} item $(4)$ and $(\ref{oneInclusion})$ in Proposition \ref{firstComputationOfCI} again, it holds
\[
\overline{\sigma}_j (h+v-u)(S)=I_j h(S)+ \overline{I}_j v(S) - \overline{I}_j u(S)=\underline{\sigma}_j h(S)+ \underline{\sigma}_j v(S) - \underline{\sigma}_j u(S).
\]
And another use of Lemma \ref{linearityOfSigma} gives the first equality in $(2)$.

$\Ee_j+ \Me_j- \Me_j \subseteq \mathcal{L}_{\overline{\sigma}_j}$ as well
as (\ref{linearEverywhere}) follow directly from $(2)$ just proven and finiteness of the involved values.
\end{proof}



\vspace{.1in}
Clearly, $\Ee_j + \Me_j -\Me_j=\Me_j -\Me_j$, using $\Ee_j + \Me_j -\Me_j$ above was meant to emphasize that $f\in \Me_j -\Me_j$  can be represented with a term $h \in \Ee_j$ with its conditional integral $\int_j f(S)$ concentrated in $\int_j h(S)$ as indicated in Remark \ref{movingMass}.

\vspace{.1in}
A hypothesis in the next theorem contains the qualifier {\it a.e.}, even if it were strengthened to every $S$, the conclusion will still hold only a.e.

\begin{theorem}[K\"{o}nig's Satz 2.9]\label{Konig2.9}
Let $f\in Q$, if for $n\ge 1$ there exist $f_n\in \mathcal{L}_{\overline{\sigma}_j}$ such that
$\lim\limits_{n\rightarrow \infty}\|f-f_n\|_j=0~~~~a.e.$, then $f\in \mathcal{L}_{\overline{\sigma}_j}$ and
$\lim\limits_{n\rightarrow \infty}\int_jf_n\doteq\int_jf$.
\end{theorem}
\begin{proof} The following holds a.e. $\lim\limits_{n\rightarrow \infty}\overline{\sigma}_j(|f-f_n|)\le\lim\limits_{n\rightarrow \infty}\|f-f_n\|_j=0$; by Proposition \ref{konig} item $(6)$ it follows that, $\overline{\sigma}_jf\le \overline{\sigma}_jf_n + \overline{\sigma}_j(|f-f_n|)$ and $\overline{\sigma}_j(-f)\le \overline{\sigma}_j(-f_n) + \overline{\sigma}_j(|-f+f_n|)$. Thus, since $f_n\in \mathcal{L}_{\overline{\sigma}_j}$ and $\overline{\sigma}_j(-f_n)=-\underline{\sigma}_j f_n$
\[
\overline{\sigma}_jf-\underline{\sigma}_jf=\overline{\sigma}_jf+\overline{\sigma}_j(-f)\le 2\overline{\sigma}_j(|f-f_n|)~~a.e.,
\]
from which we obtain $f\in \mathcal{L}_{\overline{\sigma}_j}$. Now from Proposition \ref{Konig2.1} item $(4)$ and Corollary \ref{good_properties} item $(2)$,
\[
|\int_j f - \int_j f_n| \doteq |\overline{\sigma}_j(f - f_n)|\le \overline{\sigma}_j(|f-f_n|)\xrightarrow[n\rightarrow\infty]{{}}0~~a.e.,
\]
which means that $\lim\limits_{n\rightarrow \infty}\int_jf_n\doteq\int_jf$.
\end{proof}

\subsection{Characterization of Conditional Integrability}
The remaining of the section introduces the key sets $\mathcal{L}_j$ which we prove
satisfy $\mathcal{L}_j \subseteq \mathcal{L}^+_{\overline{\sigma}_j}$ and so $\mathcal{L}_j - \mathcal{L}_j \subseteq \mathcal{L}_{\overline{\sigma}_j}$. Under an strengthening of property $(L_j)$ (by means of a
condition introduced in \cite{konig} and which we will label property $(K_j)$)
we also show that $\mathcal{L}_j = \mathcal{L}^+_{\overline{\sigma}_j}$ .
We do pursue some  results involving $\mathcal{L}_j$ and, in that way gain some generality (given that the said results will hold independently of $(K_j)$).

We want to highlight (\ref{weakerVersionOfKj}) and (\ref{functionalsOnLj}) below, they give the
equality $\overline{\sigma}_j = \overline{I}_j$ over a set of functions that later will be characterized as being
the nonnegative integrable functions. That being said, our proof of that characterization relies on having available $\overline{\sigma}_j = \overline{I}_j$  over the larger set $P$. Establishing this latter equality will be shown,  in Proposition \ref{keyForConvergenceTheoremsAndCharacterization},  to be equivalent to property $(K_j)$.

\begin{lemma}  \label{fundamentalIdentities}
Let $f: \mathcal{S} \rightarrow [0, \infty]$, fix $j \geq 0$ and $S \in \mathcal{S}$. Assume that  for any $ \epsilon >0$, ~~
$\exists$ $u, v \in \mathcal{M}_j(S)$ satisfying
$f = v-u$ a.e. on $\mathcal{S}_{(S,j)}$ and $\overline{I}_j u(S) \leq \epsilon$.
Then:
\begin{equation} \label{weakerVersionOfKj}
\underline{\sigma}_j f(S) = \overline{\sigma}_j f(S)= \overline{I}_j f(S)=
\overline{I}_j v(S)- \overline{I}_j u(S).
\end{equation}
\end{lemma}
\begin{proof} Let $f, u, v$ and $\epsilon$ as in the statement; given that $f = v-u ~~a.e.$ on $\mathcal{S}_{(S,j)}$, setting $A=\{\tilde{S}\in\mathcal{S}_{(S,j)}:f(\tilde{S})\neq v(\tilde{S})-u(\tilde{S})\}$ we have $\overline{I}_j({\bf 1}_A)(S)=0$. Notice that
$f(\tilde{S})= [\mathbf{1}_{A^c}(v-u)](\tilde{S})+ [\mathbf{1}_{A}f](\tilde{S})$ for all $\tilde{S} \in \mathcal{S}_{(S,j)}$. Since $[\mathbf{1}_{A^c}(v-u)]\le v$, and \\ $\mathbf{1}_{A}f\le \sum\limits_{m\ge 1}\mathbf{1}_{A}$, so $\mathbf{1}_{A}f$ is a conditional null function at $(S,j)$, it follows that
\begin{equation} \label{toBeUsed}
\overline{I}_j f(S) \leq \overline{I}_j[\mathbf{1}_{A^c}(v-u)](S) + \overline{I}_j[\mathbf{1}_{A}f](S)\le \overline{I}_j v(S).
\end{equation}
Given $f_m\in \mathcal{D}_{(S,j)}^+$ for $m\ge 0$, satisfying $f\leq -f_0+\sum\limits_{m\ge 1}f_m$ on $\mathcal{S}_{(S,j)}$ we can then write
$\mathbf{1}_{A^c}v-\mathbf{1}_{A^c}u \le \mathbf{1}_{A^c}(v-u)+ \mathbf{1}_{A}f\le -f_0+\sum\limits_{m\ge 0}f_m$ on $\mathcal{S}_{(S,j)}$, which gives:
\begin{equation} \label{toBeUsedShortly}
\mathbf{1}_{A^c}v+ f_0 \le  \sum\limits_{m\ge 1}f_m+\mathbf{1}_{A^c}u~~\mbox{on}~~\mathcal{S}_{(S,j)}.
\end{equation}
Furthermore, notice that  $\mathbf{1}_{A^c}v +f_0=v+f_0~~a.e.$ on $\mathcal{S}_{(S,j)}$ and $\mathbf{1}_{A^c}u\le u$. Therefore  we obtain from (\ref{toBeUsedShortly})
(where we rely on (\ref{oneInclusion}) from Proposition \ref{firstComputationOfCI}, Proposition \ref{leinertTheorem} item $(3)$  and properties of $\overline{I}_j$):
\[
\overline{I}_jv(S)+\overline{I}_jf_0(S) = \overline{I}_j[v +f_0](S)= \overline{I}_j[\mathbf{1}_{A^c}v +f_0](S)\leq \sum\limits_{m\ge 1}I_jf_m(S)+ \overline{I}_j u(S).
\]

Combining the above inequality with (\ref{toBeUsed})  we obtain:
\[
\overline{I}_jf(S) \leq \overline{I}_jv(S)\le -I_jf_0(S)+ \sum\limits_{m\ge 1}I_jf_m(S)+ \epsilon.
\]
It then follows that $\overline{I}_j f(S) \leq \overline{\sigma}_j f(S) + \epsilon$ which implies $\overline{I}_jf(S) \leq \overline{\sigma}_jf(S)$ but the reverse inequality holds for any $f \geq 0$, and any $S$, from Proposition \ref{konig} item $(\ref{sigma_le_Iup})$ and so,  $\overline{I}_jf(S) = \overline{\sigma}_jf (S)$ is established.

From Proposition \ref{konig} item $(4)$, $\overline{\sigma}_j f(S)= \overline{\sigma}_j (v-u)(S)$ and $\underline{\sigma}_j f(S)=\underline{\sigma}_j (v-u)(S)$ holds. Lemma \ref{linearityOfSigma}, which is applicable by an appeal to Proposition \ref{firstComputationOfCI} and  our assumption $u, v \in \mathcal{M}_j(S)$,  gives $\underline{\sigma}_j(v-u)(S)=   \overline{\sigma}_j(v-u)(S)=
\overline{\sigma}_jv(S)+  \overline{\sigma}_j(-u)(S) =  \overline{I}_j v(S)- \overline{I}_j u(S)$. Therefore (\ref{weakerVersionOfKj}) is established.
\end{proof}



\begin{definition} \label{postiveIntegrablesWithoutK}
\begin{equation} \nonumber
\mathcal{L}_j \equiv \{f:\mathcal{S} \rightarrow [0, \infty]: ~\mbox{the following holds a.e. in the variable}~S \in \mathcal{S}:\\
\end{equation}
\begin{equation} \nonumber
~[ ~\mbox{for all}~~ \epsilon >0~\exists~ u, v \in \mathcal{M}_j(S)~~\mbox{such that}~
~f = v-u~~a.e. ~\mbox{on}~~\mathcal{S}_{(S,j)}~~\mbox{and}~~ \overline{I}_j u(S) \leq \epsilon]\}.
\end{equation}
\end{definition}
Notice that $\mathcal{L}_j$ is a positive cone that contains the nonnegative null functions. To check the latter claim let $g \in P$ be a null function,
it follows from item $(5)$ of Proposition \ref{leinertTheorem} that
$\overline{I}_j(|g|)=0$ a.e. i.e. $[g=0 ~a.e. ~\mbox{on}~\mathcal{S}_{(S,j)}]$ holds a.e. in the variable $S$. This shows $g \in \mathcal{L}_j$.

\vspace{.1in}
\noindent
Lemma \ref{fundamentalIdentities} then gives.
\begin{corollary} \label{restrictedPositive}
For each $f \in \mathcal{L}_j$, the following statement holds a.e.:
\begin{equation}  \label{functionalsOnLj}
\underline{\sigma}_j f=\overline{\sigma}_j f= \overline{I}_jf~~~\mbox{and so} ~~~ \mathcal{L}_j \subseteq \mathcal{L}^+_{\overline{\sigma}_j}.
\end{equation}
\end{corollary}
\begin{proof} The equalities are clear, and since by (\ref{weakerVersionOfKj}), $0 \leq \overline{\sigma}_j f= \underline{\sigma}_j f   = \overline{I}_j v -\overline{I}_j u < \infty ~~~a.e.$, it follows that $\overline{\sigma}_j f-\underline{\sigma}_j f=0$ $a.e.$
\end{proof}

\begin{theorem}[Beppo-Levi for $\mathcal{L}_j$]
Let $f_k \in \mathcal{L}_j$, then:
\begin{equation} \nonumber
\underline{\sigma}_j(\sum_k f_k) = \overline{\sigma}_j(\sum_k f_k)= \overline{I}_j(\sum_k f_k) = \sum_k \overline{I}_j f_k \dot{=} \sum_k \overline{\sigma}_j f_k \dot{=} \sum_k \underline{\sigma}_j f_k.
\end{equation}
Moreover, if $\sum_k \overline{I}_j f_k < \infty$~a.e., then
$\sum_k f_k \in \mathcal{L}_j$ and $\int_j f \dot{=} \sum_k \int_j f_k$.
\end{theorem}
\begin{proof}
\begin{equation} \label{interchangeOfSumAndIntegralWithPointwiseEquality}
\sum_{k=1}^n \underline{\sigma}_jf_k \leq \underline{\sigma}_j(\sum_k f_k) \leq \overline{\sigma}_j(\sum_k f_k)\leq \overline{I}_j(\sum_k f_k) \leq \sum_k \overline{I}_j f_k \dot{=} \sum_k \overline{\sigma}_j f_k \dot{=} \sum_k \underline{\sigma}_j f_k,
\end{equation}
where we used Lemma \ref{fundamentalIdentities} for the last two equalities. We then conclude by taking $n \rightarrow \infty$ (and relying on the fact that a countable union of null sets is also null).

Let $\epsilon_k > 0$ be such that $\sum \epsilon_k \leq \epsilon$. By assumption, for each $k$ the following holds
\begin{equation} \label{toBeAdded}
[f_k = v_k-u_k~\mbox{a.e. on}~\mathcal{S}_{(S,j)}]~\mbox{a.e.~in} ~\mathcal{S},
\end{equation}
 with $\overline{I}_ju_k(S) \leq \epsilon_k$ and $u_k, v_k \in \mathcal{M}_j(S)$.
Given that a countable union of null sets is a null set we may assume
the set of $S \in \mathcal{S}$ for which
(\ref{toBeAdded}) does not hold {\it for all} $k$ is a null set.

 \vspace{.1in}
 Let $\,S \in \mathcal{S}$ be such that (\ref{toBeAdded}) holds for all
 $k$, and $\,u \equiv \sum_k u_k$,\, given that $u_k = \sum_m f_{m, k},$ $f_{m, k} \in \mathcal{E}_j^+$.   From Proposition \ref{firstComputationOfCI},
\[
\overline{I}_j u(S) = \sum_{m,k} I_j f_{m,k}(S)= \sum_k (\sum_m I_j f_{m,k})(S) =
\sum_k \overline{I}_j u_k(S) \leq \epsilon.
\]
Therefore $u < \infty$ a.e.
on $\mathcal{S}_{(S,j)}$, this implies that $f = v-u$ holds a.e.
on $\mathcal{S}_{(S,j)}$,  where $f \equiv \sum_k f_k, ~v \equiv \sum_k v_k$. Notice that $u \in \mathcal{M}_j(S)$.

We have $\overline{I}_j v \leq \sum_k \overline{I}_j v_k$, also,
$\sum_{k=1}^n \overline{I}_j v_k = \sum_{k=1}^n \underline{\sigma}_j v_k
\leq \overline{\sigma}_j \sum_{k=1}^n v_k \leq \overline{I}_j v$,  where we have relied on (\ref{oneInclusion}). From $v_k = f_k+u_k$ being valid a.e. on $\mathcal{S}_{(S,j)}$ we obtain $\overline{I}_j v_k(S) \leq \overline{I}_j f_k(S)+ \overline{I}_j u_k(S)$. For convenience, let us call $S \in \mathcal{S}$ {\it admissible} if (\ref{toBeAdded}) holds for all
 $k$ (as we have been considered thus far) and $\sum_k \overline{I}_j f_k(S) < \infty$ as well. Therefore the set of non admissible $S$ is a null set; if $S$ is admissible:
\begin{equation}  \label{tobeUsedSoon}
\overline{I}_j v(S) \leq \sum_k \overline{I}_j v_k(S) \leq \sum_k \overline{I}_j f_k(S)+\sum_k \overline{I}_j u_k(S) < \infty.
\end{equation}
Given that $v_k = \sum_m g_{m, k},~~ g_{m, k} \in \mathcal{E}_j^+$,
$v = \sum_{m,k} g_{m,k}$ from
and (\ref{tobeUsedSoon}) we conclude that $v \in \mathcal{M}_j(S)$.
We then have proved that for admissible $S$:
\begin{equation} \label{holdsAe}
f= v-u~ \mbox{a.e. on}~ \mathcal{S}_{(S,j)}~\mbox{with}~ u, v \in \mathcal{M}_j(S),\overline{I}_ju \leq \epsilon.
\end{equation}
From our selection of $S$ we conclude that (\ref{holdsAe}) holds a.e.
on $\mathcal{S}$. Therefore $f \in \mathcal{L}_j$ and
$\int_j f \dot{=} \sum_k \int_j f_k$ follows from (\ref{interchangeOfSumAndIntegralWithPointwiseEquality}).
\end{proof}

From $\mathcal{L}_j \subseteq \mathcal{L}_{\overline{\sigma}_j}^+$, a fact established above, it follows that $\mathcal{L}_j - \mathcal{L}_j \subseteq \mathcal{L}_{\overline{\sigma}_j}$ given that the latter is a vector space. As we pointed out, $\mathcal{L}_j$ is a positive cone so
$\alpha \mathcal{L}_j + \beta \mathcal{L}_j \subseteq \mathcal{L}_{\overline{\sigma}_j}$, $\alpha, \beta \in \mathbb{R}$ follows.

\begin{definition}
\begin{equation} \nonumber
\mathcal{L}^K_j \equiv \{f:\mathcal{S} \rightarrow [- \infty, \infty]:~ ~\mbox{the following holds a.e. in the variable}~  S \in \mathcal{S}:
\end{equation}
\begin{equation} \nonumber
[\mbox{for all}~ \epsilon >0~\exists~ u, v \in \mathcal{M}_j(S),~~h \in \mathcal{E}_{(S,j)}, f = (h+v-u)~a.e. ~~\mbox{on}~~\mathcal{S}_{(S,j)}~~\mbox{and} ~~\overline{I}_j u(S) \leq \epsilon, ~\overline{I}_j v(S) \leq \epsilon]\}.
\end{equation}
\end{definition}
Define also
\begin{equation} \nonumber
\tilde{\mathcal{L}}^K_j \equiv \{f:\mathcal{S} \rightarrow [- \infty, \infty]:~
~ ~\mbox{the following holds a.e. in the variable}~  S \in \mathcal{S}:
\end{equation}
\begin{equation} \nonumber
[\mbox{for all}~ \epsilon >0~~\exists~ u, v \in \mathcal{M}_{j}(S),~~h \in \mathcal{E}_j,~
f = (h+v-u)~a.e. ~~\mbox{on}~~\mathcal{S}_{(S,j)}~~\mbox{and}
\end{equation}
\begin{equation} \nonumber
~~ -\epsilon \leq \underline{\sigma}_j(f-h)(S) \leq \overline{\sigma}_j(f-h)(S) \leq \epsilon]\}.
\end{equation}

\begin{proposition}[See proof in Appendix \ref{someProof1}]  \label{goodToKnow}
$\mathcal{L}^K_j = \tilde{\mathcal{L}}^K_j$.
\end{proposition}

The following result will give us the tools to extend functions
defined by local conditions (i.e. in terms of the conditional spaces $\mathcal{S}_{(S,j)}$)
to functions defined directly on $\mathcal{S}$ (i.e. {\it globally} defined).

\vspace{.1in}
A subset $\mathcal{C} \subseteq \mathcal{S}$ will be called \emph{admissible} if it satisfies: whenever $S^1,S^2\in\mathcal{C}$ then $\Se_{(S^1,j)}\cap\Se_{(S^2,j)}=\emptyset$. For a fixed $j \geq 0$ and a given admissible subset $\mathcal{C} \subseteq
\mathcal{S}$ we will use the notation $w(\cdot)$ for a choice function of the following type:

 $w: \widetilde{\mathcal{C}}\equiv\cup_{S\in \mathcal{C}} \Se_{(S,j)}\rightarrow \mathcal{C}$ such that if we let  $S \equiv w(\tilde{S})$ we have  $\tilde{S} \in \mathcal{S}_{(S, j)}$. \\
Notice that $\Se_{(w(\tilde{S}),j)}= \Se_{(\tilde{S},j)}$ for any $\tilde{S} \in \widetilde{\mathcal{C}}$, and
that $w(\widetilde{\mathcal{C}})=\mathcal{C}$.

\begin{lemma} \label{globalizing}
Fix $j\ge 0$ and let $\mathcal{C}$ be admissible as defined above. Let $(h^{S})_{S \in \mathcal{C}}$ and $(v^{S})_{S \in \mathcal{C}}$ be families of functions where $h^{S} \in \Ee_{(S,j)}$, $v^{S} \in \Me_j(S)$ for $S \in \mathcal{C}$ and set $\widetilde{\mathcal{C}}\equiv \cup_{S\in \mathcal{C}}\Se_{(S,j)}$ . Then, there exist $h\in \Ee_j$, $v\in  \cap_{S \in \widetilde{\mathcal{C}}}\Me_j(S)$ such that
the following holds for all $S \in \widetilde{\mathcal{C}}$:
\begin{equation}\label{global_E_j}
h=h^{w(S)}\;\; on \; \Se_{(S,j)} ~~\quad \mbox{and}\quad I_jh(S)= I_jh^{w(S)}(S),
\end{equation}
\begin{equation}\label{global_M_j}
v=v^{w(S)}\;\; on \; \Se_{(S,j)} ~~  \quad \mbox{and}\quad \overline{I}_jv(S)= \overline{I}_jv^{w(S)}(S).
\end{equation}
\end{lemma}
\begin{proof}
We will prove (\ref{global_M_j}), the proof of (\ref{global_E_j}) is similar.
For each $S\in\mathcal{C}$, $v^S = \sum\limits_{m=1}^{\infty} f^S_m,~~~~f^S_m \in \mathcal{D}_j^+$, so by (\ref{conditionalElementarySpace})\;  ${f^S_m}|_{\Se_{(S,j)}}\!\!\in \mathcal{D}_{(S,j)}^+$. Keeping in mind that
$\mathcal{C} \subseteq \widetilde{\mathcal{C}}$,
define for $m \geq 1$:
\[
f_m(S)=\left\{\begin{array}{ccc}
f^{w(S)}_m(S) & if & S\in \widetilde{\mathcal{C}} \\
0 & otherwise
\end{array}\right.
\quad \mbox{and} \quad v \equiv\sum_{m=1}^{\infty} f_m.
\]
Then, for $S \in \widetilde{\mathcal{C}}$,
${f_m}|_{\Se_{(S,j)}}={f^{w(S)}_m}|_{\Se_{(S,j)}}\in \mathcal{D}_{(S,j)}^+$ and so $v = v^{w(S)}$ on $\mathcal{S}_{(S, j)}$. 
On the other hand, if $S \notin \widetilde{\mathcal{C}}$ we have $f_m|_{\mathcal{S}_{(S, j)}}=0 \in \mathcal{D}_{(S,j)}^+$, therefore,
$f_m \in \mathcal{D}_j^+$. 

Moreover, if $S \in \widetilde{\mathcal{C}}$, from Proposition \ref{firstComputationOfCI}
\[
\overline{I}_jv(S) = \sum_{m=1}^{\infty} I_jf_m(S) = \sum_{m=1}^{\infty} I_jf^{w(S)}_m(S) = \overline{I}_jv^{w(S)}(S)= \overline{I}_jv^{w(S)}(w(S))<\infty,
\]
where we used Corollary \ref{representationUniqueness} for the second equality and the fact that $\overline{I}v^{w(S)}$ is constant on
$\mathcal{S}_{(w(S), j)}$, and $S \in \mathcal{S}_{(w(S), j)}$ for the last equality. The last inequality follows from the assumption
$v^{S} \in \mathcal{M}_j(S)$ for any $S \in \mathcal{C}$. We then conclude that $v \in \mathcal{M}_j(S)$
for each $S \in \widetilde{\mathcal{C}}$.
\end{proof}

The next proposition establishes global formulations
for the sets $\mathcal{L}_j$
and $\mathcal{L}_j^K$. The result will be used in Theorem \ref{punctualCharacterization} to obtain that the global formulation applies to integrable functions as well.
\begin{proposition} \label{globalFormuylations}
\begin{equation} \label{globalFormulationLj}
\mathcal{L}_j = \mathcal{L}_j^G \equiv
\{f:\mathcal{S} \rightarrow [0, \infty]:
~\mbox{for each}~~ \epsilon >0~\exists~ u, v \in \mathcal{M}_j~
~\mbox{such that the following}
\end{equation}
\begin{equation} \nonumber
~\mbox{holds a.e. in $S$:}~~[f = v-u~~a.e. ~\mbox{on}~~\mathcal{S}_{(S,j)}~~\mbox{and}~~\overline{I}_j u(S) \leq \epsilon, ~\overline{I}_j v(S) < \infty\}.
\end{equation}
\begin{equation} \label{globalFormulationL_j^K}
\mathcal{L}_j^K = \mathcal{L}_j^{K, G} \equiv
\{f:\mathcal{S} \rightarrow [- \infty, \infty]:
~\mbox{for each}~~ \epsilon >0~\exists~ u, v, \in \mathcal{M}_j,~h \in \mathcal{E}_j~~\mbox{such that}
\end{equation}
\begin{equation} \nonumber
\mbox{ the following holds a.e. in $S$}:~[f = h+v-u~~a.e. ~\mbox{on}~~\mathcal{S}_{(S,j)}~~\mbox{and}~~\overline{I}_j u(S), \overline{I}_j v(S) \leq \epsilon]\}.
\end{equation}
\end{proposition}
\begin{proof}
We prove (\ref{globalFormulationLj}), the proof of (\ref{globalFormulationL_j^K}) is analogous.

Let $f \in \mathcal{L}_j^G$ and $\epsilon >0$, then there exist $u, v \in \mathcal{M}_j$ such that for $a.e.~~S\in\Se$ it holds that $f= v- v~~~ a.e. ~\mbox{on}~~ \mathcal{S}_{(S,j)}$ and $\overline{I}u(S) \leq \epsilon$, $\overline{I}v(S) < \infty$. Consequently,  for those $S$, $u, v \in \mathcal{M}_j(S)$. Thus, for the referred $S\in\Se$, $u, v$ (the same for any $S$) verify the conditions which gives that $f \in \mathcal{L}_j$.

Let now $f \in \mathcal{L}_j$ then, there exists $\hat{\mathcal{S}}\subset\Se$, with $\hat{\mathcal{S}}^C$ a null set, such that for each $S \in \hat{\mathcal{S}}$ we have: for any $\epsilon>0$  there exist
$u^S, v^S\in \mathcal{M}_j(S)$ satisfying
\[
f=v^S-u^S \;\; a.e. \;\;\mbox{on} \; \Se_{(S,j)} \;\; \mbox{and} \;\; \overline{I}_ju^S(S) < \epsilon.
\]
We can construct  $\mathcal{C} \subset \hat{\mathcal{S}}$ that satisfies  $\tilde{\mathcal{C}} \equiv \cup_{S \in \mathcal{C}}\Se_{(S,j)}=\cup_{S\in \hat{\mathcal{S}}}\Se_{(S,j)}$  \, and  $\Se_{(S^1,j)}\cap\Se_{(S^2,j)}=\emptyset$ for $S^1,S^2\in \mathcal{C}$. Thus, the families of functions $(h^{S})_{S \in \mathcal{C}}$,\,
$(u^{S})_{S \in \mathcal{C}}$\, and \,$(v^{S})_{S \in \mathcal{C}}$ satisfy the corresponding hypothesis of Lemma \ref{globalizing}, which gives  the existence of  $u,v \in \cap_{S \in \tilde{\mathcal{C}}} \mathcal{M}_j(S)$ such that for any $S \in \tilde{\mathcal{C}}$ (where $w(S)$ below is as introduced in  Lemma
\ref{globalizing})
\begin{equation}\nonumber
f= v^{w(S)}- u^{w(S)}= v- u~~\mbox{a.e. on}~\mathcal{S}_{(S,j)}~~\mbox{and}~~
\end{equation}
\begin{equation}\nonumber
\overline{I}u(S) = \overline{I}u^{w(S)}(S) \leq \epsilon, ~
\overline{I}v(S) = \overline{I}v^{w(S)}(S) < \infty.
\end{equation}
Notice that $\hat{\mathcal{S}} \subseteq
\cup_{S \in \hat{\mathcal{S}}} \mathcal{S}_{(S,j)}$, it then follows that
$\tilde{\mathcal{C}}^C= \cap_{S \in \hat{\mathcal{S}}} \mathcal{S}^C_{(S,j)}$ is a null set; therefore, $u,v \in \mathcal{M}_j$ which allows us to conclude $f \in \mathcal{L}_j^G$.

\end{proof}

Notice that below we establish $\mathcal{L}_j - \mathcal{L}_j  = \mathcal{L}^K_j$ which easily implies
$\alpha \mathcal{L}_j + \beta \mathcal{L}_j  = \mathcal{L}^K_j$, for arbitrary $\alpha, \beta \in \mathbb{R}$.

\begin{proposition}\label{characterizationLjK}
Whenever $f \in ~\mathcal{L}^K_j$~~we have
$\underline{\sigma}_j f \dot{=} \overline{\sigma}_j f$. Moreover,
\begin{equation} \nonumber
 \mathcal{L}_j - \mathcal{L}_j  = \mathcal{L}^K_j \subseteq \mathcal{L}_{\overline{\sigma}_j}.
\end{equation}
\end{proposition}
\begin{proof}
$\underline{\sigma}_j f \dot{=} \overline{\sigma}_j f~~\mbox{for}~~f \in \mathcal{L}^K_j$ follows from Corollary \ref{linearityE+M-M}. In particular, $\mathcal{L}^K_j \subseteq \mathcal{L}_{\overline{\sigma}_j}$
then follows.

To establish $\mathcal{L}_j - \mathcal{L}_j \subseteq \mathcal{L}^K_j$, consider $f=f_1-f_2\in\mathcal{L}_j - \mathcal{L}_j$. Let $S \in \mathcal{S}$ be such that for all $\epsilon >0$ there exist $u, v, u', v'$ in $\mathcal{M}_j(S)$, such that: $f_1 = (v-u),~f_2=(v'-u')$ a.e. on $\mathcal{S}_{(S,j)}$ with $\overline{I}_j u(S) \leq \epsilon/2, \overline{I}_j u'(S) \leq \epsilon/2$. Given that $\overline{I}_j v(S) < \infty, \overline{I}_j v'(S) < \infty$ as well, we can assume (by relying on $(2)$ from Proposition \ref{leinertTheorem}) that $u, v, u', v'$ are finite on the set where the equalities $f_1 = (v-u), ~f_2=(v'-u')$ take place  a.e. on $\mathcal{S}_{(S,j)}$.


From last assertion of Proposition \ref{firstComputationOfCI} we can find $h,h'\in \mathcal{D}_j^+$, $v_1,v'_1\in \mathcal{M}_j(S)$ such that $v=h+v_1, v'=h'+v'_1$ and $\overline{I}_jv_1(S), \overline{I}_jv'_1(S) \le \epsilon/2$, therefore $f = \tilde{v} - \tilde{u} +h$ a.e. on $\mathcal{S}_{(S,j)}$ where  $\tilde{v} \equiv v_1 + u'$ , $\tilde{u} \equiv  v'_1+ u$ and $h \equiv h_1 - h_2\in \mathcal{E}_j$. Given that the above properties hold a.e. in the originally chosen $S \in \mathcal{S}$ we have established $f \in \mathcal{L}^K_j$.

Let now $f \in \mathcal{L}_j^K$ hence, $f \in \mathcal{L}_j^{K, G}$
by Proposition \ref{globalFormuylations}. Therefore there exist $\mathcal{A} \subseteq \mathcal{S}$ with $\overline{I}({\bf 1}_{\mathcal{A}^c})=0$ and such that for each $\epsilon >0$ there exist $u, v \in \mathcal{M}_j$ and $h \in \mathcal{E}_j$ satisfying
\begin{equation} \nonumber
f = h + v -u~~a.e.~~ \mbox{on}~~\mathcal{S}_{(S,j)}~~~\mbox{and}~~~ \overline{I}_ju(S), \overline{I}_jv(S) \le \epsilon ~~~\mbox{for all}~~S \in \mathcal{A}.
\end{equation}

Observe that if $S\in \mathcal{A}$ and $\tilde{S}\in \mathcal{S}_{(S,j)}$ then also
\[f = h + v -u~~a.e.~~ \mbox{on}~~\mathcal{S}_{(\tilde{S},j)}~~~\mbox{and}~~~ \overline{I}_ju(\tilde{S}), \overline{I}_jv(\tilde{S}) \le \varepsilon,\]
this property justifies the use of
 $\mathcal{A} \subseteq \tilde{\mathcal{A}} \equiv \cup_{S\in \mathcal{A}}\mathcal{S}_{(S,j)}$, so $\tilde{\mathcal{A}}^c$ is a null set as well. Below, we will consider an arbitrary choice function $w:\tilde{\mathcal{A}} \rightarrow \mathcal{A}$ satisfying $ \tilde{S} \in \mathcal{S}_{(w(\tilde{S}), j)}$.

Then for any $\tilde{S} \in \tilde{\mathcal{A}}$ and letting $S \equiv w(\tilde{S})$ (a notation that remains in effect for the remaining of the proof), there exists $B^{S}\subset \mathcal{S}_{(S,j)}$ such that $B^{S}_c \equiv (B^{S})^c \cap \mathcal{S}_{(S,j)}$ is a conditionally null set at $(S,j)$, and
\[
f = (h + v -u){\bf 1}_{B^{S}}+ f {\bf 1}_{B^{S}_c}~~ \mbox{on}~~\mathcal{S}_{(S,j)}=\mathcal{S}_{(\tilde{S},j)}.
\]
We remark, in passing, that we can choose the same set $B^S$ for all $\tilde{S} \in \tilde{\mathcal{A}}$ such that $w(\tilde{S})=S$.

Since $h =h_1-h_2, h_1, h_2 \in \mathcal{D}_j^+$, define $\mbox{on}~~\mathcal{S}_{(S, j)}$,
\begin{equation} \nonumber
f^{S}_1 = (h_1-u){\bf 1}_{B^{S}}+ f {\bf 1}_{B^{S}_c}\quad \mbox{and}\quad f^{S}_2=(h_2-v){\bf 1}_{B^{S}}.
\end{equation}
Therefore,  for any $\tilde{S} \in \tilde{\mathcal{A}}$:
\begin{equation} \label{f_1}
f^{S}_1 = h_1-u ~~a.e.~~ \mbox{on}~~\mathcal{S}_{({S},j)},~ h_1, u \in \mathcal{M}_j(S), \overline{I}_ju(S)\le \varepsilon ~~ \mbox{and}
\end{equation}
\begin{equation} \label{f_2}
f^{S}_2 = h_2-v ~~a.e.~~ \mbox{on}~~\mathcal{S}_{(S,j)},~ h_2, v \in \mathcal{M}_j(S), \overline{I}_jv(S)\le \varepsilon.
\end{equation}
We define now
\[
f_1(\tilde{S})=\left\{\begin{array}{ccc}
f_1^{S}(\tilde{S}) & if & \tilde{S} \in \tilde{\mathcal{A}} \\
f(\tilde{S}) & if & \tilde{S} \notin \tilde{\mathcal{A}}
\end{array}\right.,
\;\; \mbox{and}\;\;
f_2(\tilde{S})=\left\{\begin{array}{ccc}
f_2^{S}(\tilde{S}) & if & \tilde{S} \in \tilde{\mathcal{A}} \\
0 & if & \tilde{S} \notin \tilde{\mathcal{A}}.
\end{array}\right.
\]
Then for $ \tilde{S} \in \tilde{\mathcal{A}}$
\[
f_1-f_2=f_1^{S}-f_2^{S}=(h_1-u){\bf 1}_{B^{S}}+ f {\bf 1}_{B^{S}_c} - (h_2-v){\bf 1}_{B^{S}}=
\]\[
=(h+v-u){\bf 1}_{B^{S}}+ f {\bf 1}_{B^{S}_c}=f,
\]
and clearly $f=f_1-f_2$ on $\tilde{\mathcal{A}}^c$. Moreover, by means of  (\ref{f_1}) and (\ref{f_2}) it follows that $f_1,f_2\in \mathcal{L}_j$.

\end{proof}

Under the assumption that property $(K_j)$ holds, we  characterize the spaces of integrable functions in Theorem \ref{punctualCharacterization}. Proposition \ref{proofOfK_j}, in  Appendix \ref{proofOfK}, provides sufficient conditions that imply property $(K_j)$:
\begin{equation} \nonumber
\overline{I}_j(f^+)= I_j(f) + \overline{I}_j(f^-),~\forall f \in \mathcal{E}_j.
\end{equation}
Proposition \ref{keyForConvergenceTheoremsAndCharacterization} below shows that property $(K_j)$ is equivalent to $\overline{\sigma}_j=\overline{I}_j$ on $P$. This latter property is the one used in further results (i.e. $(K_j)$ is not used
directly).

(\ref{KonigCondition}) $\Rightarrow$ (\ref{sigma=I}) in Proposition \ref{keyForConvergenceTheoremsAndCharacterization} below can be proven as in
\cite[Behauptung 1.8]{konig}. Observe also that $\overline{\sigma}_j0=\overline{I}_j0=0$ hence $(L_j)$ holds if $(K_j)$ is valid (a direct proof is given in Proposition \ref{KjImpliesLj} in Appendix \ref{proofOfK}).


\begin{proposition}[See proof in Appendix \ref{someProof1}]  \label{keyForConvergenceTheoremsAndCharacterization}
The following assertions are equivalent:
 \begin{enumerate}
    \item \label{KonigCondition}$(K_j)$.
    \item \label{sigma=I}$\overline{I}_jf=\overline{\sigma}_jf\;\;$ for every $f\in P$.
    \item \label{finiteMat}$\overline{I}_jf^-=\overline{\sigma}_jf^-\;\;$ for every $f\in \mathcal{E}_j$.
 \end{enumerate}
\end{proposition}

The proof of Theorem \ref{punctualCharacterization} below follows the steps of the one in \cite[Satz 2.4]{konig} but taking care of the a.e. condition in Definition \ref{conditional_integrable}.

\begin{theorem} \label{punctualCharacterization}
Assume $\overline{\sigma}_j f = \overline{I}_j f~~\mbox{for all}~~f \geq 0$, then:
\begin{equation} \label{L+Characterization}
\mathcal{L}^+_{\overline{\sigma}_j}= \mathcal{L}_j= \mathcal{L}_j^G.
\end{equation}
\begin{equation} \label{LsigmaCharacterization}
\mathcal{L}_{\overline{\sigma}_j}= \mathcal{L}^K_j = \mathcal{L}_j^{K, G}= \mathcal{L}^+_{\overline{\sigma}_j} - \mathcal{L}^+_{\overline{\sigma}_j}.
\end{equation}
\end{theorem}

\begin{proof}
The second inequalities, in both displays above, were proven in Proposition \ref{globalFormuylations} and are reproduced here for emphasis.
The third equality in (\ref{LsigmaCharacterization}) follows directly from Proposition \ref{characterizationLjK} and the first equality in (\ref{L+Characterization}).
The inclusion
$\mathcal{L}_j \subseteq \mathcal{L}^+_{\overline{\sigma}_j}$
is established in Corollary \ref{restrictedPositive}
and
$\mathcal{L}_j^K \subseteq \mathcal{L}_{\overline{\sigma}_j}$ in
Proposition \ref{characterizationLjK}. All these results only require the validity of $(L_j)$ and we only need the hypothesis $\overline{\sigma}_j f = \overline{I}_j f~~\mbox{for all}~~f \geq 0$ for the remaining two inclusions.

In order to complete the proof of (\ref{L+Characterization}), let $f\in \mathcal{L}^+_{\overline{\sigma}_j}$ and $\epsilon>0$. For $m\ge 1$, choose $\epsilon_m>0$ such that $\epsilon=\sum\limits_{m\ge 1}\epsilon_m$. Define also $\mathcal{N}=\{S\in\Se: \overline{\sigma}_j f(S)-\underline{\sigma}_j(S)\neq 0\}$, which is a null set by definition of $\mathcal{L}_{\overline{\sigma}_j}$.

We show below that for a fixed $S\in \mathcal{N}^c$, there are sequences of functions $(h_m)_{m\ge 1}$ in $\Me_j(S)$ and $(f_m)_{m\ge 1}$ in $\mathcal{L}^+_{\overline{\sigma}_j}$ with $\overline{\sigma}_jf_m(S)<\epsilon_m$, such that $~f_m = h_m-f_{m-1}~$ on $\Se_{(S,j)}$.

Rename $f_0=f$; from hypothesis $0\le\overline{I}_j f_{0}(S) =\overline{\sigma}_j f_{0}(S)<\infty$ (Lemma \ref{linearityOfSigma} applies to our particular $S$ providing $\overline{\sigma}_j f_{0}(S)<\infty$).
By definition of $\overline{I}_j$ there exists $h_1\in \mathcal{M}_j(S)$ with
\begin{equation}\label{firstStep}
f_{0}\le h_1\;\; \mbox{and} \;\;\overline{I}_jh_1(S)<\overline{I}_jf_{0}(S) + \epsilon_1.
\end{equation}
Define $f_1\equiv h_1-f_{0}\ge 0$ on $\Se_{(S,j)}$, and $f_1\equiv 0$ outside $\Se_{(S,j)}$. By Proposition \ref{konig} item $(4)$ we obtain that outside $\Se_{(S,j)}$: $\overline{\sigma}_j f_1=0=\overline{\sigma}_j (-f_1)$ and then $\overline{\sigma}_j f_1-\underline{\sigma}_j f_1=0$. To obtain the same result on $\Se_{(S,j)}$ observe that $\overline{\sigma}_j h_1(S)-\underline{\sigma}_jh_1(S)=0$ follows from the fact that $\overline{\sigma}_j h_1(S)=\overline{I}_jh_1(S)<\infty$  from Proposition \ref{firstComputationOfCI}. Then, applying Lemma \ref{linearityOfSigma} to $f_0$ and $h_1$ we obtain $\overline{\sigma}_j f_1(S)-\underline{\sigma}_j f_1(S)=0$, wich is equivalent to $\overline{\sigma}_j f_1-\underline{\sigma}_j f_1=0$ on $\Se_{(S,j)}$. Noticing that $f_1\ge 0$ we have shown that $f_1\in \mathcal{L}^+_{\overline{\sigma}_j}$. By means of Lemma \ref{linearityOfSigma},
the hypothesis $\overline{\sigma}_j = \overline{I}_j$ on $P$, and (\ref{firstStep}) we derive:
\[
\overline{\sigma}_jf_1(S)=\overline{\sigma}_jh_1(S)-\overline{\sigma}_jf_0(S)<\overline{I}_jf_{0}(S) + \epsilon_1-\overline{I}_jf_0(S)=\epsilon_1.
\]

\vspace{.1in}
For $m\ge 2$ we can then proceeded inductively. Whenever $f_{m-1}\in \mathcal{L}^+_{\overline{\sigma}_j}$  has been constructed satisfying: $\overline{\sigma}_j f_{m-1}(S)-\underline{\sigma}_jf_{m-1}(S)=0$ and $\overline{\sigma}_jf_{m-1}(S)<\epsilon_{m-1}$; it then follows from
$\overline{I}_jf_{m-1}(S)=\overline{\sigma}_jf_{m-1}(S)<\infty$ that  there exist $h_m\in \mathcal{M}_j(S)$ with
\begin{equation}\label{m_Step}
f_{m-1}\le h_m\;\;  \mbox{on}\;\; \Se_{(S,j)},\;\; \mbox{and}\;\; \overline{I}_jh_m(S)<\overline{I}_jf_{m-1}(S) + \epsilon_m.
\end{equation}
We can then define
\begin{equation}\label{f_m}
f_m\equiv h_m-f_{m-1}\;\;  \mbox{on}\;\; \Se_{(S,j)},\;\; \mbox{and}\;\; f_m\equiv 0\;\; \mbox{on}\;\; \Se\setminus \Se_{(S,j)}.
\end{equation}
As we have argued for the case of $f_1$, it follows that $\overline{\sigma}_j f_{m}(S)-\underline{\sigma}_jf_{m}(S)=0$, $f_m\in\mathcal{L}^+_{\overline{\sigma}_j}$ and $\overline{\sigma}_jf_{m}(S)<\epsilon_{m}$. In fact, by Proposition \ref{konig} item $(4)$ we obtain that outside $\Se_{(S,j)}$: $\overline{\sigma}_j f_m=0=\overline{\sigma}_j (-f_m)$ and then $\overline{\sigma}_j f_m-\underline{\sigma}_j f_m=0$. Moreover, $\overline{\sigma}_j h_m(S)-\underline{\sigma}_j h_m(S)=0$
follows from  the fact that $\overline{\sigma}_j h_m(S)=\overline{I}_jh_m(S)<\infty$ from Proposition \ref{firstComputationOfCI}. Then, applying Lemma \ref{linearityOfSigma} to $h_m$ and $f_{m-1}$, gives $\overline{\sigma}_j f_m(S)-\underline{\sigma}_j f_m(S)=0$, wich is equivalent to $\overline{\sigma}_j f_1-\underline{\sigma}_j f_1=0$ on $\Se_{(S,j)}$. Noticing that $f_m\ge 0$, we have then shown: $f_m\in \mathcal{L}^+_{\overline{\sigma}_j}$. Furthermore, by Lemma \ref{linearityOfSigma}, (\ref{m_Step}) and (\ref{f_m}), $\overline{\sigma}_jf_m(S)<\epsilon_m$.

\vspace{.1in} Observing that $\sum\limits_{m\ge 1}f_m\ge 0$ and
\[
\overline{\sigma}_j[\sum\limits_{m\ge 1}f_m](S) = \overline{I}_j[\sum\limits_{m\ge 1}f_m](S)
\le \sum\limits_{m\ge 1}\overline{I}_jf_m(S) <\infty,
\]
then by Proposition \ref{konig} item $(8)$, $\sum\limits_{m\ge 1}f_m<\infty \;\; a.e.~~on~~\Se_{(S,j)}$, from where $\lim\limits_{m\rightarrow \infty}f_m = 0 \;\; a.e.~~on~~\Se_{(S,j)}$.

On the other hand, from (\ref{f_m}), $h_m=f_{m-1}+f_m$ (if $f_{m-1}(\hat{S})=\infty$ for some $\hat{S}\in\Se$, it must be, from (\ref{m_Step}) that $h_m(\hat{S})=\infty$).  Then, using Lemma \ref{linearityOfSigma}, for $m\ge 2$,
\[
\overline{I}_jh_m(S)=\overline{I}_j(f_m+f_{m-1})(S)=
\overline{\sigma}_j (f_m+f_{m-1})(S)=\overline{\sigma}_jf_{m-1}(S)+\overline{\sigma}_jf_m(S)< \epsilon_{m-1}+\epsilon_m.
\]
Set $u=\sum\limits_{m\ge 1}h_{2m}$ and $v=\sum\limits_{m\ge 0}h_{2m+1}$, then $u,v\in \mathcal{M}_j(S)$, since
\[
\|u\|_j(S)=\overline{I}_ju(S)\le\overline{I}_j[\sum\limits_{m\ge 1}h_{2m}](S) < \sum\limits_{m\ge 1}(\epsilon_{2m-1}+\epsilon_{2m})=\epsilon,
\]
and
\[
\overline{I}_jv(S)\le \sum\limits_{m\ge 0}\overline{I}_j h_{2m+1}(S)< \overline{I}_jf+\epsilon_1+\sum\limits_{m\ge 1}(\epsilon_{2m}+\epsilon_{2m+1})=\overline{I}_jf+\epsilon<\infty.
\]
Thus, on $\Se_{(S,j)}$, for $n\ge 1$
\[
\sum_{m=1}^n(-1)^{m-1}h_m = \sum_{m=1}^n ((-1)^{m-1}f_{m-1}-(-1)^mf_m) = f-(-1)^nf_n,
\]
from where
\[
f= \sum_{m=1}^{2n}(-1)^{m-1}h_m +f_{2n} = \sum_{m=0}^{n-1}h_{2m+1} - \sum_{m=1}^{n}h_{2m}+f_{2n}.
\]
Finally, taking the limit  $n\rightarrow\infty$, it follows that $f=v-u \;\; a.e.~~on~~\Se_{(S,j)}$, which gives $f \in \mathcal{L}_j$.

\vspace{.15in}
In order to complete the proof of (\ref{LsigmaCharacterization}), let $f\in\mathcal{L}_{\overline{\sigma}_j}$ and so $\mathcal{N}\equiv \{S\in\Se: \overline{\sigma}_j f(S)-\underline{\sigma}_j f(S)\neq 0\}$ is a null set and $\overline{\sigma}_jf(S)<\infty$ if $S\in \mathcal{N}^c$.
For any  $S \in \mathcal{N}^c$, by definition of $\overline{\sigma}_j$ there exist $\tilde{h}\in \mathcal{E}_{(S,j)}, \tilde{w}\in \mathcal{M}_j(S)$ such that $f\le \tilde{h}+\tilde{w}$ on $\Se_{(S,j)}$, with $I_j\tilde{h}(S)+\overline{I}_j\tilde{w}(S)<\overline{\sigma}_jf(S)+\epsilon$. Set $\tilde{f}\equiv \tilde{h}+\tilde{w}-f$ on $\Se_{(S,j)}$ and $0$ on $\Se\setminus\Se_{(S,j)}$. By Proposition \ref{Properties_L} item $(4)$ $\overline{\sigma}_j\tilde{h}(S)+ \overline{\sigma}_j(-\tilde{h})(S)=0$, then reasoning as in the case of $f_m$ before, $\overline{\sigma}_j\tilde{f}(S)- \underline{\sigma}_j\tilde{f}(S)=0$ and $\tilde{f}\in \mathcal{L}_{\overline{\sigma}_j}^+$.

So, by (\ref{L+Characterization}) $\tilde{f}=\tilde{v}-\tilde{u}  ~~a.e.~~on~~\Se_{(S,j)}$, with $\tilde{u},\tilde{v}\in \mathcal{M}_j(S)$ and $\|\tilde{u}\|_j(S)<\epsilon$. Consequently $f=\tilde{h}+\tilde{w}+\tilde{u}-\tilde{v} \;\; a.e.~~on~~\Se_{(S,j)}$.

Then, by Proposition \ref{firstComputationOfCI}, $\tilde{w}+\tilde{u} = \tilde{h}_1+\tilde{v}_1$, $\tilde{v}= \tilde{h}_2+\tilde{v}_2$ with ~~ $\tilde{h}_1, \tilde{h}_2\in \mathcal{E}_{(S,j)}^+, \tilde{v}_1,\tilde{v}_2\in \mathcal{M}_{(S,j)}$, such that $\overline{I}_j\tilde{v}_1(S),\overline{I}_j\tilde{v}_2(S)<\epsilon$, and
\[
f= \tilde{h}+\tilde{h}_1-\tilde{h}_2 + \tilde{v}_1-\tilde{v}_2  \;\; a.e.~~on~~\Se_{(S,j)},
\]
notice that $h\equiv (\tilde{h}+\tilde{h}_1-\tilde{h}_2) \in \mathcal{E}_{(S,j)}$ and so, we have established $f \in \mathcal{L}_j^K$.
\end{proof}

\vspace{.2in}
\noindent Fixed any node $(S,j)$  let $\Fe_{(S,j)} \equiv \{f|_{\Se_{(S,j)}}: f\in Q,~ \|f\|_j(S)<\infty\}$,  where the functions that are equal a.e. conditionally are identified. $(\Fe_{(S,j)},\|.\|_j(S))$ with pointwise operations, becomes a linear normed  space thanks to Propositions \ref{propertiesOfIBarra}  and \ref{leinertTheorem}. Moreover from \cite[Theorem 3]{ferrando} it is also complete. Also define $\Fe_j \equiv \{f\in Q : f|_{\Se_{(S,j)}}\in \Fe_{(S,j)}\}$.


\begin{corollary} \label{normProperties} Consider any  node $(S,j)$.
\begin{itemize}
\item[a)] $\overline{\sigma}_j$ is linear a.e., continuous and positive on $\mathcal{L}_{\overline{\sigma}_j}$.
\end{itemize}
\noindent
{\mbox The following statements require the property: $\overline{\sigma}_j f = \overline{I}_j f~~\mbox{for all}~~f \geq 0$. Then:}

\begin{enumerate}
\item $\mathcal{L}_{\overline{\sigma}_j}|_{\Se_{(S,j)}}\subset \Fe_{(S,j)}$ is a complete subspace under $\|.\|_j(S)$.
\item $\mathcal{L}_{\overline{\sigma}_j}|_{\Se_{(S,j)}}$ is the $\|.\|_j(S)$-closure of $\Ee_{(S,j)}$.
\item $\mathcal{L}_{\overline{\sigma}_j}^+|_{\Se_{(S,j)}}$ is the $\|.\|_j(S)$-closure of $\Ee_{(S,j)}^+$.
\end{enumerate}
 Observe that $(2)$, Proposition \ref{Properties_L} item (5) and Corollary \ref{good_properties} item $(4)$ show that $\overline{\sigma}_j$ is the unique linear extension of $I_j$ from $\Ee_{(S,j)}$ onto  $\mathcal{L}_{\overline{\sigma}_j}|_{\Se_{(S,j)}}$.
\end{corollary}
\begin{proof}
Linearity $a.e.$ in item $a)$ follows from Proposition \ref{Konig2.1}, continuity and positivity from Corollary \ref{good_properties} and general hypothesis.

Let $f\in \mathcal{L}_{\overline{\sigma}_j}$, then by Theorem \ref{punctualCharacterization}, for any $n\ge 1$ there exist $h_n\in\Ee_{(S,j)}$, $u_n,v_n \in \mathcal{M}_j(S)$ such that $f=h_n+v_n-u_n  ~~a.e.~~at~~\Se_{(S,j)}$, and $\overline{I}_ju_n(S),\overline{I}_jv_n(S)<\frac1n$. Since $h_n=h^1_n-h^2_n$ with $h^1_n,h^2_n\in\Ee^+_{(S,j)}$, it follows that
\[
\|f\|_j(S)\le \|h^1_n\|_j(S)+\|h_n^2\|_j(S)+\|v_n\|_j(S)+\|u_n\|_j(S) < \infty.
\]
This shows that $f|_{\Se_{(S,j)}}\in \Fe_{(S,j)}$. Moreover by Theorem \ref{Konig2.9}, $\mathcal{L}_{\overline{\sigma}_j}|_{\Se_{(S,j)}}$ is a closed subspace of $(\Fe_{(S,j)},\|.\|_j(S))$, and so complete. This gives item $(1)$.

For item $(2)$ observe that $f-h_n=\mathbf{1}_{A^c}[v_n-u_n]+\mathbf{1}_{A}[f-h_n]$ on $\Se_{(S,j)}$ holds for the null set $A=\{f\neq h_n+v_n-u_n \}\subset\Se_{(S,j)}$. This implies that
\[
|f-h_n|\le\mathbf{1}_{A^c}|v_n-u_n|+\mathbf{1}_{A}|f-h_n|\quad\mbox{from where}\quad \|f-h_n\|_j(S)\le {\small\mbox{$\frac 2n$}},
\]
where Proposition \ref{leinertTheorem} item $(3)$ was used because $\mathbf{1}_{A^c}|v_n-u_n|=|v_n-u_n|$ a.e. on $\Se_{(S,j)}$.

The proof of item $(3)$ is similar. Any $f\in \mathcal{L}_{\overline{\sigma}_j}^+$ can be written as $f=v_n-u_n  ~~a.e.~~at~~\Se_{(S,j)}$, $u_n,v_n \in \mathcal{M}_j(S)$, with $\overline{I}_ju_n(S)<\frac1n$. Also, from Proposition \ref{firstComputationOfCI}, $v_n=h_n+\tilde{v}_n$ with $h_n\in \Ee_j^+$, $\tilde{v}_n \in \mathcal{M}_j$ and $\overline{I}_j\tilde{v}_n(S)<\frac1n$.
\end{proof}

The classical Beppo-Levi and monotone convergence theorems holds in $\mathcal{L}_{\overline{\sigma}_j}$.
\begin{proposition}\label{BeppoLevi} Assume for $j\ge 0$, $\overline{\sigma}_j f = \overline{I}_j f~~\mbox{for all}~~f \geq 0$. Let $\{f_n \}_{n\ge 1} \subseteq\mathcal{L}_{\overline{\sigma}_j}$ such that $\sum\limits_{n\ge 1}\|f_n\|_j<\infty$; then $\sum\limits_{n\ge 1}^m f_n$ converges a.e. and in the norm $\|.\|_j(S)$ on each conditional space $\Se_{(S,j)}$ to $f \equiv \sum\limits_{n\ge 1}^{\infty} f_n$. Moreover, $f \in \mathcal{L}_{\overline{\sigma}_j}$ and
\[\int_j\sum_{n\ge 1} f_n \doteq \sum_{n\ge 1}\int_j f_n.\]
\end{proposition}
\begin{proof}
From hypothesis, since by Corollary \ref{normProperties}, $\mathcal{L}_{\overline{\sigma}_j}|_{\Se_{(S,j)}}\subset\Fe_{(S,j)}$, item $(2)$ of \cite[Proposition 5]{ferrando} gives that $f|_{\Se_{(S,j)}} \equiv \sum\limits_{n\ge 1} f_n|_{\Se_{(S,j)}}$ exists pointwise a.e. on ${\Se_{(S,j)}}$, and converges to $f|_{\Se_{(S,j)}}\in\Fe_{(S,j)}$ in the norm $\|.\|_j(S)$. So completeness implies $f \in \mathcal{L}_{\overline{\sigma}_j}$.

Linearity $a.e.$ and continuity of $\overline{\sigma}$,
given by Corollary \ref{normProperties} item $a)$, imply
\[|\int_j f - \sum_{n\ge 1}^m \int_j f_n| \doteq |\overline{\sigma}_j(f - \sum_{n\ge 1}^m f_n)|\le \overline{\sigma}_j|f - \sum_{n\ge 1}^m f_n|\le \|f - \sum_{n\ge 1}^m f_n\|_j\xrightarrow[m\rightarrow\infty]{} 0, \]
where the inequalities are $a.e.$
\end{proof}
\begin{proposition}[Monotone Convergence Theorem]\label{mct}
Assume $\overline{\sigma}_j f = \overline{I}_j f~~\mbox{for all}~~f \geq 0$. Let $\{f_n \}_{n\ge 1} \subseteq \mathcal{L}_{\overline{\sigma}_j}$, $f_n \nearrow f\in Q$. If for a.e. $S\in\Se$, $\int_j f_n \leq C = \mbox{constant}<\infty$~~ a.e. on $\Se_{(S,j)}$. Then $\|f - f_n\|_j \rightarrow 0$, $f \in \mathcal{L}_{\overline{\sigma}_j}$ and $\int_j f \doteq \lim_{n \rightarrow \infty} \int_j f_n$.
\end{proposition}
\begin{proof}
Define for $n\ge 1$, $g_n \equiv f_{n+1} - f_n \geq 0$, then for a.e. $S\in\Se$
\[
\sum_{n=1}^m\|g_n\|_j=\sum_{n=1}^m\overline{\sigma}_jg_n\doteq\overline{\sigma}_jf_{m+1}-\overline{\sigma}_jf_1=\int_jf_{m+1}-\int_jf_1<\infty~~ ~~a.e.~~~~ \mbox{on}~~~~\Se_{(S,j)}.
\]
The result then follows from Proposition \ref{BeppoLevi}.
\end{proof}
\section{Basic Properties of Conditional Integrals}
\label{basicPropertiesOfConditionalIntegrals}

\begin{lemma}\label{integProperty}  Let $0 \leq j \le k $ and $f_V\in Q$ satisfy: $f_V(S)=f_V(S_0,...,S_j)=V(S_0,...,S_j)$. Then $\overline{\sigma}_k f_V(S) = \underline{\sigma}_k f_V(S)=  V(S)$. In particular for $f\in Q$, since $\overline{\sigma}_j f$ depends on $S$ just through $S_0, \ldots,S_j$, it follows that,
\[
 \overline{\sigma}_k [\overline{\sigma}_j f] = \overline{\sigma}_j f = \underline{\sigma}_k [\overline{\sigma}_j f]\quad \mbox{and} \quad \overline{\sigma}_k [\underline{\sigma}_j f]= \underline{\sigma}_j f = \underline{\sigma}_k [\underline{\sigma}_j f].
\]
Furthermore if $f\in \mathcal{L}_{\overline{\sigma}_j}$, then $~\int_j f\in \mathcal{L}_{\overline{\sigma}_k}$ and
\[
\int_k \int_ j f~~\equiv \int_k[\int_j f]= \int_j f.
\]
\end{lemma}
\begin{proof}
Writing $f_V = V = V^+ - V^-$ we see that $f_V \in \mathcal{E}_j \subseteq
\mathcal{E}_k$, an application of Proposition \ref{Properties_L} item $(4)$ gives
$\overline{\sigma}_k f_V= \underline{\sigma}_k f_V= I_k f_V= V$.
\end{proof}

\begin{proposition}  \label{mainDirectionOfTowerProperty}
Let $f \in Q$, and $0 \leq j \leq k$, then
\begin{equation}  \label{tower_ineq}
\underline{\sigma}_j f \le \underline{\sigma}_j[\underline{\sigma}_k f]\le \underline{\sigma}_j[\overline{\sigma}_k f] \leq \overline{\sigma}_j[\overline{\sigma}_k f]\le \overline{\sigma}_j f
\end{equation}
and
\begin{equation}  \label{tower_ineq2}
\underline{\sigma}_j f \le \underline{\sigma}_j[\underline{\sigma}_k f]\le \overline{\sigma}_j[\underline{\sigma}_k f] \leq \overline{\sigma}_j[\overline{\sigma}_k f]\le \overline{\sigma}_j f.
\end{equation}
Therefore,  if $f \in \mathcal{L}_{\overline{\sigma}_j}$ then $\underline{\sigma}_k f,\overline{\sigma}_k f \in \mathcal{L}_{\overline{\sigma}_j}$ and
\begin{equation}  \label{int_tower_ineq}
\int_j [\underline{\sigma}_k f] \doteq \int_j f \doteq \int_j [\overline{\sigma}_k f].
\end{equation}
\end{proposition}
\begin{proof} Assume $f\le \sum_{m\ge 0} f_m$, with $f_0\in \mathcal{D}^-_j,~~f_m\in \mathcal{D}^+_j, m\ge 1$ and $f_m = \Pi^{V^m, H^m}_{j,n_m}$. Recall that we have assumed that $(H^m_i|_{\Se_{(S,k)}})_{i\ge k}\in\He_{(S,k)}$ for any $S\in\Se$.

Then also $f_0\in \mathcal{D}^-_k,~~f_m\in \mathcal{D}^+_k, m\ge 1$, and it results that
\[
\overline{\sigma}_k f (S) \le \sum_{m\ge 0} \left(V^m(S)+\sum_{i=j}^{k-1}H^m_i(S)\Delta_iS\right) = \sum_{m\ge 0} I_kf_m(S).
\]
Observing that by Corollary \ref{representationUniqueness} $I_k f_0\in \mathcal{D}^-_j,~ I_k f_m \in \mathcal{D}^+_j~~ m\ge 1$, and $V^m(S)=V^m(S_0,...,S_j)$, it follows that $\overline{\sigma}_j[\overline{\sigma}_k f](S) \le \sum_{m\ge 0} V^m(S)$, and consequently
\begin{equation}  \label{tower_ineq_1}
\overline{\sigma}_j[\overline{\sigma}_k f] \le \overline{\sigma}_j f.
\end{equation}
Applying (\ref{tower_ineq_1}) to $-f$, taking into consideration definitions of $\underline{\sigma}_{\ast}$ and established inequalities, it follows that
\begin{equation} \nonumber 
\underline{\sigma}_j f \le - \overline{\sigma}_j[\overline{\sigma}_k(-f)]= \underline{\sigma}_j[\underline{\sigma}_k f] \le \underline{\sigma}_j[\overline{\sigma}_k f] \le \overline{\sigma}_j[\overline{\sigma}_k f] \leq \overline{\sigma}_j f.
\end{equation}
Similarly
\begin{equation}  \nonumber 
\underline{\sigma}_j f \le - \overline{\sigma}_j[\overline{\sigma}_k(-f)]= \underline{\sigma}_j[\underline{\sigma}_k f] \le \overline{\sigma}_j[\underline{\sigma}_k f] \le \overline{\sigma}_j[\overline{\sigma}_k f] \leq \overline{\sigma}_j f.
\end{equation}
So, (\ref{tower_ineq}) and (\ref{tower_ineq2}) hold and the remaining statements follow.
\end{proof}

Notice that the hypothesis $\overline{I}= \overline{\sigma}$ (on non-negative functions)
required in the next corollary follows from property $(K_0)$ (as proven in Appendix \ref{proofOfK}). A similar remark applies to Corollary \ref{towerPropertyCorollary}.

The next corollary refers to the notion of trajectorial martingale, more details are provided in Section \ref{definitionsOfTM}. We recall that we use the notation $\overline{\sigma}= \overline{\sigma}_0$.

\begin{corollary}\label{integrableImpliesConditionalIntegrable}
Assume $\overline{I}= \overline{\sigma}$ on non-negative functions. If $f \in \mathcal{L}_{\overline{\sigma}}$ then, $f \in  \mathcal{L}_{\overline{\sigma}_k}$ for all $k \geq 0$ and $\{f_k \equiv \int_k f\}$ is a trajectorial martingale i.e.
$\int_k[\int_{k+1} f]\doteq \int_k f_{k+1} \doteq f_k$.
\end{corollary}
\begin{proof}
We compute
\begin{equation} \nonumber
\overline{\sigma}(\overline{\sigma}_k f - \underline{\sigma}_k f)=   \overline{\sigma}(\overline{\sigma}_k f)+ \overline{\sigma}(-\underline{\sigma}_k f)= \overline{\sigma}(\overline{\sigma}_k f)- \underline{\sigma}(\underline{\sigma}_k f),
\end{equation}
where we have used Lemma  \ref{linearityOfSigma} which is applicable because $\overline{\sigma}(\overline{\sigma}_kf)(S)-
\underline{\sigma}(\overline{\sigma}_kf)(S)=0$ and
$\underline{\sigma}(\underline{\sigma}_kf)(S)-
\overline{\sigma}(\underline{\sigma}_kf)(S)=0$ hold, for all $S$, by taking $j=0$ in
(\ref{tower_ineq}) and (\ref{tower_ineq2}) respectively (actually we could have replaced $S$ by $S_0$, i.e. the quantities are constant).
Our hypothesis $f \in \mathcal{L}_{\overline{\sigma}}$  together with (\ref{tower_ineq}) and (\ref{tower_ineq2}) imply that
$\overline{\sigma}(\overline{\sigma}_k f)- \underline{\sigma}(\underline{\sigma}_k f)=0$. Therefore, given that $\overline{\sigma}_k f - \underline{\sigma}_k f \geq 0$ we have
$\overline{I}(\overline{\sigma}_k f - \underline{\sigma}_k f)=
\overline{\sigma}(\overline{\sigma}_k f - \underline{\sigma}_k f)=0$ which implies $\overline{\sigma}_k f - \underline{\sigma}_k f= 0$ ~~a.e., this means $f \in \mathcal{L}_{\overline{\sigma}_k}$.
Moreover from Proposition \ref{mainDirectionOfTowerProperty}:\\
\[
\int_k f_{k+1}\doteq \int_k[\int_{k+1}f] \doteq \int_k f = f_k.
\]
\end{proof}

\begin{corollary}  \label{towerPropertyCorollary}
Assume $\overline{I}_j = \overline{\sigma}_j$ on non-negative functions. Let $f\in \mathcal{L}_{\overline{\sigma}}$ and $0\le j \le k$; then
\begin{equation}\label{Tower_property}
\int_k[\int_j f] = \int_j f \doteq \int_j[\int_k f].
\end{equation}
Moreover, assuming also that $|f|\in \mathcal{L}_{\overline{\sigma}}$ it follows that
\begin{equation}\label{conditionallyNullImpliesNull}
\mbox{If $||f||_k ~\dot{=} ~0$ then $f$ is a null function.}
\end{equation}
\end{corollary}
\begin{proof}
Assume first that $f\in \mathcal{L}_{\overline{\sigma}}$ then $f \in \mathcal{L}_{\overline{\sigma}_p}$ for all $p\ge 0$, by Corollary \ref{integrableImpliesConditionalIntegrable}. Thus the left equality of (\ref{Tower_property}) holds by Lemma \ref{integProperty}, while the right hand one follows from (\ref{int_tower_ineq}) in Proposition \ref{mainDirectionOfTowerProperty}.

For (\ref{conditionallyNullImpliesNull}), by Proposition \ref{konig} item $(\ref{sigma_le_Iup})$: $0 \leq \overline{\sigma}_k(|f|) \leq  \|f\|_k \dot{=} 0$. Then by Proposition \ref{mainDirectionOfTowerProperty}, and (\ref{sigma=I}),
\[
0= \overline{\sigma}_0[\overline{\sigma}_k(|f|)]= \overline{\sigma}_0(|f|)=\|f\|.
\]

\end{proof}

\section{Trajectorial Martingales}  \label{definitionsOfTM}


\begin{definition}
Consider a sequence of functions $f_n:\mathcal{S} \rightarrow [- \infty, \infty]$
satisfying $f_n(S)= f(S_0, \ldots, S_n)$. $\{f_n\}$ is a supermartingale if
\begin{equation} \nonumber
\overline{\sigma}_j f_{j+1} \leq f_j~~a.e.
\end{equation}
$\{f_n\}$ is a submartingale if
\begin{equation} \nonumber
f_j \leq \underline{\sigma}_j f_{j+1}~~a.e.
\end{equation}
$\{f_n\}$ is a martingale if
\begin{equation} \nonumber
\underline{\sigma}_j f_{j+1}= \overline{\sigma}_j f_{j+1}= f_j~~a.e.
\end{equation}
\end{definition}

\begin{proposition}
Let $f: \mathcal{S} \rightarrow [- \infty, \infty]$.

\begin{itemize}
\item $f_j \equiv \underline{\sigma}_j f$ is a submartingale
sequence.

\item $f_j \equiv \overline{\sigma}_j f$ is a supermartingale
sequence.

\item Assume $\overline{I}= \overline{\sigma}$ on non-negative functions and $f \in \mathcal{L}_{\overline{\sigma}}$. Then
$f_j \equiv \overline{\sigma}_j f$ is a martingale sequence.
\end{itemize}
\end{proposition}
\begin{proof}

\noindent
Let $f_j \equiv \underline{\sigma}_j f$, taking $k = j+1$ in (\ref{tower_ineq})
from Proposition  \ref{mainDirectionOfTowerProperty} gives
$f_j \leq \underline{\sigma}_j f_{j+1}$.

\vspace{.1in}
\noindent Let $f_j \equiv \overline{\sigma}_j f$, taking $k = j+1$ in (\ref{tower_ineq})
from Proposition  \ref{mainDirectionOfTowerProperty} gives
$\overline{\sigma}_j f_{j+1} \leq f_j$.

\vspace{.1in}
\noindent
Assume $\overline{I}= \overline{\sigma}$ on non-negative functions, then  Corollary \ref{integrableImpliesConditionalIntegrable} implies
 $f \in \mathcal{L}_{\overline{\sigma}_j}~~\mbox{for all}~~j \geq 0$.
In particular $-\infty < \underline{\sigma}_j f = \overline{\sigma}_j f < \infty$ holds a.e.

Then, taking $k =j+1$ in (\ref{tower_ineq}) 
from Proposition  \ref{mainDirectionOfTowerProperty} gives
$-\infty < f_j = \overline{\sigma}_j f_{j+1} < \infty$ holds a.e.
\end{proof}

\section{Conclusions} \label{conclusions}

The paper is based on the observation that hedging prices of simple portfolios, i.e. the elementary integral, can be extended to a larger class of functions under the idealization of taking the limit of  increasing the number of portfolios used in the approximations. This approach leads to upper and lower approximations providing a natural class of integrable functions by the method of exhaustion. We show that the theory mimics the classical theory of integration but in this case providing a pricing interpretation for the integral.
The meaning of the constructed integral follows from the elementary integrable functions representing initial portfolio values.

The integration framework naturally leads to the concept of trajectorial martingale. There is also the possibility to introduce a further natural idealization into the foundations of the basic constructions namely, introducing another, independent, limit by increasing the number of rebalancing times. We plan to pursue these research possibilities in future work.

\appendix

\section{Some Proofs for Sections \ref{integralConstruction}-\ref{IntegrabilityCharacterization}} \label{someProof1}
\noindent
{\bf Proof of Proposition \ref{propertiesOfIBarra}}:
\begin{proof} We only provide the proof of countable subadditivity.

We may assume that $\sum\limits_{k\ge 1}\Iup_j g_k(S) <\infty$, which leads to $\Iup_j g_k(S)<\infty$ for $k\ge 1$. Therefore, for a fixed $\epsilon>0$ and for any $k\ge 1$, by definition of $\Iup_j$ there exist $H^{m,k}\in\He_{(S,j)}$ and $V^{m,k}\in \mathbb{R}$, $m\ge 1$, such that:
\[g_k(\tilde{S})\le  \sum_{m=1}^{\infty} [V^{m,k} + \sum_{i=j}^{n_{m,k}-1} H^{m,k}_i(\tilde{S}) \Delta_i \tilde{S}]\quad \forall~ \tilde{S} \in \Se_{(S,j)},\]
with
\[
V^{m,k} + \sum_{i=j}^{n_{m,k}-1} H^{m,k}_i(\tilde{S}) \Delta_i \tilde{S} \ge 0 \quad \forall \tilde{S} \in\Se_{(S,j)},\quad \mbox{and}\quad  \sum_{m=1}^{\infty}V^{m,k}\le \overline{I}_jg_k(S) + \frac{\epsilon}{2^{k}}.
\]
Then
\begin{equation} \nonumber
\sum_{k=1}^{\infty}g_k(\tilde{S})\le \sum_{k=1}^{\infty}\sum_{m=1}^{\infty}~[V^{m,k} + \sum_{i=j}^{n_{m,k}-1} H^{m,k}_i(\tilde{S}) \Delta_i \tilde{S}] \quad \forall \tilde{S} \in\Se_{(S,j)},
\end{equation}
noticing that the double sum of nonnegative terms can be reordered into a single sum, we can then deduce that
\[
\overline{I}_j g(S) \leq \sum_{k=1}^{\infty}\sum_{m=1}^{\infty}V^{m,k}\le \sum_{k=1}^{\infty} \overline{I}_jg_k(S)+ \epsilon.
\]
\end{proof}

\noindent
{\bf Proof of Proposition \ref{leinertTheorem}}:
\begin{proof}
$(1)$ Assume $||g||_j(S)=0$, consider $A=\{\tilde{S} \in \Se_{(S,j)}: g(\tilde{S}) \neq 0\}$.\\ From ${\bf 1}_A(\tilde{S}) \leq \sum_{k \geq 1} |g(\tilde{S})|$ it follows that
$||{\bf 1}_A||_j(S) \leq \sum_{k \geq 1} ||g||_j(S)=0$. Therefore, $A$ is a conditionally null set at $(S,j)$ and so $g(\tilde{S})=0$ holds conditionally a.e. at $(S,j)$.

\noindent For the converse of $(1)$, by assumption, there exists $B \subseteq \Se_{(S,j)}$ such that $||{\bf 1}_B||_j(S)=0$ and
$g(\tilde{S}) =0 \;\; \forall \tilde{S} \in \Se_{(S,j)} \setminus B$.
Given that $|g(\tilde{S})| \leq \sum_{k \geq 1} {\bf 1}_{B}(\tilde{S}) \;\; \forall \tilde{S} \in \Se_{(S,j)}$ we obtain $||g||_j(S)=0$.

\vspace{.1in}
\noindent
$(2)$  Let  $A \equiv \{\tilde{S} \in \mathcal{S}_{(S,j)}: g(\tilde{S})= \infty\}$, then $n ~{\bf 1}_A(\tilde{S}) \leq |g(\tilde{S})|, \; \forall \tilde{S} \in \Se_{(S,j)},\; \forall n\ge 1$, thus
\[
n ~||{\bf 1}_A||_j(S)  \leq ||g||_j(S)\quad \mbox{and so}\quad ||{\bf 1}_A||_j(S)=0.
\]

\noindent $(3)$ Let $N\equiv \{\tilde{S} \in \Se_{(S,j)}: |f(\tilde{S})| > |g(\tilde{S})|\}$.

Then $|f|(\tilde{S}) \leq |g|(\tilde{S}) + \sum\limits_{m\ge 1} {\bf 1}_N(\tilde{S})$ for $\tilde{S} \in \Se_{(S,j)}$. Therefore:

$\|f\|_j(S)=\overline{I}_j|f|(S)  \leq  \overline{I}_j |g|(S) +  \overline{I}_j( \sum\limits_{m\ge 1} {\bf 1}_N )(S) =  \overline{I}_j|g|(S)=\|g\|_j(S)$.

\noindent $(4)$ Follows from the countable subadditivity of $\overline{I}_j$.

\noindent $(5)$  Assume $f\le \sum_{m\ge 1} f_m$, with $f_m = \Pi^{V^m, H^m}_{j,n_m}\in \mathcal{D}^+_j, m\ge 1$.

Since also $f_m\in \mathcal{D}^+_k, m\ge 1$, it results that for any $S\in\Se$

\[
\overline{I}_k f (S) \le \sum_{m\ge 1} \left(V^m(S)+\sum_{i=j}^{k-1}H^m_i(S)\Delta_iS\right) = \sum_{m\ge 1} I_kf_m(S).
\]
Since $ I_k f_m \in \mathcal{D}^+_j,~ m\ge 1$, and $V^m(S)=V^m(S_0,...,S_j)$, it conducts to $\overline{I}_j[\overline{I}_k f](S) \le \sum_{m\ge 0} V^m(S)$ and consequently
\begin{equation}  \label{tower_ineqForIBar}
\overline{I}_j[\overline{I}_k f] \le \overline{I}_j f.
\end{equation}
Assume now that $g \in Q$ is conditionally null at $\mathcal{S}_{(S,j)}$
i.e. $\overline{I}_j(|g|)(S)=0$. It then follows from (\ref{tower_ineqForIBar}) and and item $(1)$ that
$\overline{I}_k(|g|)=0 $ a.e. on $\mathcal{S}_{(S,j)}$.

\end{proof}

\noindent
{\bf Proof of Proposition \ref{konig}}:

\begin{proof}

\noindent
$(1)$ follows from the definitions and $(2)$ is clear taking $f_0=0$.

\noindent For $(3)$, we can write $f=f_1+f_0$ with $f_0\in \mathcal{D}_{(S,j)}^-=- \mathcal{D}_{(S,j)}^+$ and $f_1\in \mathcal{D}_{(S,j)}^+$, taking $f_m=0$ for $m\ge 2$ the result follows.

\noindent
To establish $(4)$, notice that $f \leq g,~~a.e. ~\mbox{on}~ \mathcal{S}_{(S,j)}$ allows to write $f(\tilde{S}) \leq g(\tilde{S}) + \infty ~ {\bf 1}_N(\tilde{S})$ for all $\tilde{S} \in \Se_{(S,j)}$  with $N \subseteq \Se_{(S,j)}$ and
$\overline{I}_j{\bf 1}_N(S) =0$. Therefore $\overline{\sigma}_j f(S) \leq
\overline{\sigma}_j g(S) + \overline{\sigma}_j(\infty ~ {\bf 1}_N)(S) \leq
\overline{\sigma}_j g(S) + \overline{I}_j(\infty ~ {\bf 1}_N)(S) \leq \overline{\sigma}_j g(S)$.

\noindent
$(5)$ Assume $gf(\tilde{S}) \le \sum\limits_{m\ge 0}\Pi_{j, n_m}^{V^m, H^m}(\tilde{S}),\; \tilde{S}\in\Se_{(S,j)}$,
with $\Pi_{j, n_0}^{V^0, H^0}\in \mathcal{D}_{(S,j)}^-$ and, for $m\ge 1$, $\Pi_{j, n_m}^{V^m, H^m}\in \mathcal{D}^+_{(S,j)}$.
For each $\tilde{S} \in \Se_{(S,j)}$ and ~~ $m\ge 0$ ~~~~ define
\[
U^m(\tilde{S}) = \frac {V^m}{g(S)},\;\;\mbox{and}\;\; G^m_i(\tilde{S})=\frac {H^m_i(\tilde{S})}{g(S)},\;\mbox{for}\; i\ge j.
\]
It follows that $f(\tilde{S})\le \sum\limits_{m\ge 0}\Pi_{j, n_m}^{U^m, G^m}(\tilde{S}),\; \tilde{S}\in\Se_{(S,j)}$ with $\Pi_{j, n_0}^{U^0, G^0}\in \mathcal{D}_{(S,j)}^-$, and for $m\ge 1$, $\Pi_{j, n_m}^{U^m, G^m}\in \mathcal{D}^+_{(S,j)}$. Thus
\[
\overline{\sigma}_jf(S) \le \frac{\overline{\sigma}_j[gf](S)}{g(S)}.
\]
The reverse inequality follows similarly.

\noindent
$(6)$ Observe that $f \le |f-g|+g$~ (as $f(\hat{S})=\infty$ implies $|f-g|(\hat{S})=\infty$), from where, $(6)$ holds.

\noindent
$(7)$ Taking $f_m\equiv 0$ for any $m\ge 0$, it follows that $\overline{\sigma}_j 0\le 0$.  In case that $\overline{\sigma}_j 0 \neq 0$, it then follows that there exist at least some $S \in \mathcal{S}$ satisfying  $\overline{\sigma}_j 0 (S) < 0$.
Therefore, there exist $\tilde{f}_0=\Pi_{j, \tilde{n}_0}^{\tilde{V}^0, \tilde{H}^0}\in \mathcal{D}_{(S,j)}^-$ and, for $m\ge 1$, $\tilde{f}_m=\Pi_{j, \tilde{n}_m}^{\tilde{V}^m, \tilde{H}^m}\in \mathcal{D}_{(S,j)}^+$, such that $0 \le \sum\limits_{m\ge 0} \tilde{f}_m(\tilde{S})$ for any $\tilde{S}\in \Se_{(S,j)}$ and $\sum\limits_{m\ge 0} \tilde{V}^m=r < 0$.

Consider now $f\in Q$, if there exists $f_0=\Pi_{j, n_0}^{V^0, H^0}\in \mathcal{D}_{(S,j)}^-$ and $f_m=\Pi_{j, n_m}^{V^m, H^m}\in \mathcal{D}_{(S,j)}^+$ for $m\ge 1$, such that $f(\tilde{S}) \le \sum\limits_{m\ge 0} f_m(\tilde{S})$ for any $\tilde{S}\in \Se_{(S,j)}$ with $\sum\limits_{m\ge 0} V^m$ finite, then (since also $f(\tilde{S}) \le \sum\limits_{m\ge 0} f_m(\tilde{S})+ \gamma~ \tilde{f}_m(\tilde{S})$, for any given $\gamma >0$):
\[
\overline{\sigma}_j f (S) \leq \sum_{m\ge 0} (V^m+ \gamma ~\tilde{V}^m)= \sum_{m\ge 0} V^m + \gamma ~r,
\]
thus $\overline{\sigma}_jf(S)=-\infty$. On the other hand if such family of functions $f_m$ does not exist, $\overline{\sigma}_jf(S)=\infty$.



\noindent
$(8)$ Let $A=\{\hat{S}\in\Se_{(S,j)}:f(\hat{S})=\infty\}$ and assume that $\overline{\sigma}_jf(S)<\infty$. There exist  $f_0=\Pi_{j, n_0}^{V^0, H^0}\in \mathcal{D}_{(S,j)}^-$, and for $m\ge 1$, $f_m=\Pi_{j, n_m}^{V^m, H^m}\in \mathcal{D}_{(S,j)}^+$, such that $f(\tilde{S}) \le \sum_{m\ge 0} f_m(\tilde{S})$ for any $\tilde{S}\in \Se_{(S,j)}$ and
\[
\sum_{m\ge 0} V^m < \overline{\sigma}_jf(S)+1.
\]
For $\hat{S}\in A$ results that $~~\infty=f(\hat{S})\le\sum_{m\ge 1} f_m(\hat{S})$, because $|f_0|<\infty$.
Thus, for any $n > 0$ and $\tilde{S}\in\Se_{(S,j)}$,  it follows that
\[
n \mathbf{1}_A(\tilde{S})\le\sum_{m\ge 1} f_m(\tilde{S}).
\]
From which
\[n\|\mathbf{1}_A\|_j(S) \le \sum_{m\ge 1} V^m < \infty.\]

\noindent
$(9)$ From hypothesis $0=\overline{\sigma}_j 0\le \overline{\sigma}_jf(S)+\overline{\sigma}_j(-f)(S)$, the two terms in the right hand side being finite by hypothesis, can then conclude that $-\infty < \overline{\sigma}_j(-f)(S),  -\infty <\overline{\sigma}_jf(S)$.
The last statement follows from item $(8)$.
\end{proof}

\noindent
{\bf Proof of Proposition \ref{Properties_L}}:
\begin{proof} $(2)$ follows from $(1)$ as follows.

Let $f=h-h' \in \mathcal{E}_j$ with $h,h' \in \mathcal{D}_j^+$, and $h_m \in \mathcal{D}_j^+$ such that $f \leq \sum_{m \geq 1} h_m$. Then $0 \leq -h + h' + \sum_{m \geq 1} h_m$, thus (taking $f_0 \equiv -h\in \mathcal{D}_j^-$, $f_1=h'$ and $f_m=h_{m-1}$ for $m\ge 2$) by Definition \ref{cond_integ_def}, $0 = \overline{\sigma}_j(0) \leq \sum\limits_{m \geq 0} I_jf_m$, which leads to $I_jf \leq \sum\limits_{m \geq 1} I_j h_m$ as required.

\vspace{.1in} \noindent Assume now that property $(L_j)$ holds, $f\in \mathcal{E}_j$ and $f^+\le \sum\limits_{m \geq 1} f_m,\; f_m \in \mathcal{D}_j^+$. Since $f\le f^+$, then
\[
I_j f \le \sum_{m \geq 1} I_jf_m,\quad \mbox{which implies}\quad I_j f \le \overline{I}_j f^+.
\]

\vspace{.1in} \noindent Let us now show that $(3)$ implies $(4)$. It follows the proof of \cite[Behauptung 1.9]{konig}, but here its (*) condition in page 449 is not needed. Fix $f_0\in \mathcal{D}^-_j$ and for $m\ge 1$, $f_m\in\mathcal{D}_j^+$ such that $f\le\sum\limits_{m\ge 0}f_m$, then from linearity of $I_j$ and $(3)$,
\[
I_jf - I_jf_0 = I_j(f - f_0) \le \overline{I}_j(f - f_0)^+.
\]
Moreover, since $(f - f_0)^+\le \sum\limits_{m\ge 1} f_m$, countable subadditivity of $\overline{I}_j$ and Proposition \ref{propertiesOfIBarra} item $(4)$, gives $\overline{I}_j(f - f_0)^+\le\sum\limits_{m\ge 1} \overline{I}_jf_m \le\sum\limits_{m\ge 1} I_jf_m$.

Therefore $I_jf\le \sum\limits_{m\ge 0} I_jf_m$, so $I_j f \le \overline{\sigma}_j f$, thus $I_j f = \overline{\sigma}_j f$ from Proposition \ref{konig} item $(\ref{sigma_le_Iup})$. Consequently $(4)$ holds, since $\underline{\sigma}_j f = -\overline{\sigma}_j (-f) = I_j f$.

\vspace{.1in} \noindent From $(4)$ $\overline{\sigma}_j 0 = I_j 0 = 0$.

\vspace{.1in} \noindent
We now incorporate items $(5)$, $(6)$ and $(7)$ into the set of equivalences. Let us see that $(L_j)$ implies $(5)$. Fix $f\in \mathcal{E}_j$ and for $m\ge 1$, $f_m\in \mathcal{D}_j^+$ such that $f\le|f|\leq \sum\limits_{m \geq 1} f_m$. Then by $(L_j)$,
$I_jf\le \sum\limits_{m\ge 0} I_jf_m$, so $I_j f \le \overline{I}_j |f|=\|f\|_j$. The same analysis gives $-I_jf=I_j[-f] \le \overline{I}_j |f|=\|f\|_j$, and so $-\|f\|_j\le I_j f \le \|f\|_j$.

$(6)$ follows from $(5)$, since for $f\in \mathcal{D}_j^+$, $0\le I_jf=|I_jf|\le \overline{I}_j |f|=\overline{I}_j f\le I_jf$, last inequality by item $(4)$ of Proposition \ref{propertiesOfIBarra}.

$(6)$ implies $(L_j^+)$. For $f, f_m\in \mathcal{D}_j^+$, such that $f\leq \sum\limits_{m \geq 1} f_m$ it follows from $(6)$ and countable subadditivity of $\overline{I}_j$ that
\[
I_jf=\overline{I}_jf \le \sum_{m \geq 1} \overline{I}_jf_m = \sum_{m \geq 1} I_jf_m.
\]

\noindent
Finally, we will prove that $(L_j^+)$ implies $\overline{\sigma}_j 0 =0$.  Clearly $\overline{\sigma}_j 0 \leq 0$ so we need only to prove the reverse inequality. So take $f_0 \in \mathcal{D}^-_j, f_m \in \mathcal{D}_j^+$ with $0 \leq \sum\limits_{m \geq 0} f_m$, since $-f_0 \in \mathcal{D}_j^+$ and $-f_0 \leq \sum\limits_{m \geq 1} f_m$. An application of $(L_j^+)$ implies
$0 \leq I_j f_0 + \sum\limits_{m \geq 1} I_j f_m$, which in turn gives $0 \leq \overline{\sigma}_j 0$.
\end{proof}

{\bf Proof of Proposition \ref{goodToKnow}}:
\begin{proof}
 Take
$f \in \mathcal{L}^K_j$, and consider $S \in \mathcal{S}$ for which the following holds: for all $\epsilon/2 >0$ there exist $u, v ,h$
(as appearing in the definition of $\mathcal{L}^K_j$)
and  $\overline{I}_j u(S) \leq \epsilon/2,~~ \overline{I}_j v(S) \leq \epsilon/2$. We compute
\begin{equation} \nonumber
- \epsilon \leq - \overline{\sigma}_j v(S)-  \overline{\sigma}_j u(S)=
\underline{\sigma}_j (-v)(S)+  \underline{\sigma}_j (-u)(S) \leq
\underline{\sigma}_j v(S)+  \underline{\sigma}_j (-u)(S)
\leq
\end{equation}
\begin{equation} \nonumber
\underline{\sigma}_j (v-u)(S) = \underline{\sigma}_j (f-h)(S)
\leq  \overline{\sigma}_j (f-h)(S)=  \overline{\sigma}_j (v-u)(S)
\leq \overline{I}_j v(S)+ \overline{I}_j u(S) \leq \epsilon
\end{equation}
where we have used $(4)$ of Proposition \ref{konig} for  two of the equalities. It follows that $f \in \tilde{\mathcal{L}}_K$.

Conversely, take now $f \in \tilde{\mathcal{L}}_K$ and consider $S \in \mathcal{S}$ for which the following holds:   for all $\epsilon >0$
let $u= \sum_{k} u_k, v= \sum_{k} v_k$ and $h$ as appearing in the definition of $\tilde{\mathcal{L}}_K$. Therefore,  we can find $n$ such that $\overline{I}_j(\sum_{k \geq n} u_k)(S) \leq \epsilon$ and $\overline{I}_j(\sum_{k \geq n} v_k)(S) \leq \epsilon$. Define $\tilde{h}\equiv h +\sum_{k=1}^{n-1} v_k- \sum_{k=1}^{n-1} u_k$ and notice that $\tilde{h} \in \mathcal{E}_j$ and  set $\tilde{u} \equiv \sum_{k \geq n} u_k$, $\tilde{v} \equiv \sum_{k \geq n} v_k$. Therefore $f = \tilde{v}- \tilde{u}+ \tilde{h}~a.e.$ on $\mathcal{S}_{(S,j)}$ and
$\overline{I}_j \tilde{u}(S) \leq \epsilon$ and $\overline{I}_j \tilde{v}(S) \leq \epsilon$ (clearly $\tilde{u}, \tilde{v} \in \mathcal{M}_{(S,j)}$).
\end{proof}

{\bf Proof of Proposition \ref{keyForConvergenceTheoremsAndCharacterization}}:
[the proof of (\ref{KonigCondition}) $\Rightarrow$ (\ref{sigma=I}) is from \cite[Behauptung 1.8]{konig} but in the conditional setting]
\begin{proof}
(\ref{KonigCondition}) $\Rightarrow$ (\ref{sigma=I}). By Proposition \ref{konig} item $(\ref{sigma_le_Iup})$, it is enough to prove that $\overline{I}_j f\le \overline{\sigma}_j f$ for $f\in P$. Fix $S\in\Se$ and $f_0\in \mathcal{E}_{(S,j)}$ and for $m\ge 1$, $f_m\in \mathcal{E}^+_{(S,j)}$, such that $f\le\sum\limits_{m\ge 0}f_m$ on $\Se_{(S,j)}$. It is enough to prove that $\overline{I}_j f(S) \le \sum\limits_{m\ge 0}I_jf_m(S)$, we can then assume $\sum\limits_{m\ge 0}I_jf_m(S)<\infty$.
For $n>0$ it holds that
\[
f \le[\sum_{m=0}^nf_m]^+ + \sum_{m\ge n+1}f_m \;\; \mbox{on}\;\; \Se_{(S,j)},\;\; \mbox{and}\;\; \overline{I}_j f(S) \le \overline{I}_j[\sum_{m\ge 0}^n f_m]^+(S) + \overline{I}_j[\sum_{m\ge n+1}f_m](S).
\]
Since $ \sum\limits_{m=0}^nf_m\in \mathcal{E}_{(S,j)}$, by property $(\ref{KonigCondition})$
\[
\overline{I}_j[\sum_{m=0}^n f_m]^+(S) = \sum_{m=0}^n I_j f_m(S) + \overline{I}_j[\sum_{m=0}^n f_m]^-(S).
\]
On the other hand, from $\sum\limits_{m\ge 0}f_m\ge f\ge 0$ it follows that
\[
(\sum_{m=0}^n f_m)^- = (-\sum_{m=0}^n f_m)^+ \le (\sum_{m\ge n+1} f_m)^+ = \sum_{m\ge n+1} f_m.
\]
Summarizing we have
\[
\overline{I}_j f(S) \le \sum_{m=0}^n I_j f_m(S) + 2\overline{I}_j[\sum_{m\ge n+1} f_m](S)\le \sum_{m=0}^n I_j f_m(S) + 2\sum_{m\ge n+1} I_j f_m(S),
\]
where countable subadditivity of $\overline{I}_j$ and Proposition \ref{Properties_L} item $(6)$ were used. Taking limits for $n\rightarrow \infty$, (\ref{sigma=I}) follows.

The implication  (\ref{sigma=I}) $\Rightarrow$ (\ref{finiteMat}) is trivial. 
It remains to prove the implication (\ref{finiteMat}) $\Rightarrow$ (\ref{KonigCondition}). Then for fixed $f\in\Ee_j$ we may assume that $\overline{I}_jf^+ <\infty$. Consider $f_m\in \mathcal{D}_j^+$ such that $f^+\le \sum\limits_{m\ge 1}f_m$. As $f_0\equiv-f\in\Ee_j$, it follows that $f^-\le \sum\limits_{m\ge 0}f_m$ thus: $\overline{\sigma}_jf^- \leq -I_j(f)+\overline{I}_jf^+$.
Since by (\ref{finiteMat}) $\overline{\sigma}_jf^-= \overline{I}_jf^-$, one inequality of $(\ref{KonigCondition})$ is obtained. The other inequality follows by applying the said inequality to $-f$.
\end{proof}

\section{Properties   $(L_j)$, $(K_j)$ and Examples} \label{proofOfK}

The property $(L_j)$ is necessary in order to define  conditional integrals and to have a non trivial theory for the given definitions. As indicated by Proposition \ref{konig}, item $(7)$,  if $\overline{\sigma}_j0(S) \neq 0$ for some $S\in\Se$ we will have $\overline{\sigma}_jf(S) = \pm \infty$ for all $f \in Q$. On the other hand, property
$(K_j)$ is equivalent to statement $(2)$ in  Proposition \ref{keyForConvergenceTheoremsAndCharacterization} and the latter can be used to obtain Lebesgue's monotone convergence theorem  and also allows to prove a characterization result (Theorem \ref{punctualCharacterization}) for our integrable functions. The financial meaning of property $(K_j)$ and its importance are discussed in \cite{bender}.

Recall, that property $(L_{(S,j)})$ (this notation was introduced below Proposition \ref{konig}), for $j\ge 0$ and $S\in\Se$, is equivalent to $\overline{\sigma}_j0(S)=0$. Analogously, we will say that property $(K_{(S,j)})$ holds at node $(S,j)$  if $\overline{I}_j f^{+}(S)= I_j f(S)+ \overline{I}_j f^{-}(S)$  for any $f \in \Ee_{(S,j)}$.

In this section, after remarking that for node  $(S,j)$ $(K_{(S,j)}) \implies (L_{(S,j)})$, we  provide sufficient conditions implying $(K_j)$ and $(L_j)$ for all $j \geq 0$.

\begin{proposition} \label{KjImpliesLj}
For $j\ge 0$ and $S\in\Se$ fixed, $(K_{(S,j)}) \implies (L_{(S,j)})$.
\end{proposition}
\begin{proof}
Consider $f \in \mathcal{E}_{(S,j)}$, recall that $(K_{(S,j)})$ is the equality:
$\overline{I}_j f^+(S)= I_j f(S)+\overline{I}_j f^-(S)$. Therefore,
\begin{equation} \label{firstUseOfKj}
I_j f(S) \leq I_j f(S)+ \overline{I}_j f^-(S)= \overline{I}_j f^+(S).
\end{equation}

Let $f_m \in \mathcal{D}_{(S,j)}^+, ~~m \geq 1$ satisfying $f \leq \sum\limits_{m \geq 1} f_m$, on $\Se_{(S,j)}$, then \\ $f^+ \leq \sum\limits_{m \geq 1} f_m$ on $\Se_{(S,j)}$. It follows from (\ref{firstUseOfKj})  that $I_j f(S) \leq \overline{I}_j f^+(S) \leq \sum_{m \geq 1} I_j f_m(S)$. Therefore, $(L_{(S,j)})$ holds.
\end{proof}

The proofs of  $(L_j)$ and $(K_j)$ rely on the following lemma
which is a reformulation of \cite[Lemma 3]{ferrando}
for the conditional setting. Eventually, for the actual proof of $(L_j)$
we will need to assume specific properties for nodes in $\mathcal{S}$ but, for preliminary results we do not need to assume
any property for $\mathcal{S}$. We notice that the proof
of $(L_j)$  simplifies  if all nodes are assumed to be only no-arbitrage nodes.

When we use the word {\it convergent} we mean that the corresponding limit exists in $(-\infty, \infty)$.

\begin{lemma}\label{convergenceOfPortfolioCoordinates}
 For any $m\ge 0$, let $H^m=\{H_i^m\}_{i\ge j}$, be sequences of non-anticipative functions on $\Se$, and $V^m$ functions defined on $\Se$, depending for each $S$ only on $S_0,...,S_j$, $j \ge0$ fixed. We fix a node $(S,j)$ and  assume:
\[
\Pi_{j,n}^{V^m, H^m}(\tilde{S})=V^m(S)+\sum\limits_{i=j}^{n_m-1}H_i^m(\tilde{S})\Delta_i\tilde{S}\ge 0,~~\tilde{S}\in\Se_{(S,j)},\; n_m\ge j~~\mbox{and}
\]
$~~\sum\limits_{m\ge 1}V^m(S) <\infty$.
Define, for any $\tilde{S} \in \mathcal{S}_{(S,j)}$ and $k \geq j$:
$$
H_k(\tilde{S}) \equiv
\begin{cases}
\sum_{m\ge 1} H_k^m(\tilde{S})~~~~\mbox{whenever this series is convergent}\\
*,~~\mbox{an arbitrary real number otherwise}.
\end{cases}
$$
Then, for any $\hat{S} \in \mathcal{S}_{(S,j)}$  satisfying
\begin{equation} \label{conditionOnTrajectory}
[\mbox{if $(\hat{S},p)$ is an arbitrage node, $j \leq p$, then it is of type I and $\hat{S}_{p+1} = \hat{S}_p$}],
\end{equation}
the following holds for all $k$ s.t. $j \leq k$:
\begin{equation} \label{fromManyPortfoliosToOnePortfolio}
\sum_{m=1}^{\infty} [H_k^m(\hat{S})\Delta_k\hat{S}] = H_k(\hat{S})\Delta_k\hat{S}
\end{equation}
and
\begin{equation} \label{convergentCase}
\sum_{m\ge 1} H_k^m(\hat{S})~\mbox{converges whenever $(\hat{S}, k)$ is an up-down node}.
\end{equation}
\end{lemma}

\begin{proof} The defined functions $H_k$ are non-anticipative given that the functions $H^m_k$ are non-anticipative. Fix $\hat{S}$ satisfying (\ref{conditionOnTrajectory}), 
from this hypothesis $(\hat{S}, k)$ $j \leq k$, is an arbitrage node of type I, a flat node or an up-down one. 
Whenever it is flat or an arbitrage node of type I, (\ref{fromManyPortfoliosToOnePortfolio}) holds because both sides of the equality equal $0$. Therefore, to conclude our proof it is enough to establish (\ref{convergentCase}) in case $(\hat{S}, k)$ is up-down, since it also implies (\ref{fromManyPortfoliosToOnePortfolio}).

Consider $\hat{S} \in \mathcal{S}_{(S,j)}$
satisfying (\ref{conditionOnTrajectory}) and let $j \leq k^{\ast}$ be the smallest integer such that $(\hat{S}, k^{\ast})$ is an up-down node and $\sum_{m\ge 1} H_{k^{\ast}}^m(\hat{S})$ does not converge, if such $k^{\ast}$ does not exist the lemma holds.  Therefore, by the non-convergence, there exists $\epsilon >0$, with the property that for any $M\in \mathbb{N }$ there exist $m_2>m_1\ge M$ such that
\begin{equation}\label{nonConvegentSerie}
|\sum_{m=m_1+1}^{m_2} H_{k^{\ast}}^m(\hat{S})| \ge \epsilon.
\end{equation}
Since the node $(\hat{S}, k^{\ast})$ is up-down, let
\[\theta^-=\frac 12\inf_{\tilde{S} \in \Se(\hat{S}, k^{\ast})}(\tilde{S}_{k^{\ast}+1}- \hat{S}_{k^{\ast}}) <0\;\;\mbox{and}\;\;\theta^+=\frac 12\sup_{\tilde{S} \in \Se(\hat{S}, k^{\ast})}(\tilde{S}_{k^{\ast}+1}-\hat{S}_{k^{\ast}}) >0.\]
Set $\epsilon^* \equiv \epsilon \min\{-\theta^-,\theta^+\}$.

By the choice of $(\hat{S},k^{\ast})$, it follows that if $j\le p< k^{\ast}$ then $(\hat{S},p)$ is a flat node or an arbitrage node of type I, so $\Delta_p \hat{S}=0$ in either case, or $(\hat{S}, p)$ is an up-down node and $\sum\limits_{m\ge 1} H_p^m(\hat{S})$ is convergent. In any case, there exists $M_0$ such that for any $j \le p < k^{\ast}$, and $m''> m'\ge M_0$
\begin{equation}\label{i_menorque j}
|\sum_{m=m'+1}^{m''} H_p^m(\hat{S})||\Delta_p \hat{S}| < \frac{\epsilon^*}{2^{p+2}}.
\end{equation}
Having in mind the convergence of $\sum\limits_{m\ge 1}V_m(\hat{S})=\sum\limits_{m\ge 1}V_m(S)$ there also exists $M_1\ge M_0$ such that for any  $m''> m'\ge M_1$
\begin{equation}\label{V}
\sum_{m=m'+1}^{m''}V_m(\hat{S})  < \frac{\epsilon^*}{2^{k^{\ast}+2}}.
\end{equation}
Given that $(\hat{S}, k^{\ast})$ is an  up-down node, for $M=M_1$ and the corresponding $m_1<m_2$ as in (\ref{nonConvegentSerie}), there exists $S^{1}\in\Se_{(\hat{S}, k^{\ast})}$ with $\Delta_{k^{\ast}} S^{1}\le \theta^-$ or $\theta^+\le \Delta_{k^{\ast}} S^{1}$,  such that
\begin{equation}\label{j}
\sum_{m=m_1+1}^{m_2} H_{k^{\ast}}^m(S^{1})\Delta_{k^{\ast}} S^{1} \le -\epsilon|\Delta_{k^{\ast}} S^{1}| \le -\epsilon^*.
\end{equation}

Consequently from (\ref{i_menorque j}), for $j \le p < k^{\ast}$,
\[
|\sum_{m=m_1+1}^{m_2} H_p^m(S^{1})\Delta_p S^{1}| = |\sum_{m=m_1+1}^{m_2} H_p^m(\hat{S}) \Delta_p \hat{S}| <  \frac{\epsilon^*}{2^{p+2}}.
\]
Thus, since $V_m(S^{1})=V_m(\hat{S})$, from the last equation, (\ref{V}) and (\ref{j}) it results the following contradiction:
\[
0\le \sum_{m=m_1+1}^{m_2}\Pi_{j, k^{\ast}+1}^{V_m,H^m}(S^{1})=
\]
\[
=\sum_{m=m_1+1}^{m_2}V_m(S^{1}) + \sum_{m=m_1+1}^{m_2}~~\sum_{p=j}^{k^{\ast}} H_p^m(S^{1}) \Delta_p S^{1} < -\epsilon^*(1 - \sum_{p=j}^{k^{\ast}}\frac1{2^{p+2}}) \leq \frac{- \epsilon^{\ast}}{2} < 0.
\]
We have thus established the convergence  of $\sum_{m\ge 1} H_k^m(\hat{S})$ whenever~~~
$(\hat{S}, k)$ is an up-down node.
This then completes the proof of the lemma.
\end{proof}
\subsection{Establishing Property $(L_j)$}\label{proofOfLj}
Let $S \in \mathcal{S}$; fix a node $(S,j)$ and $k \geq j$; define:
\begin{equation} \label{typicalNullSet}
\overline{N}_k(S,j) \equiv \{\overline{S} \in \overline{\mathcal{S}_{(S,j)}}:  (\overline{S}, k)~~\mbox{is an arbitrage node and}~\overline{S}_{k+1} \neq \overline{S}_k\},~\mbox{also let}
\end{equation}
\begin{equation} \nonumber
 ~\overline{N}(S,j) \equiv \cup_{k \geq j} \overline{N}_k(S,j).
\end{equation}


The next lemma is key to our arguments.

\begin{lemma} \label{twoTermsUpperBound}
Fix a node $(S,j)$ and assume $\mathcal{S}_{(S,j)}$ satisfies the RFP. \\Let  $f(\tilde{S}) \leq \sum_{m \geq 1} f_m(\tilde{S}), ~\forall~\tilde{S} \in \mathcal{S}_{(S,j)}$ where $f \in \mathcal{E}_{(S, j)} $ and $f_m \in \mathcal{D}^+_{(S, j)}$, $m \geq 1$ (we will use the notation $f_0\equiv f$ and
$f_m(\tilde{S})= V^m(S)+ \sum_{k=j}^{n_m-1} H_k^m(\tilde{S}) \Delta_k\tilde{S}$ for $m \geq 0$.)

Then:
\begin{equation} \label{twoTerms:PortfoliASndInfinity}
V^0(S) \leq \sum_{m \geq 1} V^m(S)+ \liminf_{p \rightarrow \infty} \sum_{k=j}^p
G_k(\overline{S}) \Delta_k \overline{S} + \infty ~{\bf 1}_{\overline{N}(S,j)}(\overline{S}),~~~~\forall~\overline{S} \in \overline{\mathcal{S}_{(S,j)}},
\end{equation}
where,  for any $\hat{S} \in \mathcal{S}_{(S, j)}$, $G_k(\hat{S}) \equiv H_k(\hat{S})- H^0_k(\hat{S}), j \leq k \leq n_0 -1$,
$G_k(\hat{S}) \equiv H_k(\hat{S}), k \geq n_0$ and
\begin{equation}  \label{definedWhenConvergentA}
H_k(\hat{S}) \equiv
\begin{cases}
\sum_{m\ge 1} H_k^m(\hat{S})~~~~\mbox{whenever the series is convergent,}\\
0,~~\mbox{otherwise}.
\end{cases}
\end{equation}
\end{lemma}
\begin{proof}
Consider $\overline{S} \in \overline{\mathcal{S}_{(S,j)}}$ with
$\overline{S} = \lim_{n \rightarrow \infty} S^n$ and $S^n$ as in (\ref{rFP}), then
\begin{equation} \label{firstApplyRFP}
f(\overline{S})= \limsup_{n \rightarrow \infty} f(S^n) \leq \limsup_{n \rightarrow \infty} \sum_{m \geq 1} f_m(S^n) \leq  \sum_{m \geq 1} f_m(\overline{S}).
\end{equation}

The defined functions $H_k$ are non-anticipative given that the functions $H_k^m$ are non-anticipative and, therefore, so are the functions $G_k$ as well. Notice that our hypothesis on $\mathcal{H}_{(S,j)}$ allows us
to assume $H_i^m(\tilde{S})=0$ for all
$i \geq n_m$ and for all $\tilde{S} \in \mathcal{S}_{(S,j)}$ (this can be achieved despite
that $H^0$, being the difference of two elements of  $\mathcal{H}_j$, does not necessarily belong to $\mathcal{H}_j$). These comments
also apply to the case when the portfolios are extended to act
on $\overline{\mathcal{S}_{(S,j)}}$ as we do below.
From the standard Fatou's lemma
and (\ref{firstApplyRFP}) we obtain:
\begin{equation} \label{fromFatou}
V^0(S) \leq \sum_{m \geq 1} V^m(S)- \sum_{k=j}^{n_0-1} H_k^0(\overline{S})~ \Delta_k \overline{S} +  \liminf_{p \rightarrow \infty} \sum_{m \geq 1} \sum_{k=j}^p
[H^m_k(\overline{S}) ~\Delta_k \overline{S}],~~\forall~\overline{S} \in \overline{\mathcal{S}_{(S,j)}}.
\end{equation}
Whenever $\overline{S} \in \overline{N}(S,j)$ the r.h.s. of (\ref{twoTerms:PortfoliASndInfinity}) equals infinity and so upperbounds the r.h.s. of (\ref{fromFatou}).
Let then $\overline{S} \in  \overline{\mathcal{S}_{(S,j)}} \setminus \overline{N}(S,j)$; it then follows that
for any $k \geq j$ such that $(\overline{S}, k)$ is an arbitrage node it can only be so if it is of type I and $\overline{S}_{k+1} = \overline{S}_k$ and so condition (\ref{conditionOnTrajectory}) in Lemma \ref{convergenceOfPortfolioCoordinates} holds (notice that
we are applying that lemma to $\overline{\mathcal{S}_{(S,j)}}$).  We can then conclude from that lemma that the only possibility for $\sum_{m \geq 1}H^m_k(\overline{S})$ to diverge is that $(\overline{S}, k)$ is flat or an arbitrage node of type I but in any of those cases $\sum_{m \geq 1}[H^m_k(\overline{S})~ \Delta_k \overline{S}] =0$. Therefore  $\sum_{m \geq 1}
[H^m_k(\overline{S}) ~\Delta_k \overline{S}] = H_k(\overline{S}) ~\Delta_k \overline{S},~~\forall k \geq j$, where the $H_k$ are given by (\ref{definedWhenConvergentA}).
It then follows that, whenever $\overline{S} \in  \overline{\mathcal{S}_{(S,j)}} \setminus \overline{N}(S,j)$, (\ref{twoTerms:PortfoliASndInfinity})
equals (\ref{fromFatou}).
\end{proof}

Next we define contrarian trajectories; they were originally introduced in \cite[Def. 12]{ferrando} and here we modify their use somehow.

\begin{definition}[Contrarian Trajectories, CT] \label{def:contrarianTrajectory}
We will say that a trajectory set $\Se$ has the {\it contrarian trajectory} (CT) property  if  the following holds: for any $S \in \Se, \; n_0 \ge 0$,  any  $F=\{F_i\}_{i\ge n_0}$, a sequence of non anticipative functions on $\mathcal{S}_{(S, n_0)}$, and $\epsilon >0$,  there exists $\overline{S}^{\epsilon} \in \overline{\Se_{(S,n_0)}}$ satisfying
\begin{equation} \nonumber
\liminf_{n \rightarrow \infty} \sum_{i=n_0}^{n} F_{i}(\overline{S}^{\epsilon}) ~~\Delta_i \overline{S}^{\epsilon} \leq \epsilon.
\end{equation}
Such $\overline{S}^{\epsilon}$ will be called a contrarian trajectory (CT) ({\it to $F$ at $n_0$} and beyond).
\end{definition}
\begin{remark}
If $\Se$ is locally $0$-neutral, then $\mathcal{S}$ has the CT property. This follows from (\ref{sequenceInequality}).
\end{remark}

Lemma \ref{CTAvoidsN} below provides more detailed information on contrarian trajectories under an assumption on nodes.


\begin{lemma} \label{CTAvoidsN}
Assume  nodes in $\mathcal{S}$ are only of two types: no-arbitrage nodes or arbitrage nodes of type I. Let $(S,n_0)$ be any node and  $\epsilon >0$ and $F = \{F_i\}_{i \geq n_0}$ be non-anticipative functions. Then, there exists a CT\hspace{.05in} (to $F$ at $n_0$ and beyond) ~~~$\overline{S}^{\epsilon} \in \overline{\mathcal{S}_{(S, n_0)}}$ such that
 $\overline{S}^{\epsilon} \notin \overline{N}(S,n_0)$ (the set $\overline{N}(S, n_0)$ was introduced in (\ref{typicalNullSet})).
\end{lemma}



\begin{proof}
From our assumption on nodes, it follows that each node $(S,j)$ is $0$-neutral. Therefore, Lemma \ref{scapingSequence} applies
and the sequence $S^n$ in (\ref{sequenceInequality}) gives the desired
$\overline{S}^{\epsilon} \equiv \lim_{n \rightarrow \infty} S^n$.
The property $\overline{S}^{\epsilon} \notin \overline{N}(S,n_0)$
follows from Remark \ref{contrarianAtArbitrageNodeOfTypeI}.
\end{proof}

\vspace{.1in}\noindent
{\bf Proof of Theorem \ref{proofOfL}}:
\vspace{-.055in}
\begin{proof}
We do the analysis for a fixed, but arbitrary, node $(S,j)$;
consider then (\ref{twoTerms:PortfoliASndInfinity}) in Lemma \ref{twoTermsUpperBound},
for a given $\epsilon >0$ and the non-anticipative functions $G_k$ introduced in that lemma, it follows from Lemma \ref{CTAvoidsN} that there exists
$\overline{S}^{\epsilon} \in \overline{\mathcal{S}_{(S,j)}}$\\
 such that $\liminf_{p \rightarrow \infty} \sum_{k=j}^p
G_k(\overline{S}^{\epsilon})~ \Delta_k \overline{S}^{\epsilon} \leq \epsilon$ and ${\bf 1}_{\overline{N}(S,j)}(\overline{S}^{\epsilon})=0$. Therefore $V_0(S) \leq \sum_{m \geq 1} V_m(S) + \epsilon$ which concludes the proof.
\end{proof}

It is important that we managed to prove property $(L_j)$ under hypothesis that allow for arbitrage nodes as they provide
non-trivial examples of null sets as indicated in the next proposition.
As the portfolio which exploits the arbitrage may not be included
in minimal portfolio sets satisfying (H.1)--(H.4), we state the proposition, for sake of simplicity, under the assumption that there
are no trading restrictions.

\begin{proposition}
 Suppose that each portfolio set $\He_{(S',j)}$, $S'\in \mathcal{S}$, $j\geq 0$, consists of all nonanticipating sequences of functions.
 Given an arbitrage node $(S,k)$ in $\mathcal{S}$,  define: $A \equiv \mathcal{S}_{(S,k)} \setminus \{\tilde{S} \in \mathcal{S}_{(S,k)}: \tilde{S}_{k+1} = \tilde{S}_k\}$. Then $\overline{I}_j({\bf 1}_{A})=0$ for all $0 \leq j \leq k$.
\end{proposition}
\begin{proof}
Without loss of generality assume $\tilde{S}_{k+1} \geq \tilde{S}_k$ for all $\tilde{S} \in \mathcal{S}_{(S,k)}$. Let, for $m \geq 1$, $H_k(\hat{S}) =1$ if $\hat{S}_i = S_i$ $ 0 \leq i \leq k$ and $H^m_k(\hat{S}) =0$  otherwise. Set also $H_n^m(\hat{S})=0$ for all $n \neq k$ for all $\hat{S} \in \mathcal{S}$ and all $m \geq 1$. Then
\begin{equation} \label{boundForNullSet}
{\bf 1}_A(\tilde{S}) \leq \sum_{m \geq 1} H_k^m(\tilde{S}) \Delta_k \tilde{S}
\end{equation}
holds for all $\tilde{S} \in \mathcal{S}_{(S,j)}$. Notice that
$(S,k)$ being an arbitrage node, the right hand side of (\ref{boundForNullSet}) equals $\infty$ whenever the left hand side equals $1$.
\end{proof}


The following proposition clarifies how one can check the validity
of the hypothesis in Theorem \ref{proofOfL}.

\begin{proposition} \label{validityOfHypothesisInLjTheorem}
If $\mathcal{S}$ is complete, then it satisfies GRFP. If, moreover, $\mathcal{S}$ satisfies GRFP, then
RFP holds at any node $(S,j)$ and for any choice of the portfolio sets saytisfying (H.1)--(H.4).
\end{proposition}

The proof relies on the following lemma.

\begin{lemma} \label{rFPGlobalImpliesRFPConditional}
Assume $\mathcal{S}$ satisfies the RFP and $\mathcal{H}_0$ satisfies
the following property:

\begin{itemize}
\item If $k \geq 0$, $\tilde{S} \in \mathcal{S}$
and $(H_i)_{i \geq k} \in \mathcal{H}_{(\tilde{S},k)}$
then $(H_i~{\bf 1}_{i \geq k}~{\bf 1}_{\mathcal{S}_{(\tilde{S}, k)}})_{i \geq 0} \in \mathcal{H}_{0}$ (i.e. it is admissible to do nothing until time $k$ and, if
a trajectory in node $(\tilde{S}, k)$ realizes, apply
strategy $(H_i)_{i \geq k}$, otherwise do nothing).
\end{itemize}
 Then,  for any node $(S,j)$ with $j \geq 0$,  $\mathcal{S}_{(S,j)}$
also satisfies the RFP.
\end{lemma}
\begin{proof}
Fix an arbitrary node $(S, j)$, $j \geq 0$ and take $\overline{S}
\in \overline{\mathcal{S}_{(S,j)}}$; therefore, 
$\overline{S}\in \overline{\mathcal{S}}$ as well. Consider $g_m(\hat{S}) = U^m(S)+ \sum_{i=j}^{n_m-1} G_i^m(\hat{S}) ~\Delta_i \hat{S}$, $m \geq 0$, $\hat{S} \in \mathcal{S}_{(S,j)}$, such that
$g_0 \in \mathcal{E}_{(S, j)}, g_m \in \mathcal{D}^+_{(S, j)}$, $m \geq 1$, and $\sum_{m \geq 1} U^m(S) < \infty$.

Define, for all $m \geq 0$ and $\tilde{S} \in \mathcal{S}$ :  $(i)$  $H^m_i(\tilde{S}) =0$ whenever $0 \leq i < j$,
 $(ii)$
 let $H_i^m=G_i^m$ on $\mathcal{S}_{(S,j)}$
  and $H^m_i=0$ otherwise, whenever $j \leq i$. Also set the constants $V^m \equiv U^m(S)$ for all $m \geq 0$. Define
  $f_m(\tilde{S}) \equiv V^m + \sum_{i=0}^{n_m-1} H_i^m(\tilde{S})~\Delta_i \tilde{S}$ for $\tilde{S} \in \mathcal{S}$ and notice
  $f_m \in \mathcal{D}^+$, $m \geq 1$, $f_0 \in \mathcal{E}$
  and $\sum_{m \geq 1} V^m< \infty$.

Given that $\mathcal{S}$ satisfies the RFP and $\overline{S} \in \overline{\mathcal{S}}$, for any $n\ge 0$ there exists $S^n \in \mathcal{S}_{(S,j)}$ such that
$\overline{S} = \lim_{n \rightarrow \infty} S^n $, satisfying
\begin{equation} \label{rFPWithoutInitialValue}
\sum_{m \geq 0} \limsup_{n \rightarrow \infty} \sum_{i=0}^{n_m-1} H^m_i(S^n)~\Delta_i S^n \geq
\limsup_{n \rightarrow \infty} \sum_{m \geq 0}  \sum_{i=0}^{n_m-1} H^m_i(S^n)~\Delta_i S^n.
\end{equation}

Define now $\tilde{S}^n \equiv S^{j+n}$ for $n \geq 0$, notice that since $S^{j+n}\in \mathcal{S}_{(S,j)}$ then $\tilde{S}^n \in \mathcal{S}_{(S,j)}$, and $\overline{S}= \lim_{n \rightarrow \infty} \tilde{S}^n $, this follows because
$\overline{S}_n = S^n_n= S^{n+j}_n= \tilde{S}^n_n$.

Using (\ref{rFPWithoutInitialValue}), we compute as follows
\begin{eqnarray*} \nonumber
&&\sum_{m \geq 0} \limsup_{n \rightarrow \infty} \sum_{i=j}^{n_m-1} G^m_i(\tilde{S}^n)~\Delta_i \tilde{S}^n =
\sum_{m \geq 0} \limsup_{n \rightarrow \infty} \sum_{i=j}^{n_m-1} G^m_i(S^{n+j})~\Delta_i S^{j+n}
\\
&=&\sum_{m \geq 0} \limsup_{n \rightarrow \infty} \sum_{i=0}^{n_m-1} H^m_i(S^{n+j})~\Delta_i S^{n+j}= \sum_{m \geq 0} \limsup_{n \rightarrow \infty} \sum_{i=0}^{n_m-1} H^m_i(S^n)~\Delta_i S^n
\\
&\geq& \limsup_{n \rightarrow \infty} \sum_{m \geq 0}  \sum_{i=0}^{n_m-1} H^m_i(S^n)~\Delta_i S^n=
\limsup_{n \rightarrow \infty} \sum_{m \geq 0}  \sum_{i=0}^{n_m-1} H^m_i(S^{j+n})~\Delta_i S^{j+n}
\\
&=&\limsup_{n \rightarrow \infty} \sum_{m \geq 0}  \sum_{i=j}^{n_m-1} G^m_i(\tilde{S}^n)~\Delta_i \tilde{S}^n,
\end{eqnarray*}
we have then completed the proof.
\end{proof}

\begin{proof}[Proof of Proposition \ref{validityOfHypothesisInLjTheorem}]
 First note that GRFP is obviously equivalent to the RFP
 at the initial nodes $\mathcal{S}=\mathcal{S}_{(S,0)}$, if there
 are no trading restrictions, i.e. if each of the portfolio sets
  $\He_{(S',i)}$, $S'\in \mathcal{S}$, $i\geq 0$, consists of all nonanticipating sequences of functions. Then, on the one hand,
  by Remark \ref{completenessImpliesRFP}, completeness of $\mathcal{S}$ implies GRFP.
  On the other hand, Proposition \ref{validityOfHypothesisInLjTheorem}
  implies RFP at every node $\mathcal{S}_{(S,j)}$, when there are no trading restrictions. The RFP clearly is still satisfied, if one reduces the portfolio set from the set of all nonanticipating sequences of functions to any subset  $\He_{(S,j)}$ satisfying (H.1)--(H.4).
\end{proof}

Below we present two examples, in the first one $(S,0)$ is an arbitrage node of type I but
we have failure of $(L_{(S,0)})$, this is explained by the failure of RFP.

\subsubsection{Example 1}
Let $\Se$ the trajectory space composed by the following trajectories:
\[S^{n},\; n\in\mathbb{N}:\quad S^{n}_i=1,\;\;0\le i<n;\quad \quad S^{n}_i=2\;\;i\ge n.\]

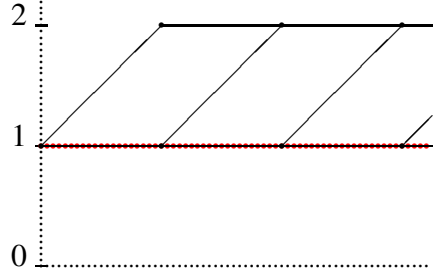
\[
{\setlength{\unitlength}{.08 cm} \label{grafica}
\begin{picture}(80,50)
\multiput(0,5)(0,1){45}{\circle*{.15}}\multiput(0,5)(1,0){65}{\circle*{.15}}
\put(-5,5){0}\put(-1,5){\line(1,0){2}}\put(-5,25){1}\put(-1,25){\line(1,0){2}}\put(-5,45){2}\put(-1,45){\line(1,0){2}}
{\color{red}\multiput(0,25)(1,0){65}{\circle*{.7}}}\put(0,25){\line(1,0){65}}
\put(20,45){\line(1,0){45}}
\put(0,25){\line(1,1){20}}\put(20,25){\line(1,1){20}}\put(40,25){\line(1,1){20}}
\put(0,25){\circle*{1}}\put(20,25){\circle*{1}}\put(40,25){\circle*{1}}
\put(20,45){\circle*{1}}\put(40,45){\circle*{1}}\put(60,45){\circle*{1}}
\put(60,25){\circle*{1}}
\put(60,25){\line(1,1){5}}
\end{picture}}
\]

Here $\overline{\sigma}_0 0= - \infty$, we will show that the reverse Fatou condition fails  (which is what should happen by Theorem \ref{proofOfL}).

We note that $\overline{\mathcal{S}}\setminus \mathcal{S}= \{\overline{S}: \overline{S}_i \equiv 1,\; i\ge 0\}$, and we will show that the RFP fails for
the portfolios: $V^m =0$, $H^m_m(S)= 1$ if $S_i=1$ for $0 \leq i \leq m$, $H_i^m(S)= 0$ otherwise, therefore $n_m= m+1$.
Then $f_m(S^n)= (S^{n}_{n+1}- S^{n}_{n})= 1$ if $m =n$ and
$f_m(S^n)= 0$ if $m \neq n$. Therefore 
$\sum_{m \geq 1} f_m(S) = 1$ for all $S\in \mathcal{S}$, whereas 
$\sum_{m \geq 1} f_m(\bar S) = 0$. Thus, for every sequence 
$\{S^{(k)}\}_{k\in \mathbb{N}}$ in $\mathcal{S}$ satisfying $\lim_{k\rightarrow \infty} S^{(k)} =\bar S$, we obtain
$$
\limsup_{k\rightarrow \infty} \sum_{m \geq 1} f_m(S^{(k)}) = 1>0=
\sum_{m \geq 1} f_m(\bar S) =\sum_{m \geq 1} \limsup_{k\rightarrow \infty} f_m(S^{(k)}),
$$
i.e. RFP fails.

\vspace{.1in}
In the following example all nodes are no-arbitrage nodes or arbitrage nodes of type I and we show that the RFP holds, it then follows from Theorem \ref{proofOfL} that $(L_{(S,j)})$ holds
at all nodes. The point of this example is that $\mathcal{S}$ is not complete.

\subsubsection{Example 2} Consider the example with $\Se$ given in the following figure,
here we know by direct computation that $\overline{\sigma}_{j, \mathcal{S}} 0= 0$ at all nodes $(S,j)$. We will verify this result by an application of
Theorem \ref{proofOfL} by  checking that the reversed Fatou condition  holds for this example. The trajectories are
\[S^{u,n},\; n\in\mathbb{N}:\quad S^{u,n}_i=2,\;\;0\le i<n;\quad \quad S^{u,n}_i=3\;\;i\ge n,\;\;\mbox{and}\]
\[S^{d,n},\; n\in\mathbb{N}:\quad S^{d, n}_i=2,\;\;0\le i<n;\quad \quad S^{d, n}_i=1\;\;i\ge n.\]
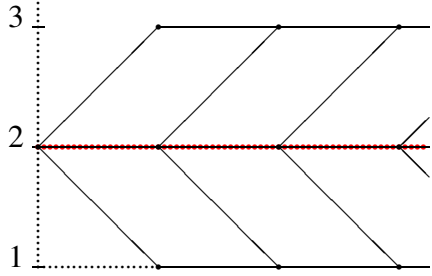
\[
{\setlength{\unitlength}{.08 cm} \label{grafica}
\begin{picture}(80,50)
\multiput(0,5)(0,1){45}{\circle*{.15}}\multiput(0,5)(1,0){20}{\circle*{.15}}
\put(-5,5){1}\put(-1,5){\line(1,0){2}}\put(-5,25){2}\put(-1,25){\line(1,0){2}}
\put(-5,45){3}\put(-1,45){\line(1,0){2}}
{\color{red}\multiput(0,25)(1,0){65}{\circle*{.7}}}\put(0,25){\line(1,0){65}}
\put(20,5){\line(1,0){45}}\put(20,45){\line(1,0){45}}

\put(0,25){\line(1,1){20}}\put(0,25){\line(1,-1){20}}

\put(0,25){\circle*{1}}
\put(20,25){\circle*{1}}\put(40,25){\circle*{1}}\put(60,25){\circle*{1}}
\put(20,25){\line(1,1){20}}\put(20,25){\line(1,-1){20}}
\put(20,5){\circle*{1}}\put(40,5){\circle*{1}}\put(60,5){\circle*{1}}
\put(20,45){\circle*{1}}\put(40,45){\circle*{1}}\put(60,45){\circle*{1}}
\put(40,25){\line(1,1){20}}\put(40,25){\line(1,-1){20}}
\put(60,25){\line(1,1){5}}\put(60,25){\line(1,-1){5}}

\end{picture}}
\]
We note that $\overline{\mathcal{S}}\setminus \mathcal{S}= \{\overline{S}: \overline{S}_i=2,\;i\ge 0\}$, we will consider
the sequence \\ $S^n \equiv (2,2, \ldots, 2,\ast, \ast, \ast, \ldots)$, where the last appearing $2$ is the nth coordinate of the sequence and the entry value $\ast$ will be equal to $3$ or $1$ (i.e. $S^n_i =3$ for all $i \geq n+1$ or  $S^n_i =1$ for all $i \geq n+1$), which value to be chosen, $3$ or $1$, will depend on $n$ as well as on a given sequence $f_m \in \mathcal{E}^+,\;\;m\ge 1$ and $f_0 \in \mathcal{E}$. The sequence $S^n$ to be constructed is a variation from
the usual CT construction, now, a contrarian move is used to define  $S^n$ but discarded in $S^{n+1}$ given that in fact
we keep the flat trajectory (which, of course, also acts as a CT)
converging to $\overline{S}$. The construction of $S^n$ is completed  next.

Given $f_0 \in \mathcal{E},\; f_m \in \mathcal{D}^+\;m\ge 1$  with $\sum_{m \geq 1} V^m < \infty$;
notice that
$H^m_i(S^n) (S^n_{i+1}- S^n_{i}) = 0$ if $i \neq n$, so $\sum_{m \geq 1} H^m_i(S^n) (S^n_{i+1}- S^n_{i}) = 0$ if $i \neq n$ and $f_m(S^n)=V^m$ if $n \ge n_m$.

Set $H_i(S) \equiv H^0_i(S)+\sum_{m \geq 1} H^m_i(S)$ which exists
by an application lemma \ref{convergenceOfPortfolioCoordinates} (in fact we only need
to apply the lemma along $\overline{S}$, i.e. we only need $H_i(\overline{S})=\sum_{m \geq 0} H^m_i(\overline{S})$).
Depending if $H_n(\overline{S}) \leq 0$ choose $S^n_i = \ast= 1$ for $i \geq n+1$ or if $H_n(\overline{S}) > 0$ choose $S^n_i = \ast= 3$ for $i \geq n+1$. (e.g. if $H_3(\overline{S}) \leq 0$, $S^3\equiv S^{d,4}=(2,2,2,2, 1,1,1, \ldots)$,  if $H_4(\overline{S}) > 0$, $S^4\equiv S^{u,5}=(2,2,2,2,2, 3,3,3, \ldots)$ etc).
Therefore, using the standard Fatou's lemma,
\begin{equation} \label{integralBound}
\sum_{m \geq 0} f_m(S^n) = \sum_{m \geq 0}[V^m + \sum_{i=0}^{n_m-1}
H_i^m(S^n) \Delta_i S^n] \leq \liminf_{p \rightarrow \infty}[
\sum_{m \geq 0}V^m + \sum_{i=0}^{p} \sum_{m \geq 0}
H_i^m(S^n) \Delta_i S^n] =
\end{equation}
\begin{equation} \nonumber
\liminf_{p \rightarrow \infty}[ \sum_{m \geq 0}V^m+  {\bf 1}_{p \geq n} \sum_{m \geq 0}  H_n^m(S^n) (S^n_{n+1} - S^n_n))]= \sum_{m \geq 0}V^m +H_n(S^n) (S^n_{n+1}- S^n_n)
\leq \sum_{m \geq 0}V^m .
\end{equation}
Inequality (\ref{integralBound}), combined with the fact that $\limsup\limits_{n \rightarrow \infty} f_m(S^n)=V^m$, 
allows us to check the reversed Fatou condition directly:
\begin{equation} \nonumber
\limsup_{n \rightarrow \infty} \sum_{m \geq 0} f_m(S^n)
\leq \sum_{m \geq 0}V^m = \sum_{m \geq 0}\limsup_{n \rightarrow \infty}  f_m(S^n).
\end{equation}
Then, Theorem \ref{proofOfL} is applicable to this example.


\subsection{Establishing property $(K_j)$} \label{establishingKj}

We note that the hypotheses in the following result are stronger
that the ones in Theorem \ref{proofOfL} needed to prove $(L_j)$.
In particular, we require the following closedness property for the accumulation of
countably many portfolios:
\begin{itemize}
 \item [(H.5)] For any sequence $\{H^m\}_{m\geq 1}$ in $\mathcal{H}_{(S,j)}$ the portfolio $H$ given by
  $$
  H_i(\tilde S):=\left\{\begin{array}{cl} \sum_{m=1}^\infty H_i^m(\tilde S ), & \textnormal{ if the series is convergent in } \mathbb{R}\\  0, & \textnormal{otherwise} \\ \end{array} \right.
  $$
  belongs to $\mathcal{H}_{(S,j)}$.
\end{itemize}

\begin{theorem}\label{proofOfK_j}
  (i) Assume $\mathcal{S}_{(S,j)}$ satisfies the RFP for each $S \in \mathcal{S}$ and each $j \geq 0$, $\mathcal{S}$ is free of arbitrage nodes, and the portfolio sets satisfy (H.1)--(H.5).  Then $(K_j)$ is valid.  \\[0.2cm]
 (ii) Assume $\mathcal{S}_{(S,j)}$ satisfies the RFP for each $S \in \mathcal{S}$ and each $j \geq 0$, $\mathcal{S}$ is free of arbitrage nodes of type II, and there are no trading restrictions, i.e. the portfolio sets consist of all sequences of nonanticipating functions. Then $(K_j)$ is valid.
\end{theorem}

The proof is a consequence of Theorem \ref{proofOfL} and the following result due to \cite{bender}.

\begin{theorem}  \label{lAndKEquivalentIfValidAtAllNodes}
Under each of the sets of assumptions (i) or (ii) in Theorem  \ref{proofOfK_j}, $(L_j)$ implies $(K_j)$.
\end{theorem}
Note that \cite{bender} only considers portfolio sets without trading restrictions, and  the case (ii) is directly covered by Theorem 5.1 in \cite{bender}. The proof of Theorem 5.1 in \cite{bender} consists of two key steps. In the first step it is shown that all involved portfolios may be assumed to have the same finite maturity. The corresponding manipulations make use of (H.3) and (H.4) only. In the second step, a variant of Lemma \ref{convergenceOfPortfolioCoordinates}
is applied to accumulate portfolios. Condition (H.5) is required
to ensure that the accumulated portfolio is still admissible. In the absence of arbitrage nodes of type I, the technical ramifications of the proof of Theorem 5.1 in \cite{bender} involving (null) sets $\mathcal{N}_n$ and a random time $\tau$ become trivial, and so the identical proof works in the case (i) as well.


\begin{thebibliography}{99}

\bibitem{acciaio} B. Acciaio, M. Beiglb\"ock, F. Penkner and  W. Schachermayer (2016),
\emph{A Model-Free Version of the Fundamental Theorem of Asset Pricing and the Super-Replication Theorem}. Mathematical Finance, {\bf 26} (2), 233-251.



\bibitem{bender} C. Bender, S.E. Ferrando and A.L. Gonzalez (2021),  \emph{Model-Free Finance and Non-Lattice Integration}; submitted for publication. {\tt arXiv:2105.10623 [q-fin.MF]}

\bibitem{biagini} S. Biagini  and R. Cont, R. (2007). \emph{Model-free representation of pricing rules as conditional expectations}.  Stochastic Processes and Applications to Mathematical Finance, Proceedings of the 6th Ritsumeikan International Symposium, World Scientific (2007), pp. 53-66

   \bibitem{burzoni} M. Burzoni, M. Frittelli, and M. Maggis (2017), \emph{Model-free superhedging duality}. Annals of Applied Probability, {\bf 27} (3), 1452-1477.



\bibitem{dalang} R.C. Dalang, A. Morton, W. Willinger, (1990), \emph{Equivalent Martingale
measures and no-arbitrage in stochastic securities market model}. Stochastics and Stochastic Reports, {\bf 29}, 185-201.




\bibitem{ferrando} S.E. Ferrando and A.L. Gonzalez (2018), \emph{Trajectorial martingale transforms.
Convergence and integration}. New York Journal of Mathematics, {\bf 24}, 702-738.


\bibitem{ferrando2} S.E.~Ferrando, A. Fleck, A. Gonzalez and A. Rubtsov (2019),
\emph{Trajectorial Asset Models with Operational Assumptions}.
Quantitative Finance and Economics, {\bf 3} (4), 661-708.


\bibitem{kassberger} S. Kassberger and T. Liebmann (2016),  \emph{An Alternatice Axiomatic Characterisation of pricing Operators}. Journal of Applied Probability, {\bf 53}, 1257-1264.


\bibitem{konig} H. K\"{o}nig (1982). \emph{Integraltheorie ohne Verbandspostulat}. Mathematische Annalen {\bf 258}, 447-458.


\bibitem
{leinert} M. Leinert (1982),  \emph{Daniell-Stone integration without the lattice condition}. Archiv der Mathematik, {\bf 38}, 258-265.

\bibitem{schachermayer} W. Schachermayer, (1994), \emph{Martingale Measures for Discrete time Processes with Infinite Horizon}. Mathematical Finance {\bf 4}, 25-56.



\end{thebibliography}
\end{document}